\documentclass{article}
\usepackage{amsmath}
\usepackage{amsthm}
\usepackage{amssymb}
\usepackage{mathrsfs}
\usepackage{graphicx}
\usepackage{verbatim}
\usepackage{tikz}
\usepackage[all]{xy}
\usepackage{graphics}
\usepackage[colorlinks=true]{hyperref}
\hypersetup{urlcolor=blue, citecolor=red}

\newtheorem{theorem}{Theorem}[section]
\newtheorem{lemma}[theorem]{Lemma}
\newtheorem{corollary}[theorem]{Corollary}
\newtheorem{proposition}[theorem]{Proposition}
\theoremstyle{definition}

\newtheorem{defin}{Definition}[section]

\newtheorem{remark}{Remark}

\newtheorem{Oquestion}{Open question}

\newcommand{\R}{\mathbb{R}}

\newcommand{\N}{\mathbb{N}}

\newcommand{\esse}{\mathbb S}

\def\beq{\begin{equation}}
\def\eeq{\end{equation}}

\def\eps{\varepsilon}

\begin{document}
{\large \textit{The mathematical content that follows dates from 9 March 2016. I began writing it in early December 2014, whilst I was visiting Scuola Normale Superiore in Pisa, Italy. I would like to thank Luigi Ambrosio, Andrea Malchiodi, Carlo Mantegazza e Gérard Tronel who, between January 2015 and March 2016, read this version of the manuscript and some earlier drafts. Their feedback and constant support have been a source of encouragement.}

\textit{The manuscript has been at rest since March 2016. Neither its content nor its status will undergo changes in the future; it asks merely to be considered through the eyes of that moment.}

\vskip1truecm

Paris, October 23, 2026 \hfill Alberto Farina\footnote{In the forthcoming Note titled  \textit{The Eleven Questions : March 2016 - October 2026} we will review the current status of the eleven open questions posed in Appendix H of the present manuscript.}}

\vspace{2em}

\title{\sc On the classification of entire solutions to the Allen-Cahn equation in any dimension}
\date{}
\author{Alberto Farina }
\numberwithin{equation}{section}
\maketitle
{\footnotesize
\centerline{LAMFA, CNRS UMR 7352} 
\centerline {Universit\'e de Picardie Jules Verne}
\centerline{33 rue Saint-Leu, 80039 Amiens, France}}

\begin{abstract} 
We study and classify entire solutions of the Allen-Cahn equation in any space dimension.  After introducing the notion of energy, and the one of energy levels, we give complete or partial answers to open questions about minimizers, quasiminimizers, stable solutions or solutions which are stable outside a compact set, as well as monotone solutions. We also prove a quantization effect in the two-dimensional case, a complete characterization of one-dimensional solutions in any space dimension $N \ge 2$ (various simple and natural criteria for local minimizers and/or for monotone solutions) etc..... Those results provide a further step towards the solution of a celebrated conjecture of De Giorgi, as well as its generalization to the framework of stable/finite Morse Index solutions. In particular we prove De Giorgi's conjecture for finite energy solutions in $\R^4$.
We also provide a section devoted to some open questions. It contains further discussions, comments and remarks about the validity of the considered questions and conjectures.

\end{abstract}


\tableofcontents

\section{General structure of the set of solutions and the definition of energy} \label{sec: Gene+energy}

\medskip

The aim of this work is to study and classify entire solutions to the Allen-Cahn equation, 

\beq\label{eqAC}
-\Delta u = u - u^3  \qquad \text{in} \qquad \R^N , \qquad N \ge 1.
\eeq

This is the natural continuation of the project we have started in our earlier works \cite{FarGibbons, Fa-BUMI, FarLincei, FarDisc, FarDisc-Cras, FarHab, Fa, Fa3, FMV, Far-Dino-Enrico-SNS, FVquasimin, FV, FVstate, FVtams,FVcomPDE}.

\medskip

Throughout the paper $B(x,R) = B^N(x,R)$ will denote the $N-$dimensional open ball of center $x \in \R^N$ and radius $R>0$.

\medskip

Our first result is 

\medskip

\begin{theorem} {\textbf{(Regularity and universal bounds)}} \label{Reg-univBound}

\noindent Every distribution solution $u \in L^3_{loc}(\R^N)$ of the Allen-Cahn equation \eqref{eqAC} is smooth, actually real-analytic, and satisfies 

\item i)  \quad $ \vert u \vert \le 1 $ \quad on \quad $ \R^N $.

\medskip

\noindent Moreover, 

\smallskip

\item ii) \, either \quad $ \vert u \vert <1$ \quad on \quad $\R^N$ \qquad or \qquad $ u \equiv \pm1$,

\item iii)  if $ N \ge 2$ and $u$ is not constant, then $u$ must change sign and the set $  \{ \, u = 0 \, \} $ is not empty and unbounded. 

\end{theorem}

\begin{proof} [Proof of Theorem \ref{Reg-univBound}]  

By Proposition 1.9 in \cite{FarFEQL} we know that {\it every distribution solution} $u \in L^3_{loc}(\R^N)$ of the Allen-Cahn equation \eqref{eqAC} is smooth and satisfies $ \vert u \vert \le 1$ on $ \R^N$. Therefore the real-analiticity is a consequence of the fact that the r-h-s of the Allen-Cahn equation is a polynomial in the variable $u$ (see \cite{Bern-analitico, Hopf-analitico, Mor-Nir, Mor-analitico}). The content of item ii) and the first claim of iii) are an immediate consequence of the strong maximum principle. If the set $  \{ \, u = 0 \, \} $ were bounded, we would have $ u \equiv \pm1$ by Theorem 2.1 in \cite{Fa} (see Theorem \ref{Teo-superGibbons} in Section \ref{sec : Soluzioni-controllate}), a contradiction. 
\end{proof}

In view of the above result, in the sequel we always consider classical solutions. 

\bigskip

Next we recall the fundamental contributions of  L. Modica about bounded entire solutions of the Allen-Cahn equation. The deep and elegant results proved by L. Modica in his seminal works \cite{ModicaGradBound, ModicaMonot} will play a crucial role in our analysis.

\medskip

\begin{theorem}{\textbf{(L. Modica \cite{ModicaGradBound, ModicaMonot} )}} \label{Modica} 

\noindent Let $u \in \mathcal{C}^2(\R^N)$ be a bounded solution of \eqref{eqAC}. Then we have 

\smallskip

\item i) {\sl Gradient bound : }

\beq 
\frac{|\nabla u|^2}{2}  \le \frac{(1-u^2)^2}{4} \qquad {\rm on} \qquad \R^N. 
\eeq
 
\medskip

\item ii) {\sl Monotonicity formula : }

\smallskip

\noindent the function $ \phi (R) : = \frac{1}{R^{N-1}}  \Big [ \int_{B(x,R)} \frac{\vert \nabla u\vert^2}{2} + \frac{(1-u^2)^2}{4} \Big ] $ is non-decreasing on $(0, +\infty)$.

\end{theorem}

\medskip

The above gradient bound has been refined, in a sharp way, and extended to entire solutions of a general semilinear equation of the form $ -\Delta u = f(u) $, with $ f \in C^{0,1}_{loc}(\R)$, by the present author, B. Sciunzi and E. Valdinoci in \cite{Far-Dino-Enrico-SNS}. We state here only the special case concerning the Allen-Cahn equation. The interested reader can find the general case in \cite{Far-Dino-Enrico-SNS} (cf. Lemma 4.11 therein), where also the quasilinear case is studied.    

\medskip

\begin{theorem}\label{raffin-Modica} 

\noindent Let $u \in \mathcal{C}^2(\R^N)$ be a solution of \eqref{eqAC}. Then 

\smallskip


\beq 
\frac{|\nabla u|^2}{2}  \le \frac{(1-u^2)^2}{4} + c_u \qquad {\rm on} \qquad \R^N, 
\eeq
\noindent where 

\beq 
c_u : = \sup_{t \in  [\inf u, \sup u]} \Big [ - \frac{(1-t^2)^2}{4} \Big ] \le 0. 
\eeq

\end{theorem}

\bigskip

We are now ready to define the notion of {\it energy} of an entire solution of the Allen-Cahn equation \eqref{eqAC}. This will be the natural outcome of the following 

\medskip

\begin{theorem} {\textbf{(The energy)}} \label{Energia}

\noindent Let $u \in \mathcal{C}^2(\R^N)$ be a  solution of \eqref{eqAC}. 

\item[i)] For every $ x \in \R^N$, the function $ \phi_{u,x} (R) : = \Big [ \frac{1}{R^{N-1}}  \int_{B(x,R)} \frac{\vert \nabla u\vert^2}{2} + \frac{(1-u^2)^2}{4} \Big ] $ is non-decreasing on $(0, +\infty)$. 

\medskip
 
\noindent Therefore, the function $ L_u : \R^N \rightarrow [0, +\infty]$ given by

\beq\label{Energia(x)}
L_u(x) : = \lim_{R\to +\infty } \phi_{u,x} (R) 
\eeq

\noindent is well-defined for every $ x \in \R^N$. 

\smallskip

\item[ii)] The function {\underline {$L_u$ is constant}}.

\smallskip

\item[iii)]  $L_u$ is translation-invariant and O(N)-invariant, i.e., if $ \rho $ is either a translation or an orthogonal transformation and $u_{\rho}(x) := u(\rho x)$ for every $ x \in \R^N, $ then $L_u = L_{u_{\rho}}$. 

\smallskip

Also, $ L_{-u} = L_u $. 

\smallskip

\item[iv)] $L_u = 0$ if and only if $ u \equiv \pm 1$. 

\end{theorem}

\medskip

In view of item ii) of the above result we give  the following 

\medskip

\begin{defin}\label{DefEnergia}
We call {\it energy of u} the quantity $ L_u$. We shall say that $u$ has {\it finite energy} if $ L_u < +\infty.$ 
\end{defin}

\medskip

Often in the sequel, we shall use also the notation $L_u(N)$ to denote the {\it energy} $ L_u$. This will be done  when we want to emphasize the role played by the Euclidean dimension $N$. 

\medskip

\begin{proof} [Proof of Theorem \ref{Energia}]  

Since $u$ is bounded by item i) of Theorem \ref{Reg-univBound}, we can apply the monotonicity formula proved by L. Modica (cf. item ii) of Theorem \ref{Modica}). This result tells us that the function $ \phi_{u,0} $ is non-decreasing on $ (0, +\infty)$. Hence, we get that $L_u(0)$ is well-defined and belongs to $ [0, +\infty]$. Moreover, for every $ y \in \R^N$ the function $v(x) := u(x+y)$ is still a solution of \eqref{eqAC} and $ \phi_{u,y} = \phi_{v,0}$, therefore $L_u(y) = L_v(0) \in  [0, +\infty]$ and this value is achieved in a monotonically non-decreasing way (again by L. Modica's result). This proves item i).

To prove item ii) we observe that for every $x,y \in \R^N$ and every $ R>0$, we have $B(x,R) \subset B(y, \vert x-y\vert +R) $ and thus 

$$ \frac{1}{R^{N-1}}  \int_{B(x,R)} \frac{\vert \nabla u\vert^2}{2} + \frac{(1-u^2)^2}{4} \le \frac{1}{R^{N-1}}  \int_{B(y, \vert x-y\vert +R) } \frac{\vert \nabla u\vert^2}{2} + \frac{(1-u^2)^2}{4} = $$

$$  \frac{({R + \vert x - y\vert})^{N-1}}{R^{N-1}} \frac{1}{({R + \vert x - y\vert})^{N-1}} \int_{B(y, \vert x-y\vert +R) } \frac{\vert \nabla u\vert^2}{2} + \frac{(1-u^2)^2}{4}.$$

Passing to the limit, as $R \to +\infty$, in the latter we obtain  $ L_u(x) \le L_u(y)$. The conclusion then follows by exchanging the roles of $x$ and $y$ in the previous argument. 

First part of item iii) is an immediate consequence of item ii), the change of variable in the integral and the fact that $\rho$ is an Euclidean isometry. The second part of item iii) follows since $\phi_{-u,x} = \phi_{u,x} $ for every $x \in \R^N$. 

If $ u \equiv \pm 1$, then clearly $L_u = 0$. If $L_u = 0$, then by the monotonicity formula of L. Modica (cf. Theorem 1.1 in \cite{ModicaMonot}) and item ii) we have $ 0 = L_u = L_u (0) $. This yields $ \frac{\vert \nabla u\vert^2}{2} = \frac{(1-u^2)^2}{4} =0$ on $\R^N$, which entails $ u \equiv \pm 1$.   

\end{proof}

Now we turn to some compactness and lower-semicontinuity results about the set of solutions to the Allen-Cahn equation. 

\medskip

\begin{theorem} {\textbf{(Compactness and lower-semicontinuity)}} \label{Compattezza}
Let $ N \ge 1$ and set

$${\cal S} : = \{ \, u \in \mathcal{C}^2(\R^N) \quad : \quad u \quad {\text {is a  solution of \eqref{eqAC} }} \, \}. $$

\noindent Then, 

\item[i)] for every $ \kappa \ge 2$ the set ${\cal S}$ is compact in $\mathcal{C}^{\kappa}_{loc}(\R^N)$,

\item[ii)]  for any sequence $(u_n)_{n\ge 1}$ of ${\cal S}$ 
there are a subsequence $(u_{n_k}) \subset (u_n) $ and $u_{\infty}\in {\cal S}$ such that  

\beq\label{cvcomp}
u_{n_k} \rightarrow u_{\infty}  \quad in \quad \mathcal{C}^{\kappa}_{loc}(\R^N)
\eeq

and

\beq\label{scicomp}
L_{u_\infty} \le \liminf_{k\to+\infty} L_{u_{n_k}} 
\eeq

\end{theorem}

\begin{proof}

By i) of Theorem \ref{Reg-univBound} the set ${\cal S}$ is contained in the unit ball of $ L^{\infty}(\R^N) $. Then, for every $m \ge 2$, standard interior elliptic estimates imply the existence of a constant $C(m,N)>0$, depending only on $m$ and $N$, such  that 

\beq\label{stime-derivate} 
\forall u \in {\cal S} \qquad \max_{ \vert \alpha \vert \le m} \Vert \partial^{\alpha} u \Vert_{L^{\infty}(\R^N)} \le C(m,N).
\eeq

From the latter, the Ascoli-Arzel\`a Theorem  and a standard diagonal procedure we get the required compactness of ${\cal S}.$

Let $(u_n)_{n\ge 1}$ be a sequence of ${\cal S}$ (i.e., such that $u_n$ belongs to ${\cal S}$, for $n \ge 1$).  By the previous step, up to a subsequence, $ u_n \rightarrow u_{\infty}$ in  $\mathcal{C}^k_{loc}(\R^N)$ for some $ u_{\infty} \in {\cal S}$. By the monotonicity formula, for every $n \ge 1$ and every $ R>0$ 

\beq
\frac{1}{R^{N-1}}  \int_{B(0,R)} \frac{\vert \nabla u_n\vert^2}{2} + \frac{(1-u_n^2)^2}{4} \le L_{u_n} ,
\eeq
hence

\beq
\liminf_{n\to+\infty} \frac{1}{R^{N-1}}  \int_{B(0,R)} \frac{\vert \nabla u_n\vert^2}{2} + \frac{(1-u_n^2)^2}{4} \le \liminf_{n\to+\infty} L_{u_n},
\eeq
and using \eqref{cvcomp} we get 
 
\beq
\frac{1}{R^{N-1}}  \int_{B(0,R)} \frac{\vert \nabla u_{\infty}\vert^2}{2} + \frac{(1-u_{\infty}^2)^2}{4} \le \liminf_{n\to+\infty} L_{u_n}.
\eeq
The desired conclusion \eqref{scicomp} then follows by applying Modica's monotonicity formula and the definition of energy. 

\end{proof}

\medskip

\begin{theorem}{\textbf{(Compactness II)}} \label{Compattezza2}
Let $ N \ge 1$, $ L \in [0,+\infty)$ and  set

$${\cal S}_L : = \{ \, u \in {\cal S} \quad : \quad L_u \le L \, \}, $$

$${\cal X}_L : = \{ \, u \in {\cal S} \quad : \quad L_u = L \, \}. $$

\smallskip

\noindent Then, for every $ k \ge 2$, the set ${\cal S}_L$ is compact in $\mathcal{C}^k_{loc}(\R^N)$ while ${\cal X}_L $ is relatively compact in $\mathcal{C}^k_{loc}(\R^N)$.
Moreover,

\smallskip

\item[i)] ${\cal S}_0 = {\cal X}_0 = \{ \, u \equiv \pm 1 \, \}.$

\smallskip

\item[ii)]  The sets ${\cal S}, {\cal S}_L$ and $ {\cal X}_L$ are translation-invariant and O(N)-invariant, i.e., if $ \rho $ is either a translation or an orthogonal transformation and $u \in {\cal S}$ (resp. ${\cal S}_L, {\cal X}_L)$, then $ u_{\rho}(x) := u(\rho x)$ for every $ x \in \R^N, $  belongs to  ${\cal S}$ (resp. ${\cal S}_L, {\cal X}_L)$ too. 

\smallskip

\item[iii)]  $ {\cal X}_L \subset {\cal S}_L \subset  {\cal S}$, for $L \in (0,  +\infty) $ and the inclusions are strict.  

\medskip

Finally, $ {\cal S}_L \subseteq  {\cal S}_{L'} \, $ if $ \,\, 0 \le L \le L' < + \infty$.

\end{theorem}

\begin{proof}

The compactness of ${\cal S}_L $ and the relative compactness of ${\cal X}_L $ are an immediate consequence  of Theorem \ref{Compattezza}.  

Item i) follows from iv) of Theorem \ref{Energia}, while ii) is implied by iii)  of Theorem \ref{Energia}. By definition we have $ {\cal X}_L \subset {\cal S}_L \subset  {\cal S}$. The inclusions are strict since the constant solution $u \equiv 0$ belongs to ${\cal S} \setminus  {\cal S}_L $ and the constant solutions  $u \equiv \pm 1$ belong to $ {\cal S}_L \setminus  {\cal X}_L $, for $L \in (0,  +\infty) $.  The last claim follows immediately from the definiton of ${\cal S}_L $. 

\end{proof}

\begin{remark}\label{Rk1}
For $ L>0$, the set ${\cal X}_L $ is not, in general, closed in the topology of $\mathcal{C}^k_{loc}(\R^N)$. Indeed, there might be a {\it loss of energy}  when passing to the limit, as illustrated by the following example. Let $h = h(s) : = \tanh( \frac{s}{\sqrt{2}})$ be the one-dimensional monotone solution of \eqref{eqAC} in $\R$ and consider the function $ u=u(x) = h(x_N), \, x \in \R^N$.  It easily seen that $u$ solves \eqref{eqAC} in $\R^N$ for any $ N\ge2$. As shown in Appendix \ref {sec : appendix A}, $u$ belongs to the set $ {\cal X}_L $, with $ L = \frac{2\sqrt{2}}{3} \omega_{N-1}$, where $ \omega_{N-1} $ is the volume of the unit ball in $\R^{N-1}$.  Set $u_n(x) = u(x_1,..., x_N + n) $, for $ n\ge 1$. Thanks to Theorem \ref{Compattezza2} and to the monotonicity of the profile $h$,  $(u_n)_{n\ge 1}$ is a  sequence of ${\cal X}_L$ converging in $\mathcal{C}^k_{loc}(\R^N)$ to the constant solution $ u\equiv1$.  On the other hand, the latter function does not belong to  ${\cal X}_L$. 
\end{remark}

In the sequel we shall focus on {\it finite energy solutions} to the Allen-Cahn equation \eqref{eqAC} (cf. Definition \ref{DefEnergia}). As we shall see in the sequel, there are many natural examples of such a kind of solutions : the one-dimensional monotone profiles $ u(x) = \pm \tanh\big (\frac{(\nu\cdot x) +\alpha}{\sqrt 2} \big)$, where $\alpha \in \R$ and $ \nu $ is a unit vector of $\R^N$ (see the Appendix \ref {sec : appendix A}),  the  local minimizers (see Section \ref{sec : locmin}), the quasi-minimizers (see \cite{FVquasimin} and Section \ref{sec : locmin}), the saddle-type solutions constructed in \cite{DFP, ACM, CT1,FLP} (see also the Appendices \ref {sec : appendix C}, \ref {sec : appendix D} and \ref {sec : appendix D3}), the two-dimensional $2k-$ended solutions  built in \cite{DKPW}\footnote{Personal communication : M. Kowalczyk and F. Pacard, Santiago (Chile), April 2015}, and many others. 
 
Furthermore, finite energy solutions naturally arise in many different areas of pure and applied sciences. For instance, they are used in the description of the process of phase separation/transition in iron alloys (cf. for instance \cite{ACahn, FVstate}). They are also used to describe the early Universe in cosmology \cite{GT, Car}.  

\smallskip

Since the complete classification of the Allen-Cahn equation in the one-dimensional case is given in  Appendix \ref {sec : appendix A}, in the rest of the paper we shall consider the higher dimensional case $ N \ge 2$.  

\medskip

Also, from now on, for any $u \in {\cal S} $ and every $\eps>0$ we denote by $u_{\eps}$ the rescaled function defined by 

\beq\label{sol-riscalata}
\forall \, x \in \R^N \qquad u_\eps(x):=u(x/\eps),
\eeq

\noindent which is a solution of the rescaled problem 

\beq\label{eq-riscalata}
-\eps\Delta u_\eps = \eps^{-1}(u_\eps - u_\eps^3)  \qquad \text{in} \qquad \R^N.
\eeq

\section{Local minimizers}\label{sec : locmin}

\medskip

Given a bounded open set~$\Omega\subset\R^N$, we 
consider the energy functional associated to~\eqref{eqAC},
namely

\beq
{\mathcal{E}}(u,\Omega):=\int_\Omega\frac{|\nabla u|^2}{2}
+  \frac{(1-u^2)^2}{4}. 
\eeq

We say that~$u$ is a {\it local minimizer} if, for any
bounded open set~$\Omega$ and any~$\varphi\in \mathcal{C}^1_c(\Omega)$,
we have that

\beq
{\mathcal{E}}(u,\Omega)\le {\mathcal{E}}(u+\varphi,\Omega). 
\eeq

\noindent Clearly, any local minimizer is a smooth solution of equation \eqref{eqAC}. 

\medskip

More generally (see \cite{FVquasimin}), a solution $u \in \mathcal{C}^2(\R^N) $ of the Allen-Cahn equation \eqref{eqAC} is said to be 
a {\it quasi-minimizer} (or a {\it Q-minimizer}) if, for some~$Q\ge1$, one has that

\beq
{\mathcal{E}}(u,\Omega)\le Q\,{\mathcal{E}}(u+\varphi,\Omega)
\eeq
for any bounded open set~$\Omega$ and any~$\varphi\in \mathcal{C}^1_c(\Omega)$. 

\medskip

{\it Clearly, for $ Q=1 $, a $Q$-minimizer is a local minimizer. }

\medskip

A classical comparison argument yields that a local minimizer $u$ satisfies the following energy bound : 

\beq \label{e bound}
\forall R >0 \qquad {\mathcal{E}}(u,B(0,R))\le CR^{N-1},
\eeq

\smallskip

\noindent for some constant $C=C(N,u)>0$ depending only on $N$ and $u$.   

In turn, the energy bound \eqref{e bound} implies that the {\it local minimizer} $u$ has {\it finite energy} $L_u$ with $ L_u \le C(N,u) <+\infty$, that is, $u$  belongs to the set ${\cal S}_C \subset {\cal S}$.

\medskip

As a matter of fact, the energy bound \eqref{e bound} is valid for quasi-minimizers too (see Lemma~10 in~\cite{FVquasimin}).  Therefore,  {\it Q-minimizers} are {\it finite energy solutions} of the Allen-Cahn equation \eqref{eqAC} too.  

\medskip

The first result of the present section is 

\bigskip

\begin{theorem} {\textbf{(Structure Theorem for local minimizers)}} \label{Teo-minloc}
Assume $N \ge 2$. Let $u$ be a non-constant local minimizer and let $(\eps_n)_{n\ge1}$ be a sequence of positive real numbers such that $ \eps_n \to 0$. Then

\smallskip

\item[i)]  there are a subsequence of $(\eps_n)_{n\ge1}$ (not relabeled) and a minimal \footnote{in this section \textit{minimal} is in the sense of De Giorgi, i.e., a local minimizer for the perimeter.} cone $ C \subset \R^N $ such that $0 \in \partial C$ and

\beq\label{convL1-min} 
u_{\eps_n} \longrightarrow \, \,  {\bf 1}_C - {\bf 1}_{\R^N \setminus C} \quad {\text {a.e. \,\, in }} \quad \R^N 
\quad {\text {and \quad in}} \quad L^1_{loc}(\R^N),
\eeq

\smallskip

\item[ii)] 
\beq\label{energiadeiminimi}
L_u = \frac{2\sqrt{2}}{3} \omega_{N-1} \Theta(C,0) = \frac{2\sqrt{2}}{3} P(C,B(0,1)),
\eeq

\noindent  where $ \Theta(C,0) $ denotes the density of the minimal cone C at the origin and $P(C,B(0,1))$ denotes the perimeter of the cone C in the open unit ball $B(0,1).$

\smallskip

\item[iii)] For $ N \le 7$,  the cone $C$ is a half-space and $ L_u = \frac{2\sqrt{2}}{3} \omega_{N-1} $.

\smallskip

\item[iv)] For $ N \ge 8 $ we have 

\beq\label{stime-energia-minloc}
\frac{2\sqrt{2}}{3} \omega_{N-1}  \le L_u \le  \frac{2\sqrt{2}}{3} \frac{N\omega_{N}}{2} =  \frac{2\sqrt{2}}{3} \omega_{N-1} 
\Big (\frac{N \omega_{N}}{2 \omega_{N-1}} \Big ),   
\eeq
where $ \omega_{k} $ is the volume of the unit ball in $\R^{k}$.  

\smallskip

\item[v)] For every unit vector $\nu \in \R^N$ and every $\alpha \in \R$ the function defined by $ v(x) = \pm \tanh\big (\frac{(\nu\cdot x) +\alpha}{\sqrt 2} \big)$ is a local minimizer with $ L_u = \frac{2\sqrt{2}}{3} \omega_{N-1} $.

\medskip

{\sl In other words, the one-dimensional monotone profiles attain the {\rm minimum positive energy level} in the set of all non constant local minimizers. 

\smallskip

In the sequel we shall denote by $ {\cal M} $ the set of non constant local minimizers.}

\end{theorem}

\bigskip

\begin{remark}\label{Rk2}
We note that in item i),  different (sub-)sequences might converge to different minimal cones. Nevertheless, item ii) tells us that all the limit minimal cones have the same density at the origin!   
\end{remark}

\bigskip

\begin{proof} [Proof of Theorem \ref{Teo-minloc}]  

Since $u$ is a local minimizer, the rescaled functions $u_{\eps} $ are local minimizers for \eqref{eq-riscalata}. Moreover, since $u$ has finite energy the rescaled energy functional 

\beq\label{energia-riscalata}
{\mathcal{E}_{\eps}}(u_{\eps},\Omega):=\int_\Omega\frac{\eps \vert \nabla u_{\eps} \vert^2}{2}
+\frac{(1 - u_\eps^2)^2}{4\eps},
\eeq

\noindent satisfies, for every $R>0$ and every $ n \ge 1$,  

\beq \label{e bound-2} 
{\mathcal{E}_{\eps_n}} (u_{\eps_n},B(0,R)) \le L_u R^{N-1}.
\eeq
As $\eps \to 0$, a classical result of L. Modica and S. Mortola \cite{ModMort} ensures that the rescaled functionals ${\mathcal{E}_{\eps}}(\cdot, \Omega)$, $\Gamma-$ converge to a constant multiple of the area functional, i.e. to 
\beq\label{area functional}
c_{\star} P(\cdot, \Omega)
\eeq
where $c_{\star} = \int_{-1}^1 \frac{1-s^2}{\sqrt 2} ds = \frac{2\sqrt{2}}{3}$ and, for any  set $E$ of locally finite perimeter in $\R^N$, $P(E,\Omega)$ denotes the perimeter of $E$ in $\Omega$. 

Thus, (as also pointed out in \cite{AAC}) as a consequence of the local minimality of the rescaled solutions $u_{\eps}$ and the energy bound \eqref{e bound-2},  there are a subsequence of $(\eps_n)_{n\ge1}$ (not relabeled) and a set of locally finite perimeter $ E $ in $ \R^N $, which is a {\it local minimizer for the  perimeter}, such that 

\beq\label{convL1-riscalata} 
u_{\eps_n} \longrightarrow \, \,  {\bf 1}_E - {\bf 1}_{\R^N \setminus E} \quad  {\text {in}} \quad L^1_{loc}(\R^N)
\eeq
and
\beq\label{convMisura}
{\mathcal{E}_{\eps_n}} (u_{\eps_n},B(0,R)) \longrightarrow c_{\star} P(E, B(0,R)),
\eeq
where the latter holds true for every $R>0$ (with the exception of an at most countable set of values).  

Now, scaling back we get 

\beq\label{cv1}
{\eps_n}^{N-1} \int_{B(0,\frac{R}{\eps_n})} \frac{\vert \nabla u \vert^2}{2} +\frac{(1 - u^2)^2}{4} = {\mathcal{E}_{\eps_n}} (u_{\eps_n},B(0,R)) \longrightarrow c_{\star} P(E, B(0,R)),
\eeq
for every $R>0$, with the exception of an at most countable set of values.
On the other hand, we also have that  

\beq\label{cv2}
{\eps_n}^{N-1} \int_{B(0,\frac{R}{\eps_n})} \frac{\vert \nabla u \vert^2}{2} +\frac{(1 - u^2)^2}{4} = R^{N-1} 
\frac{1}{(\frac{R}{\eps_n})^{N-1}} \int_{B(0,\frac{R}{\eps_n})} \frac{\vert \nabla u \vert^2}{2} +\frac{(1 - u^2)^2}{4}
\eeq
for every $R>0$ and every $n \ge 1$. 

By Theorem \ref{Energia}, as $n \to +\infty,$ the right-hand side of the latter converges to $R^{N-1} L_u$ and so, the combination of  \eqref{cv1} and \eqref{cv2} yields  

\beq\label{formula-energia}
0< L_u =  \frac{2\sqrt{2}}{3} \frac{P(E,B(0,R))}{R^{N-1}} 
\eeq
for every $R>0$, with the exception of an at most countable set of values. 

Since $E$ is minimal, the latter implies that $\frac{P(E,B(0,R))}{R^{N-1}} $ is constant for $ R>0$ and thus $E$ is a  minimal cone $C$ with $ 0 \in \partial C$ (see for instance \cite{Giusti, GMM}). These facts and \eqref{formula-energia} prove items i) and ii).  When $N\le 7$, such a minimal cone is a half-space \cite{Alm1, Simons}, which proves item iii). The energy estimates \eqref{stime-energia-minloc} are then an immediate consequence of the following  well-known density estimates for local perimeter minimizers $E$ in $\R^N$ (see for instance \cite{Giusti, GMM}) :

\beq\label{stime-dens-insiemi-minimi}
\forall \, r >0, \quad \forall \, x \in \partial E \qquad \omega_{N-1} \le \frac{P(E,B(x,r))}{r^{N-1}} \le \frac{N \omega_N}{2}.
\eeq

Finally, item v) follows from the Theorem \ref{Energia1D} in Appendix \ref{sec : appendix A}.

\end{proof}

\medskip

{\it As a consequence of the above Theorem \ref{Teo-minloc} we have the following new one-dimensional symmetry result, which  completely characterizes the one-dimensional monotone profiles. }

\medskip

Let us recall that, for $ k \in \N$, a function is said to be {\it k-dimensional} (or, for short, $kD$) if it depends only on $k$ variables, up to rotations and translations of the Euclidean space $\R^N.  $

\bigskip

\begin{theorem} {\textbf{(0ne-dimensional symmetry)}} \label{1D-minloc}

\noindent Assume $N \ge 2$ and let $u$ be a non-constant solution of the Allen-Cahn equation \eqref{eqAC}. 

Then, $u$ is a 1D monotone profile if and only if $u$ is a local minimizer attaining the {\rm minimum positive energy level. } 
Namely, if and only if $u$ is a local minimizer with energy $ L_u = \frac{2\sqrt{2}}{3} \omega_{N-1} $.

In particular, for $ N\le 7$,  we have that $u$ is a 1D monotone profile if and only if $u$ is a local minimizer. 

\end{theorem}

\medskip

\begin{remark}\label{Rk3}
The last claim of Theorem \ref{1D-minloc} (i.e. the conclusion for $ N\le 7$, and only that one) was first proved by O. Savin \cite{S}. Thus, our Theorem \ref{1D-minloc} recovers and improves upon that result, since Theorem \ref{1D-minloc}  holds true for every dimension $ N \ge 2$. 
\end{remark}

\medskip

\begin{proof}[Proof of Theorem \ref{1D-minloc}]  
From v) of Theorem \ref{Teo-minloc} we already know that any one-dimensional monotone profile  $v=v(x) = \pm \tanh\big (\frac{(\nu\cdot x) +\alpha}{\sqrt 2} \big)$ is a local minimizer with $ L_v = \frac{2\sqrt{2}}{3} \omega_{N-1} $. 

Let $u$ be a local minimizer such that $ L_u = \frac{2\sqrt{2}}{3} \omega_{N-1} $. Given a sequence $(\eps_n)_{n\ge1}$ of positive real numbers such that $ \eps_n \to 0$, Theorem \ref{Teo-minloc}  provides a minimal cone $C$, with $0 \in \partial C$, such that 

\beq\label{convL1-non} 
u_{\eps_n} \longrightarrow \, \,  {\bf 1}_C - {\bf 1}_{\R^N \setminus C} \quad  {\text {in}} \quad L^1_{loc}(\R^N),
\eeq

\beq
L_u = \frac{2\sqrt{2}}{3} P(C,B(0,1)) = \frac{2\sqrt{2}}{3} \frac{P(C,B(0,R))}{R^{N-1}}  \qquad \forall \quad R>0.
\eeq
Therefore, $0 \in \partial C$ and 

\beq\label{energiaminima}
\omega_{N-1} = \frac{P(C,B(0,R))}{R^{N-1}}  \qquad \forall \quad R>0
\eeq
and Theorem 2.1 in \cite{GMM} implies that $C$ is a half-space. 

To prove the 1D-symmetry we invoke Corollary 7 in \cite{FVquasimin}  (applied with $E = C$ and recalling that a local minimizer is a Q-minimizer with $ Q=1$). 

Item iii) of Theorem \ref{Teo-minloc} says that the energy of every non constant local minimizer $u$ must satisfy $ L_u = \frac{2\sqrt{2}}{3} \omega_{N-1}$,  when $ N \le 7$. The last claim then follows by the result established in the previous part of the proof. 

\end{proof}

\bigskip

As a by-product of the above analysis we proved that the value $\frac{2\sqrt{2}}{3} \omega_{N-1}$ gives the {\it minimum positive energy level} in the set of local minimizers (and that it is saturated  (only) by the one-dimensional monotone profiles). 

In the next result we prove that this energy level is indeed {\it isolated}. This conclusion is clearly true for $ N \le 7$, since in this case by item iii) of Theorem \ref{Teo-minloc} the only positive energy level is the minimum one.  This is why, in the next result, we restrict our attention  to dimensions 
$ N \ge 8.$ In this case, the fact that the {\it minimum positive energy level}  is isolated will be a consequence of the following {\it gap phenomenon}  for the set of local minimizers.  
 
\bigskip

\begin{theorem} {\textbf{(Gap phenomenon for local minimizers)}} \label{Gap-minloc}

\noindent Assume $N \ge 8$ and let $u$ be a non-constant local minimizer. Then 
\beq\label{Valore-gap-minloc}
\exists \quad \gamma=\gamma(N) 
\in \Big (1, \frac{N\omega_{N}}{2 \omega_{N-1}} \Big]
\eeq
such that 
\beq
L_u < \frac{2\sqrt{2}}{3} \omega_{N-1}\gamma(N) \quad \Longrightarrow \quad L_u = \frac{2\sqrt{2}}{3} \omega_{N-1}.
\eeq
Namely,  the {\it minimum positive energy level} is {\rm isolated}. 

\end{theorem}

\medskip

As a consequence of the above theorem we have the following improvement of Theorem \ref{1D-minloc}

\medskip

\begin{corollary}\label{Improved1D-minloc}
Assume $N \ge 8$ and let $u$ be a non-constant solution of the Allen-Cahn equation \eqref{eqAC}. 

\noindent Then, $u$ a 1D monotone profile if and only if $u$ is a local minimizer with energy 
 
 $$ L_u < \frac{2\sqrt{2}}{3} \omega_{N-1} \gamma(N).$$ 

\end{corollary}

\begin{proof}
We need to prove the result only for $ N \ge 8$. The assumption $L_u <\frac{2\sqrt{2}}{3} \omega_{N-1} \gamma(N)$  implies $ L_u = \frac{2\sqrt{2}}{3} \omega_{N-1}$ by Theorem \ref{Gap-minloc}. The desired conclusion  then follows from Theorem \ref{1D-minloc}.
\end{proof}

\medskip

\begin{proof} [Proof of Theorem \ref{Gap-minloc}]  
From Theorem 3.2 in \cite{GMM} we know that 
\beq\label{def-alpha(N)}
\alpha(N) : = \inf \{ \, P(C,B(0,1)) \, : \, C {\text { is a minimal cone with }} P(C,B(0,1))>\omega_{N-1} \, \}
\eeq
is a {\it minimum}.  The claim then follows from item ii) of Theorem \ref{Teo-minloc} by taking $ \gamma = \frac{\alpha(N)}{\omega_{N-1}}$ in \eqref{Valore-gap-minloc}.

\end{proof}

\medskip

We recall that it is still an open question the existence (or not) of a non-constant and non-planar local minimizer on $\R^8$, i.e., the existence (or not) of a local minimizer in $\R^8$, which is not 1D. We also recall that in Theorem 1.7 of \cite{FVtams} we proved 

\begin{theorem} {\textbf{(1D-symmetry for monotone minimizers in $\R^8$)}} \label{monot-1D-minloc}

\noindent Every monotone local minimizer in $\R^8$ must be 1D. 

\end{theorem}

\medskip


Therefore, if a non-1D local minimizer exists in $\R^8$, it must be non monotone. 

\smallskip

On the other hand, for $ N \ge 9 $ the situation is different, since  the solution built by M. del Pino, M. Kowalczyk and J. Wei in \cite{DKW} for $ N\ge9$ is a {\it monotone local minimizer} which {\it is not} 1D. 

\medskip

The following result makes the connection between the (non-)existence of non-1D local minimizers in $\R^8$ and the asymptotic behaviour of monotone solutions in  $\R^9$. 

\medskip 
 
\begin{proposition}\label{minloc8D} If the only local minimizers in $\R^8$ are 1D then every monotone local minimizer in $\R^9$ must satisfy (up to isometries) 

\begin{equation}\label{LI}
\lim_{x_9\rightarrow-\infty} u(x',x_n)=-1
\quad {\mbox{ and }} \quad
\lim_{x_9\rightarrow+\infty} u(x',x_n)=1.
\end{equation}
Therefore, one of the two following possibilities holds :

\smallskip

\item[i)] either there exists a non-1D local minimizer in $\R^8$, 

\item[ii)] or every monotone local minimizer in $\R^9$ must satisfy (up to isometries) 

\begin{equation}\label{LI2}
\lim_{x_9\rightarrow-\infty} u(x',x_n)=-1
\quad {\mbox{ and }} \quad
\lim_{x_9\rightarrow+\infty} u(x',x_n)=1.
\end{equation}
\end{proposition}

\medskip

\begin{remark}
Note that the solution $U$ built in \cite{DKW} for $ N=9$ is a {\it monotone local minimizer} satisfying \eqref{LI}. More precisely, this solution satisfies $ \frac{\partial U}{\partial x_9} >0$ in $\R^9$. 
\end{remark}
\begin{remark}
We observe that, the proof of the above Proposition \ref{minloc8D} does not exclude that both the possibilities i) and ii) may occur at the same time. 
\end{remark}

\medskip

\begin{proof} 
Let $u$ be any monotone local minimizer in $\R^9$ then, up to isometries, we can suppose that $ \frac{\partial u}{\partial x_9} >0$ in $\R^9$.  Therefore, the following two limit profiles
\begin{equation}\label{profiles}\begin{split}
&\underline u(x'):=\lim_{x_9\rightarrow-\infty} u(x',x_9), 
\\ 
&\overline u(x'):=\lim_{x_9\rightarrow+\infty} u(x',x_9),
\end{split}\end{equation}
exists and are solutions of the Allen-Cahn equation \eqref{eqAC} on $\R^8$ and by Theorem 1.3 of \cite{JM}  we  know that both $\underline u$ and $\overline u$ are local minimizers on $ \R^8$. 
Under our assumption we have that the limit profiles $\underline u$ and $\overline u$ are 1D and  therefore we can invoke Proposition 2.13 of \cite{FVtams} to deduce that $\underline u = - 1$ and $\overline u = + 1$. This proves \eqref{LI}. The second part of the claim is then a consequence of the above result. 

\end{proof}

\medskip

Now we would like to give a different proof of the second part of the previous result. This proof makes use of the notion of energy introduced in this work and the results of the present chapter.
We hope that this simple example could illustrate both the robustness and the flexibility of the concepts of energy and energy levels.  

\medskip

{\it Second proof of the second part of Proposition \ref{minloc8D}}.\begin{proof} 
Let $u$ be any monotone local minimizer in $\R^9$ then, up to isometries, we can suppose that $ \frac{\partial u}{\partial x_9} >0$ in $\R^9$. Therefore the limits profiles defined in \eqref{profiles} are solutions of the Allen-Cahn equation \eqref{eqAC} on $\R^9$ and, by the definition of local minimizer, we immediately see that  both $\underline u$ and $\overline u$ are local minimizers on $ \R^9$ (not $\R^8$) and in particular they have finite energy. 

At this point we have two possibilities : either for every monotone local minimizer $u$ we have that : $L_{\underline u}\le  \ell_{min}$ and $L_{\overline u} \le \ell_{min}$ or there is a monotone local minimizer $u$ such that at least one of the two profiles, says ${\overline u}$, satisfies $L_{\overline u} > \ell_{min}$. 

In the first case we necessarily have $L_{\underline u} \in \{\, 0, \ell_{min} \, \} $ and $L_{\overline u} \in \{\, 0, \ell_{min} \, \} $. Therefore, both the profiles are 1D, by Theorem \ref{Energia} and Theorem \ref{1D-minloc}. An application of Proposition 2.13 of \cite{FVtams} then provides : $\underline u = - 1$ and $\overline u = + 1$.  That is item ii) holds true. 

In the second case  ${\overline u}$ cannot be 1D, i.e., it cannot be of the form $u(x) = u_0 (\nu\cdot x)$ for some one-dimensional solution of the Allen-Cahn equation on $\R$. Indeed, if it were true, $u_0$ would be a one-dimensional finite-energy solution of the ode \eqref{eqAC1D} and thus, by Theorem \ref{Teo-clas1D},  necessarily of the form $u_0 \equiv \pm 1$ or $u_0(s) = \pm \tanh \Big ( \frac{s +\alpha}{\sqrt{2}} \Big )$, for every $s \in \R$ and some $\alpha \in \R$.  These two possibilities would give $L_u \in \{0, \ell_{min}\}$(cf. Theorem \ref{Energia1D}),  which is clearly impossible  because $ \ell_{min} < L_{\overline u} \le L_u$, by item ii) of Theorem \ref{Compattezza}. Therefore, ${\overline u}$ is a local minimizer in $\R^8$ which is not planar, and item i) holds true. 

\end{proof}

\medskip

Next result parallels some results in the study of the Bernstein property for entire minimal graphs. It is inspired by methods introduced by W. Fleming, E. De Giorgi and M. Miranda (cf. for instance \cite{Giusti, GMM}). It describes the  structure of the limit minimal cone provided by Theorem \ref{Teo-minloc} in the case of monotone local minimizers  which are not 1D (i.e. not planar). Again, we shall work only in dimension $N\ge9$. It will be also useful later on in the present work.

\medskip

\begin{theorem} {\bf (Structure Theorem for monotone local minimizers)} \label{monot-cono-cilindro}  Assume $N \ge 9$ and let $u$ be a local minimizer satisfying $ \frac{\partial u}{\partial x_N} >0$ in $\R^N$.  

\noindent Assume also the $u$ is not one-dimensional. 

\smallskip

\noindent Then the limit minimal cone $C$ provided by Theorem \ref{Teo-minloc} is a cylinder of the form  

\beq\label{conf-cilindro}
C = C^{'} \times \R
\eeq
where $C^{'} \subset \R^{N-1}$ is a minimal cone which is singular at the $0' \in \R^{N-1}$.  
\end{theorem}

\begin{proof}
By the monotonicity assumption the limits profiles defined by
\begin{equation}\label{N-profiles}\begin{split}
&\underline u(x'):=\lim_{x_N\rightarrow-\infty} u(x',x_N), 
\\ 
&\overline u(x'):=\lim_{x_N\rightarrow+\infty} u(x',x_N),
\end{split}\end{equation}
are solutions of the Allen-Cahn equation \eqref{eqAC} on $\R^N$ and, by the definition of local minimizer, we immediately see that  both $\underline u$ and $\overline u$ are local minimizers on $ \R^N. $

Now we have two possibilities, either $\underline u = -1 $ and $\overline u = 1$, or at least one of the two profiles is non constant, say $\overline u$. In the second case, $u$ and $\overline u$ are ordered non-constant local minimizers and therefore item i) of Proposition \ref{Teo-Ord}  yields $ C = C_{\overline{u}}$. Let us prove that $C_{\overline{u}}$ is of the form \eqref{conf-cilindro}. To this end we observe that 

\beq
\forall x^{'} \in \R^{N-1}, \, \forall s,t \in \R \qquad  {\overline u}_{\eps_n}(x^{'},s) = {\overline u}_{\eps_n}(x^{'},t) 
\eeq
and so
\beq
\lim_{n\to\infty} {\overline u}_{\eps_n}(x^{'},s) = 1 \quad \Longleftrightarrow \quad \lim_{n\to\infty} {\overline u}_{\eps_n}(x^{'},t) = 1.
\eeq
The latter immediately shows that $C_{\overline{u}}$ is a cylinder of the form \eqref{conf-cilindro} with $C^{'} \subset \R^{N-1}$ a cone with $0 \in \partial C^{'}$.  On the other hand $C_{\overline{u}} \subset \R^N $ is a minimal cone and thus the minimality of $C^{'}$ then follows by applying, for instance, Proposition 9.9 of \cite{Giusti}.  Now we observe that $C^{'}$ must be singular at the origin $0' \in \R^{N-1}$. If not, it would be an half-space of $\R^{N-1}$ and so $C$ would be an half-space of $\R^N$ thanks to \eqref{conf-cilindro}. We can then apply Theorem \ref{1D-minloc} (or alternatively Corollary 7 in \cite{FVquasimin}) to get  that $u$ is necessarily a 1D monotone profile.  A contradiction. 

\smallskip

In the first case, every level set of $u$ is a smooth entire graph in the $e_N-$direction. In particular the zero level set of $u$ satisfies this property and thus there is a smooth function $ \psi : \R^{N-1} \rightarrow \R$ such that 

\beq
\{ \,  u < 0 \, \} = \{ (x^{'},x_N) \in \R^{N-1} \times \R \, : \, x_N < \psi(x^{'}) \, \} 
\eeq
and so
\beq\label{sottografici-risc}
\{ \,  u_{\eps_n} < 0 \, \} = \{ (x^{'},x_N) \in \R^{N-1} \times \R \, : \, x_N < \psi_n(x^{'}) \, \} : = U_n
\eeq
where 
\beq
\psi_n (x^{'}) : = \eps_n  \psi \Big (\frac{x^{'}}{\eps_n} \Big).  
\eeq

By item i) of Theorem \ref{Teo-minloc} and Lebesgue dominated convergence theorem we immediately have that

\beq\label{cv-sottografici}
 {\bf 1}_{U_n} \longrightarrow {\bf 1}_{\R^N \setminus C} \quad  {\text {in}} \quad L^1_{loc}(\R^N).
\eeq
The latter and Lemma 16.3 of \cite{Giusti} tell us that $\R^N \setminus C$ is itself a subgraph of a measurable function 
$ \psi_{\infty} : \R^{N-1} \rightarrow [- \infty, + \infty]$ which is the a.e.-limit of $\psi_n$ (up to a subsequence) . In particular, according to \cite{Giusti, GMM}, we have that $ \psi_{\infty} $ is a quasi-solution of the minimal surface equation, since $C$ is a minimal cone.  Also, $ \psi_{\infty} $ is positive homogeneous of degree 1, since $C$ is a cone with vertex at the origin. Therefore the sets 
\beq\label{P-N}
P : = \{ x^{'} \in \R^{N-1} \, : \, \psi_{\infty}(x^{'}) = + \infty  \}, \qquad N : = \{ x^{'} \in \R^{N-1} \, : \, \psi_{\infty}(x^{'}) = - \infty  \}
\eeq
are both minimal cones in $\R^{N-1}$ with vertex at the origin. 

Now we distinguish two cases : either $ P = N = \emptyset$ or not. In the first case, the quasi-solution  $\psi_{\infty}$ is locally  bounded and therefore it is a classical solution.  Then, the minimal cone $C$ satisfies 
\beq
C + e_N \subset C \quad {\text and } \quad  C + e_N \not \equiv C 
\eeq
and it must be a half-space thanks to Lemma 2.3 of \cite{GMM}. But then, we can invoke Theorem \ref{1D-minloc} once again to get that $u$ is necessarily a 1D monotone profile. A contradiction. 

In the second case, up to replace $u=u(x)$ by $ -u(x^{'},- x_N)$ if necessary, we can suppose that $ N \not\equiv \emptyset. $  But then, by the maximum principle, the two minimal cones given by $N \times \R$ and $ C $ must coincide, since they satisfy $ N \times \R \subset C \subset \R^N $ and they have the origin as common vertex. This proves \eqref{conf-cilindro}. To prove that $C'$ is singular at the origin of $\R^{N-1}$ we repeat {\it verbatim} the argument used in the first part of the present proof. This concludes the proof of the theorem. 

\end{proof}

\medskip

The inspection of the previous proof immediately provides the following result for monotone solutions $u$ with constant profiles $\underline u$ and $\overline u$.

Note that we {\it do not assume} any local minimality property on the solution $u$.

\medskip

\begin{theorem} {\textbf{(Structure Theorem for monotone solutions with constants profiles)}} \label{monot-con-limiti-cost}  

\noindent Assume $N \ge 9$ and let $u$ be a non-constant solution satisfying $ \frac{\partial u}{\partial x_N} >0$ in $\R^N$ and

\begin{equation}
\lim_{x_N\rightarrow-\infty} u(x',x_n)=-1
\quad {\mbox{ and }} \quad
\lim_{x_N\rightarrow+\infty} u(x',x_n)=1.
\end{equation}

\smallskip

\noindent Then, one of the two following possibilities holds :

\smallskip

\item[i)] either $u$ is a monotone 1D profile, 

\item[ii)] or $u$ is a local minimizer and the limit minimal cone $C,$ provided by Theorem \ref{Teo-minloc} is a cylinder of the form  

\beq\label{conf-cilindro-precisato}
C = C^{'} \times \R
\eeq
where $C^{'} \subset \R^{N-1}$ is a minimal cone singular at the origin $0' \in \R^{N-1}$. 
\end{theorem}

\medskip

\begin{proof}
In view of Theorem \ref{monot-cono-cilindro} we only need to prove that $u$ is a local minimizer of the Allen-Cahn equation \eqref{eqAC} in $\R^N$. This property is a consequence of Theorem 1.3 of \cite{JM}, because $u$ is monotone and the constant functions $ \pm 1$ are local minimizers.  

\end{proof}

\begin{remark}
The above theorem applies, in particular, to the solution $U$ built in \cite{DKW} by M. del Pino, M. Kowalczyk and J. Wei. 
It will be crucial to prove that the minimal cone $C^{'}$ is exactly the (solid) minimal Simons' cone $E_{4,4}$ (cfr. Appendix \ref {sec : appendix F}) and to obtain thus the exact value of the energy of this special solution. This will be accomplished in Appendix \ref {sec : appendix E}.

\end{remark}

\section{Finite energy solutions : the general structure and the energy levels}\label{sec : finite-energy}

\medskip

In this section we study the structure of the set of all {\it finite and positive energy} solutions to the Allen-Cahn equation \eqref{eqAC}, namely the set 

\beq\begin{split}
{\cal F}_p : = \{ \, u \in {\cal S} \quad : \quad 0 < L_u < +\infty \, \} = \bigcup_{L \in (0, +\infty)} {\cal X}_L 
\end{split}\eeq

We observe that the set $ {\cal F}_p $ is not empty since it contains the set ${\cal M} $ of all non constant local minimizers. 

\medskip

One of the consequences of the analysis developed in the previous section was the existence of the  {\it minimum positive energy level} in the set ${\cal M}$.  In the first result of the present section we shall prove that this result is indeed true in the largest set ${\cal F}_p$ of  {\it finite and positive energy} solutions to the Allen-Cahn equation \eqref{eqAC}. More precisely we have

\bigskip 
 
\begin{theorem} {\textbf{(Existence of the minimum positive energy level)}}\label{Teo-esist-livellomin}

\noindent Assume $N \ge 2$ and set

\beq
{\ell}_{min} : = \inf \{ \, L_u >0 \,\, : \,\, u \in {\cal S} \, \}.
\eeq
Then, 
\beq
{\ell}_{min} >0
\eeq
and
\beq
\exists \,\,  u \in {\cal F}_p \quad : \quad L_u = {\ell}_{min}.  
\eeq

\noindent We shall call this value the  {\rm minimum positive energy level}. 

\end{theorem}

\medskip

\begin{proof}
Set $I : = \inf \{ \, L_u \,\, : \,\, u \in {\cal F}_p \, \}$ and let $(v_n)_{n\ge 1}$ be a minimizing sequence. By item iii) of Theorem \ref{Reg-univBound} for every $ n \ge 1$ there is $ x_n \in \R^N$ such that $ v_n(x_n) =0$. Therefore, $ u_n (x) :=  v_n(x+x_n) $ is still a minimizing sequence, thanks to  item iii) of Theorem \ref{Energia}.  Then, up to a subsequence, 
$ u_n \rightarrow u_{\infty}$ in  $\mathcal{C}^k_{loc}(\R^N)$ for some $ u_{\infty} \in {\cal S}$ satisfying

\beq
L_{u_{\infty}}  \le \liminf_{n\to+\infty} L_{u_n} = I
\eeq
by Theorem \ref{Compattezza}, and also
\beq\label{zero}
u_{\infty} (0) = 0.
\eeq
On the other hand \eqref{zero} implies $ u \not \equiv \pm 1$ and then item iv) of Theorem  \ref{Energia} gives $ L_{u_{\infty}} >0$. This proves that $ u_{\infty} \in {\cal F}_p$ and thus $ I = L_{u_{\infty}} >0 $, concluding the proof. 

\end{proof}

\medskip

In view of \ref{Energia1D} in Appendix \ref {sec : appendix A}, any one-dimensional monotone profile belongs to $ {\cal X}_{\frac{2\sqrt{2}}{3} \omega_{N-1}} \subset  {\cal F}_p$. Thus,  $ {\ell}_{min} \le \frac{2\sqrt{2}}{3} \omega_{N-1}$. In the next result we prove that ${\ell}_{min}$ is exactly equal to $ \frac{2\sqrt{2}}{3} \omega_{N-1}$. In particular,  the  {\it minimum positive energy level} is the same as the  {\it minimum positive energy level} in the set of local minimizers and it is achieved by any 
one-dimensional monotone profile.

\bigskip 
 
\begin{theorem} {\textbf{(Value of the minimum positive energy level)}}\label{Teo-livellomin}

\noindent Assume $N \ge 2$. Then 

\beq
{\ell}_{min} = \frac{2\sqrt{2}}{3} \omega_{N-1}.
\eeq
Moreover, any one-dimensional monotone profile attains the {\rm minimum positive energy level}.  

\end{theorem}

\bigskip

The latter result is a consequence of the following Theorem \ref{Teo-finite-energy}, which may be seen as the anologue of Theorem \ref{Teo-minloc} in the more general setting of finite energy solutions.  To state and prove Theorem \ref{Teo-finite-energy} we need to recall some standard facts about varifolds. 

Following \cite{HT} and \cite{TW}, to each rescaled solution $u_{\eps}$ we associate the varifold $V_{\eps}$ on $\R^N$ defined by 

\beq\label{varifold-eps}
V_{\eps} (\phi) : =  \frac{2}{c_{\star}} \int_{ \vert \nabla u_{\eps} \vert >0}  \phi \Big (x, I - \nu_{\eps}(x) \otimes \nu_{\eps}(x) \Big )\frac{\eps \vert \nabla u_{\eps} \vert^2}{2} dx
\eeq
for each $ \phi \in C_c(G_{N-1}(\R^N))$. 

Here $ G_{N-1}(\R^N) : = \R^N \times G(N,N-1)$, with the product topology, where $ G(N,N-1)$ denotes the Grassman manifold of unoriented $(N-1)-$dimensional planes in $ \R^N$ (and each element $S \in G(N,N-1)$ is identified with the $N \times N$ matrix corresponding to the orthogonal projection of $\R^N$ onto $S$), 
$ \nu_{\eps} = \frac{\nabla u_{\eps}(x)}{\vert \nabla u_{\eps} (x) \vert} $, $I$ is the $ N \times N$ identity matrix and 
$\otimes$ is the tensor product.

\bigskip 
 
\begin{theorem}  {\textbf{(Structure Theorem for finite energy solutions)}} \label{Teo-finite-energy}

\noindent Assume $N \ge 2$. Let $u$ be a non-constant solution of \eqref{eqAC} with finite energy and let $(\eps_n)_{n\ge1}$ be a sequence of positive real numbers such that $ \eps_n \to 0$. 

Then there are a subsequence of $(\eps_n)_{n\ge1}$ (not relabeled), a set E of locally finite perimeter in $\R^N$, a Radon measure $ \eta \neq 0 $ on $ \R^N$ and a stationary, integral ($N$-1)-varifold $V \neq 0$ on $\R^N$ such that 

\smallskip

\item[i)]  

\beq\label{convL1} 
u_{\eps_n} \longrightarrow \, \,  {\bf 1}_E - {\bf 1}_{\R^N \setminus E} \quad {\text {a.e. \,\, in }} \quad \R^N 
\qquad {\text {and \qquad in}} \quad L^1_{loc}(\R^N),
\eeq

\beq\label{misura-limite} 
{\mathcal{E}_{\eps_n}} (u_{\eps_n}, \cdot) \longrightarrow \, \,  \eta \quad  {\text {locally weakly-$\star$ as measures}}, 
\eeq

\beq\label{varifold-limite} 
V_{\eps_n} \longrightarrow \, \,  V \quad  {\text {in the varifold sense}}. 
\eeq

\smallskip

\item[ii)] 
\beq\label{massaVarifold}
\Vert V \Vert = \frac{\eta}{c_{\star}}
\eeq

\noindent  where $ \Vert V \Vert$ denotes the weight measure associated to the limit varifold $V$ and $c_{\star} = \int_{-1}^1 \frac{1-s^2}{\sqrt 2} ds = \frac{2\sqrt{2}}{3}$.

\smallskip

\item[iii)] 

\beq\label{Vconica}
0 \in supp \Vert V \Vert \qquad  {\text { and \qquad $V$ is conical about $0$}}. 
\eeq
\noindent  Here $ supp \Vert V \Vert$ denotes the support of the weight measure $\Vert V \Vert $. 

\smallskip

\item[iv)] 
\beq\label{energiaV}
L_u = \frac{2\sqrt{2}}{3} \omega_{N-1} \Theta(V,0). 
\eeq
\noindent  where $ \Theta(V,0) $ denotes the density function of the varifold $V$  at the origin $0$. 

\end{theorem}

\bigskip

\begin{remark}\label{Rk4}
In item i),  different (sub-)sequences might converge to different varifolds. Nevertheless, item iv) tells us that all the limit stationary, integral and conical varifold have the same density at the origin!   
\end{remark}

\bigskip

\begin{proof}[Proof of Theorem \ref{Teo-livellomin}]
Since $ V $ is stationary its density function $ \Theta (V, \cdot)$ is upper semicontinuous. Then the integrality of $V$ and the fact that $0 \in supp \Vert V \Vert$ imply that $ \Theta(V,0) \ge 1$. The desired conclusion follows then from \eqref{energiaV}.

The last claim follows from the Theorem \ref{Energia1D} in Appendix \ref{sec : appendix A}.

\end{proof}

\medskip

\begin{proof}[Proof of Theorem \ref{Teo-finite-energy}] 

Since $u$ is an entire solution of the Allen-Cahn equation, the rescaled functions $u_{\eps} $ are solutions to \eqref{eq-riscalata}. In particular, they are critical points of the rescaled energy functional 

\beq \label{energia-riscalata2}
{\mathcal{E}_{\eps}}(u_{\eps},\Omega):=\int_\Omega\frac{\eps \vert \nabla u_{\eps} \vert^2}{2}
+\frac{(1 - u_\eps^2)^2}{4\eps}, \qquad
\eeq
for every bounded open set $\Omega\subset\R^N$. 

Moreover, since $u \in {\cal F}_p$,  for every $R>0$ and every $ n \ge 1$ we have 

\beq \label{e bound-3} 
{\mathcal{E}_{\eps_n}} (u_{\eps_n},B(0,R)) \le L_u R^{N-1}.
\eeq

Also, from item i) of Theroem \ref{Reg-univBound} we get

\beq \label{bound-L-infinito} 
\forall \, n \ge 1 \qquad \Vert u_{\eps_n} \Vert_{L^{\infty}(\R^N)} \le 1.
\eeq

Note also that $ \Vert V_{{\eps}_n} \Vert $ corresponds to $\frac{\eps_n \vert \nabla u_{\eps_n} \vert^2}{c_{\star}}$. As a consequence of this, of \eqref{energia-riscalata2}, \eqref{e bound-3} and a standard diagonal argument, there  are a subsequence of $(\eps_n)_{n\ge1}$ (not relabeled), a set E of locally finite perimeter in $\R^N$, a Radon measure $ \eta $ on $ \R^N$ and a  ($N$-1)-varifold $V $ on $\R^N$ such that \eqref{convL1}, \eqref{misura-limite} and \eqref{varifold-limite} hold true. 

To prove that $ \eta \neq 0$ on $\R^N$ we observe that 

\beq
{\mathcal{E}_{\eps_n}} (u_{\eps_n},B(0,R)) \longrightarrow \eta(B(0,R)) 
\eeq
where the latter holds true for every $R>0$, with the exception of an at most countable set of values, since $ \eta $ is a Radon measure on $\R^N$ and $\{ \partial B(0,R) \, : \, R>0 \}$ is a disjoint family of Borel sets of $ \R^N$. 

Now, scaling back we get 

\beq\label{cv1f}
{\eps_n}^{N-1} \int_{B(0,\frac{R}{\eps_n})} \frac{\vert \nabla u \vert^2}{2} +\frac{(1 - u^2)^2}{4} = {\mathcal{E}_{\eps_n}} (u_{\eps_n},B(0,R)) \longrightarrow \eta(B(0,R)),
\eeq
for every $R>0$, with the exception of an at most countable set of values.
On the other hand, we also have that  

\beq\label{cv2f}
{\eps_n}^{N-1} \int_{B(0,\frac{R}{\eps_n})} \frac{\vert \nabla u \vert^2}{2} +\frac{(1 - u^2)^2}{4} = R^{N-1} 
\frac{1}{(\frac{R}{\eps_n})^{N-1}} \int_{B(0,\frac{R}{\eps_n})} \frac{\vert \nabla u \vert^2}{2} +\frac{(1 - u^2)^2}{4}
\eeq
for every $R>0$ and every $n \ge 1$. 

By Theorem \ref{Energia}, as $n \to +\infty,$ the right-hand side of the latter converges to $R^{N-1} L_u$ and so, the combination of  \eqref{cv1f} and \eqref{cv2f} yields  

\beq\label{supporto-eta}
0< L_u =  \frac{\eta(B(0,R))}{R^{N-1}} 
\eeq
for every $R>0$, with the exception of an at most countable set of values. In particular the origin $0 \in \R^N$ belongs to $supp  \,\eta$. Note also that \eqref{supporto-eta} holds true for every $ R>0$, since $\eta $ is a non-negative measure (and \eqref{supporto-eta} already holds true for every $R>0$, with the exception of an at most countable set of values). 

Since \eqref{energia-riscalata2}, \eqref{e bound-3} are in force and \eqref{convL1}, \eqref{misura-limite},  \eqref{varifold-limite}  hold true, we can apply Theorem 1 (see also Proposition 4.3 and Proposition 4.4)  of \cite{HT} with $ \lambda_{i} = \lambda_{\infty} =0$.  These results imply that $V$ is a stationary, rectifiable and integral varifold with $ \Vert V \Vert =  \frac{\eta}{c_{\star}}$ (observe that $ \sigma $ in \cite{HT} is exactly $ \frac{c_{\star}}{2}$ in our notations).  The latter also implies that $0 \in \R^N$ belongs to $supp \Vert V \Vert$. Moreover, the combination of \eqref{massaVarifold} and \eqref{supporto-eta}, which holds true for every $R>0$, yields that 

\beq\label{supporto-Varifold}
0< \frac{L_u}{c_{\star} \omega_{N-1}} =  \frac{\Vert V \Vert (B(0,R))}{\omega_{N-1}R^{N-1}} \qquad \forall \quad R>0
\eeq
The latter, together with the fact that $V$ is stationary and integral, implies that $V$ is conical about the origin which, in turn, immediately provides formula \eqref{energiaV}. This concludes the proof. 

\end{proof}

Next we prove that the {\it minimum positive energy level}  is {\it isolated}.  This result  will be a consequence of the following {\it gap phenomenon}  for the set of all finite energy solutions. 
 
\bigskip

\begin{theorem}  {\textbf{(Gap phenomenon for finite energy solutions)}} \label{Gap-finite-energy}

\noindent Assume $N \ge 2$ and let $u$ be a non-constant finite energy solution of \eqref{eqAC}. 
Then 
\beq\label{Valore-gap-finite-energy}
\exists \quad \beta=\beta(N)  > 1
\eeq
such that 
\beq
L_u < \frac{2\sqrt{2}}{3} \omega_{N-1}\beta(N) \quad \Longrightarrow \quad L_u = \ell_{min}.
\eeq
Namely,  the minimum positive energy level is {\rm isolated}.

\medskip

Moreover, any limit stationary, integral and conical varifold $V$ is a ($N$-1)-dimensional hyperplane through the origin $0 \in \R^N$ with multiplicity 1 (and the corresponding set $E$ is one of the two open half-spaces whose boundary is given by $ supp \Vert V \Vert$). 

\smallskip

\end{theorem}

\medskip

\begin{proof}
Let  $(\eps_n)_{n\ge1}$ be a sequence of positive real numbers such that $ \eps_n \to 0$ and let $ \eta $ and $V$ as in Theorem \ref{Teo-finite-energy}. 

By \eqref{misura-limite} and \eqref{massaVarifold}, $ \forall  x \in \R^N$ and  $ \forall  R>0$ we have  

\beq \label{f1}
\Vert V \Vert (B(x,R)) = \frac{\eta}{c_{\star}} (B(x,R)) \le  \frac{1}{c_\star}
\liminf_{n\to+\infty} {\mathcal{E}_{\eps_n}} (u_{\eps_n}, B(x,R)). 
\eeq
On the other hand 
\beq\begin{split}\label{f2}
 {\mathcal{E}_{\eps_n}} (u_{\eps_n},B(x,R))  = \int_{B(x,R)} \frac{\eps_n \vert \nabla u_{\eps_n} \vert^2}{2}
+\frac{(1 - u_{\eps_n}^2)^2}{4\eps_n} & \\
= {\eps_n}^{N-1} \int_{B(\frac{x}{\eps_n},\frac{R}{\eps_n})} \frac{\vert \nabla u \vert^2}{2} +\frac{(1 - u^2)^2}{4} & \\ =
R^{N-1} \frac{1}{(\frac{R}{\eps_n})^{N-1}}  \int_{B(\frac{x}{\eps_n},\frac{R}{\eps_n})} \frac{\vert \nabla u \vert^2}{2} +\frac{(1 - u^2)^2}{4} \le L_u R^{N-1}
\end{split} \eeq
by ii) and iii) of Theorem \ref{Energia}, therefore the combination of \eqref{f1} and \eqref{f2} leads to 
\beq\label{f3}
\forall  x \, \in \R^N, \quad \forall \,  R>0 \qquad \Vert V \Vert (B(x,R)) \le \frac{L_u R^{N-1}}{c_\star} = \Big (\frac{L_u}{\omega_{N-1} c_\star} \Big) \omega_{N-1} R^{N-1}. 
\eeq

\noindent By the {\it Allard Regularity Theorem} \cite{Al1} there exists a small 
$ \delta = \delta(N) >0$ such that 
if $V \neq 0$ is a stationary, integral ($N$-1)-varifold on $\R^N$ with $0 \in supp \Vert V \Vert$
\beq\label{f4}
\forall  x \, \in \R^N, \quad \forall \,  R>0 \qquad \Vert V \Vert (B(x,R)) \le (1 + \delta(N)) \omega_{N-1} R^{N-1}
\eeq
then $V$ is a ($N$-1)-dimensional hyperplane through the origin $0 \in \R^N$ with multiplicity 1. 

The desired conclusions then  follow from \eqref{f3}, by taking $ \beta = 1 + \delta(N) $ in \eqref{Valore-gap-finite-energy}. 

\end{proof}

To proceed further, we recall that in any dimension $ N \ge 2$, there are finite energy solutions with arbitratry large energy (see Appendix \ref {sec : appendix C}  and Appendix \ref {sec : appendix D} ). 
This result and Theorem \ref{Gap-finite-energy} lead to the existence of the {\it second positive energy level.} This is the content of the next  

\bigskip 
 
\begin{theorem}  {\textbf{(Existence of the second positive energy level)}} \label{Teo-esist-secondolivello}

\noindent Assume $N \ge 2$ and set
\beq\label{defLiv2}
{\ell}_{2} : = \inf \{ \, L_u > \ell_{min} \,\, : \,\, u \in {\cal F}_p\, \} 
\eeq
Then, 
\beq\label{Liv2>Liv1}
{\ell}_{min} < {\ell}_{2} \le  2 \ell_{min}.
\eeq

\noindent We shall call this value the  {\rm second positive energy level}. 

%

%
%

\end{theorem}

\medskip

\begin{proof}
By the above discussion the set $\{ \, L_u > \ell_{min} \,\, : \,\, u \in {\cal F}_p \, \}$ is not empty and thus $ {\ell}_{2} $ is finite. Moreover $ {\ell}_{2}  \ge c_{\star} \omega_{N-1}\beta(N) > \ell_{min}$ by Theorem \ref{Gap-finite-energy}, which proves the first inequality in \eqref{Liv2>Liv1}. The other estimate is a consequence of Theorem \ref{Teo-Energia-saddle-2D} and Theorem \ref{energia-arbitr-grande}. 

\end{proof}

\begin{remark}\label{Rk5} As we shall see in Theorem \ref{Teo-valore-secondo-livello} and in Theorem \ref{Secondo-livello4+}, in dimension $N=2,3 $ the second positive energy level is achieved by a member of  ${\cal F}_p$ and its value is $2\ell_{min}$, while for $N=4$ it is achieved by a member of  ${\cal F}_p$ and its value is $\frac{\pi}{2}\ell_{min}$. On the other hand, when $ N \ge 5$, (and also differently from the case of the minimum positive energy level) we do not know if the second positive energy level is achieved by a member of  ${\cal F}_p$. Nevertheless, we shall prove that $\ell_2 \leq \frac{\pi}{2}\ell_{min}$ (see Theorem \ref{Secondo-livello5} for the exact statement). 
\end{remark}

\bigskip

\section{Finite energy solutions in two dimensions : the quantization effect}\label{sec : finite-energy-2D}

\medskip

\begin{theorem}\label{Teo-finite-energy-2D}
Let $u$ be a non-constant solution of the Allen-Cahn \eqref{eqAC} on $\R^2$ with finite energy. Then, 

\item[i)] there is an integer $m = m(u) \ge 1$ such that 

\beq\label{quanti}
L_u = 2m \frac{2\sqrt{2}}{3}  = m \ell_{min}.
\eeq


\item[ii)] For every integer $m \ge 1$ there exists a solution $u$ of the Allen-Cahn \eqref{eqAC} on $\R^2$ 
such that \eqref{quanti} holds true. 


\medskip

In particular, all the energy levels are achieved.

\medskip

\noindent Furthermore, if $V$ is any stationary, integral 1-varifold $V \neq 0$ on $\R^2$ given by Theorem \ref{Teo-finite-energy}, then 

\item[iii)] 

\beq \label{VarCon2D-non}
V = \sum_{j=1}^k \vert v_j\vert \{ \, t v_j : t >0 \} 
\eeq
\noindent where $v_1,...,v_k$, $k\ge2$, are distinct vectors of $\R^2$ of positive integer lenght $ \vert v_j\vert$ and also satisfying the balancing condition : 
\beq \label{BalCond}
\sum_{j=1}^k v_j  = 0.
\eeq
In particular, either $supp \Vert V \Vert $ is a straight line through the origin or $supp \Vert V \Vert $ is a junction, i.e., the union of a finite number of half-lines (at least 3) meeting at the origin.

\item[iv)] $m = 1$ in \eqref{quanti} if and only if $V$ is a straight line through the origin, with multiplicity 1.

\end{theorem}

\begin{proof} Let $E$ and $V$ be as in Theorem \ref{Teo-finite-energy}. We know that $V$ is a stationary, integral 1-varifold and thus, by \cite{AlAlm},  we have that $supp \Vert V \Vert $ is locally either a line segment or a junction point. On the other hand, again thanks to Theorem \ref{Teo-finite-energy}, we also know that $V$ is conical about the origin, therefore we have that
\beq \label{VarCon2D}
V = \sum_{j=1}^k \vert v_j\vert \{ \, t v_j : t >0 \}
\eeq
\noindent where $v_1,...,v_k$, $k\ge2$, are vectors of $\R^2$ of positive integer lenght $ \vert v_j\vert$ and also satisfying the balancing condition (expressing the 1-stationarity of $V$) :
\beq \label{BalCond2}
\sum_{j=1}^k v_j  = 0.
\eeq

\noindent Hence, $supp \Vert V \Vert $ is either a straight line through the origin or a single junction,  i.e.,  either a straight line through the origin or the union of a finite number of half-lines meeting at the origin. Moreover, in the latter case, the number of half-lines is at least 3 (either by \eqref{BalCond} or by the strong maximum principle of Solomon and White \cite{SW}) and the corresponding open half-lines has integer multiplicity. 

In view of \eqref{energiaV}, to prove formula \eqref{quanti} it is enough to demonstrate that $ \omega_1 \Theta(V,0) $ is an even integer. 

If $supp \Vert V \Vert $ is a straight line through the origin, $V$ has constant and integer multiplicity. Thus $ \Theta $ is an integer and we are done (since $ \omega_1 =2)$. So, we can suppose that $supp \Vert V \Vert $ is the union of $k \ge 3$ half-lines meeting at the origin and such that each open half-lines has integer multiplicity ($= \vert v_j\vert $). 

By the integrality of $V$ and the definition of $ \Theta(V,0) $ we see that $ \omega_1 \Theta(V,0) =  \sum_{j=1}^k \vert v_j\vert $, thus $ \omega_1 \Theta(V,0) $ is an integer and we only need to prove that this integer is even. 

To this end we first observe that the open set $ (\R^2 \setminus supp \Vert V \Vert) $ has exactly $k$ connected components, and then that the sequence $(u_{{\epsilon}_n})_{n\ge 1}$ (given by Theorem \ref{Teo-finite-energy}) converges locally uniformly to $1$ (or to $-1$) in each connected component of $(\R^2 \setminus supp \Vert V \Vert) $, thanks to Theorem 1 of \cite{HT}. Let us denote by $ {\cal V}^+$ (resp. $ {\cal V}^-$) the set of the connected components of $ (\R^2 \setminus supp \Vert V \Vert) $ in which $(u_{{\epsilon}_n})_{n\ge 1}$ converges locally uniformly to $1$ (resp. to $-1$). 

We also note that, if $2$ adjacent connected components of $ (\R^2 \setminus supp \Vert V \Vert) $ belong to $ {\cal V}^+$ (resp. $ {\cal V}^-$), the open half-line separating them must belong to  $ (supp \Vert V \Vert \setminus supp \Vert \partial E \Vert) $. Thus, the multiplicity of the open half-line must be even, by item (4) of  Theorem 1 of \cite{HT}. On the other hand, if $2$ adjacent connected components of $ (\R^2 \setminus supp \Vert V \Vert) $ do not belong to the same set $ {\cal V}^+$ (resp. $ {\cal V}^-$), the open half-line separating them must belong to $ supp \Vert \partial E \Vert $ and so its multiplicity must be odd (again by item (4) of  Theorem 1 of \cite{HT}). We shall call such a kind of open half-line a {\it transition line}. 

If one of the two sets $ {\cal V}^+$, $ {\cal V}^-$ is empty, then every open half-line belonging to $supp \Vert V \Vert $ has even multiplicity and we are done. Thus, we can suppose that both $ {\cal V}^+$ and $ {\cal V}^-$ are not empty. In this case there must exist a member $V_1$ of $ {\cal V}^-$ and one of $ {\cal V}^+$,  denoted by $ V_2$, which are adjacent and such that $ V_2$ "follows" $V_1$ when moving counterclockwise around the origin (and in particular, the half-line separating them is a {\it transition line}). We can then label the $k\ge 3$ connected components of $ (\R^2 \setminus supp \Vert V \Vert) $ by $ V_1,....,V_k$ in such a way that $ V_j$ is adjacent to $ V_{j+1}$, for $ j=1, k-1$. Note also that $V_k$ is adjacent to $V_1$. With this definition in hand, we see that the number of  {\it transition lines} is necessarily even. Since the remaining open half-lines have even multiplicity and $ \omega_1 \Theta(V,0) $ is exactly the sum of the multiplicities of the $k$ open half-lines, we have proven the claim. 

\smallskip

Finally, $m=1$ if and only if $  L_u = \ell_{min}$ and the last claim follows from Theorem \ref{Gap-finite-energy}. 
 
\smallskip

item ii) is a consequence of Theorem \ref{Teo-Energia-saddle-2D} and Theorem \ref{Energia1D}.  
 
\end{proof}

As a consequence of the above analysis we have 

\medskip

\begin{corollary}\label{secondo-livello-2D}  Let $u$ be a non-constant finite energy solution of the Allen-Cahn \eqref{eqAC} on $\R^2$ and let $m$ be the integer in \eqref{quanti}.

If $supp \Vert V \Vert $ is a junction of $k$ half-lines, then $k \le 2m $ and $ k=2m$ if and only if each open half-line has multiplicity $1$. 


\end{corollary}
 
\begin{proof} 
Since the multiplicity of each open half-line is at least 1 and $V$ is integral, the conclusions immediately follow from \eqref{energiaV} and \eqref{quanti}. 

\end{proof}

\begin{remark}\label{Rkliv2}
In the next section we shall give a complete description of all the possible limit interfaces $V$, when $ L_u = 2 \ell_{min}$ (see Theorem \ref{Teo-carat-varifold-l2-2D}). 
\end{remark}

\medskip

In view of item iv) of Theorem \ref{Teo-finite-energy-2D} and the above Remark \ref{Rkliv2}, it is worthwhile to investigate all the possible limit configurations $V$ arising in Theorem \ref{Teo-finite-energy-2D} for the higher energy levels. When doing this, in the case in which $V$ is a triple junction (i.e., the union of three half-lines  meeting at the origin) we found the following remarkable relation with {\it Pythagorean triples}. Let us recall that a Pythagorean triple consists of three positive integers $a, b$ and $c$ such that $a^2 + b^2 = c^2. $ Pythagorean triples are infinite and all of them are generated by the {\it scaled Euclid's formula} \cite{Euclide} : $a = k(m^2 - n^2), \, b = k(2mn), \, c = k(m^2 + n^2)$, where $k,m,n$ are positive integers with $m>n$, $m$ and $n$ coprime and not both odd.

Given a Pythagorean triple $(a,b,c)$, we consider the vectors $ v_1 = ( a,b),$ $ v_2 = (-a, b) $ and $ v_3 = (0, -2b). $ The lenght of $v_1$ and $v_2$ is $c$ and the one of $v_3$ is $2b$, which are positive integers. Furthermore, the balancing condition \eqref{BalCond} is clearly satisfied. Thus, $v_1, v_2$ and $ v_3$ provide an admissible limit configuration $V$. 

The simplest example of a Pythagorean triple is given by : $a = 3, b = 4 $ and $c=5$.  Another one is provided by $(65, 72, 97)$, which is the largest of the $16$ primitive Pythagorean triples with $ c \le 100$ (primitive means that $a,b$ and $c$ are coprime, i.e. they do not share any prime factor).

\smallskip

On the other hand, we note that the classical triple junction : $v_1 = (\frac{\sqrt{3}}{2}, \frac{1}{2}) $, $v_2 = (-\frac{\sqrt{3}}{2}, \frac{1}{2}) $ and $ v_3 = (0, -1) $ is not admissible (i.e. it does not exist any entire solution $u$ having this limit configuration). Indeed, $ \vert v_1 \vert =  \vert v_2 \vert =  \vert v_3 \vert =1$ and so the corresponding energy would be $ L_u = \Theta(V,0) \ell_{min} = \frac{3}{2} \ell_{min}$, which contradicts the quantization formula \eqref{quanti}.  

\medskip

In view of the above discussions, the following intriguing and interesting question arises in a natural way : 

\medskip

\noindent {\bf  Open question : } {\it Do there exist entire solutions of the Allen-Cahn equation having a Pythagorean conical limit interface $V$ ? }

\medskip

\section{The second energy level}\label{sec : second-level}

\medskip

In Theorem \ref{Teo-esist-secondolivello} we proved the existence of the {\it second positive energy level}, denoted by $\ell_2$. We also proved that $ \ell_2 \le 2 \ell_{min}$. In this section we analyse the general situation  and we prove that, in dimension $N=2$ and $ N=3$,  the {\it second positive energy level} coincides exactly with the value $ 2 \ell_{min}$ while, for $ N=4 $, we have that ${\ell}_{2} = \frac{\pi}{2} \ell_{min} < 2 \ell_{min}.$ We also demonstrate that the  {\it second positive energy level} $\ell_2$ is achieved by a member of  ${\cal F}_p$ and we give a complete description of all the possible limit interfaces $V$. We underline the fact that the situation is very different depending on the value $ N=2,3$ or $N=4$. This  will require a particular (and different) treatement for each of thoses cases. 

In the second part of the present section we shall consider the case $ N \ge5 $ and we shall provide an explicit upper bound for $\ell_2$.  It will depend on Simons' cones. See Theorem \ref{Secondo-livello5} below. 

\medskip

More precisely we have 

\medskip
 
\begin{theorem}\label{Teo-valore-secondo-livello}
Assume $N = 2$ or $N=3$. We have

\beq\label{Liv2in2D} 
{\ell}_{2} = 2 \ell_{min}
\eeq

\noindent and the {\rm second positive energy level} is achieved by a member of  
${\cal F}_p.$

\end{theorem}

\medskip
 
\noindent and 

\medskip 
 
\begin{theorem}\label{Secondo-livello4+}
Assume $N = 4$. We have

\beq\label{Liv2in4D} 
{\ell}_{2} = \frac{\pi}{2} \ell_{min} 
\eeq

\noindent and the {\rm second positive energy level} is achieved by a member of  ${\cal F}_p.$ 

\end{theorem}

\bigskip

 
\begin{theorem}\label{Secondo-livello5}
Assume $N \ge 5$. We have
\begin{align}\label{Liv2in5D+} 
{\ell}_{2}(N) \leq \Theta(\mathcal{C}_{n,n},0) \ell_{min}(N)
\end{align}
where $n = \lfloor \frac{N}{2}\rfloor$  and $\Theta(\mathcal{C}_{n,n},0)$ is the density of the Simons cone $\mathcal{C}_{n,n}$ at the origin (see Appendix \ref{sec : appendix D} and \ref{sec : appendix F}). 
\end{theorem}

\begin{remark} We can compute the exact value of $\Theta(\mathcal{C}_{n,n},0)$ by using Proposition \ref{Teo-dens-coni-Simons}. For instance, for $n=2,3,4$ and $5$ we get 
$\frac{\pi}{2}, \, \frac{3}{2}, \, \frac{15}{32}\pi$ and $ \frac{35}{24}$, therefore we have
\begin{align*}\label{Liv2in5D+} 
&{\ell}_{2}(5) \leq \frac{\pi}{2} \ell_{min}(5) \\
&{\ell}_{2}(N) \leq \frac{3}{2} \ell_{min}(N) \qquad \text{if} \quad N \geq 6 \\
&{\ell}_{2}(N) \leq \frac{15}{32}\pi \ell_{min}(N) \qquad \text{if} \quad N \geq 8 \\
&{\ell}_{2}(N) \leq \frac{35}{24} \ell_{min}(N) \qquad \text{if} \quad N \geq 10\\
\end{align*}
Note that $\frac{\pi}{2} > \frac{3}{2} > \frac{15}{32}\pi > \frac{35}{24}$.

\end{remark}

\medskip

Let us start with the following general result
 
\medskip

\begin{theorem}\label{Teo-sotto-2lmin} Assume $ N \ge 3$ and let $u$ be a solution of the Allen-Cahn \eqref{eqAC} on $\R^N$ satisfying 
\beq \label{sotto-2lmin}
0 < L_u < 2 \ell_{min}.
\eeq
For every limit pair $(V,E)$ obtained in Theorem \ref{Teo-finite-energy} we have 

\medskip

\item i)  $V$ is the multiplicity-one varifold associated to the boundary of a Borel set, of locally finite perimeter in $\R^N$, equivalent to $E$.

\medskip

\item ii) If $ N = 3$, then $V$ is a plane with multiplicity 1 and $L_u = \ell_{min}.$

\medskip

\item iii) If $ N = 4$, then $ V$ is smooth outside the origin and either $ L_u = \ell_{min}$ or $L_u \ge \frac{\pi}{2} \ell_{min} . $
\end{theorem}

\medskip

\begin{proof}

Since $V$ is stationary and conical about the origin, 
formula \eqref{energiaV} and the assumption $ 0 < L_u < 2 \ell_{min}$ imply that 
\beq \label{dens0X}
0 < \Theta(V,x) \le \Theta(V,0) < 2
\eeq
for every $x \in supp \Vert V \Vert$. 

Thus, from the integrality of $V$, we infer that 

\beq \label{densXae}
\Theta(V,x)  =1 
\eeq
for ${\cal H}^{N-1}$ a.e. $x \in supp \Vert V \Vert$.

Therefore, by item (4) of  Theorem 1 of \cite{HT}, we have that $ \Vert V \Vert = \Vert \partial E \Vert$,  where  $\Vert \partial E \Vert $ denotes the perimeter of $E$ in the measure-theoretic sense. In particular,
$ supp \Vert V \Vert = supp \Vert \partial E \Vert $, and thus the varifold $V$ is the boundary of a Borel set, of locally finite perimeter in $\R^N$, equivalent to $E$. 
This proves item i). It follows that, for every $ x \in supp \Vert V \Vert $, every tangent cone $C$ for $V$ at $x$ is stationary, integral and $C$ is the boundary of a Borel set $F$ of locally finite perimeter in $\R^N$. This follows, for example, from White's compatibility results (see Theorem 1.2 and Theorem 3.6 of \cite{White2}) and that the tangent cone $C$ has unit multiplicity ${\cal H}^{N-1}$ a.e. (thanks to $ \Theta (C,0) = \Theta (V,x) <2$ by \eqref{dens0X}). Furthermore, by a classical slicing argument, if the tangent cone $C$ splits as $C_k \times \R^k$, then the cone $C_k$ is the boundary of a Borel set, of locally finite perimeter in $\R^{N-k}$.

In order to prove items ii) and iii), we recall that (see for instance \cite{Alm2, White}) tangent cones provide a stratification of a stationary varifold $V$. More precisely, for $j=0,1,\ldots, N-1$, consider the set  
\beq \label{strati}
{\Sigma_j} : = \{ \, x \in supp \Vert V \Vert \, :  \, dim(Spine(C)) \le j \quad  
{\text {for all tangent cones $C$ for $V$ at}} \, x  \, \},  
\eeq
where $Spine(C)$ is the vector subspace of $\R^N$ defined by 

\beq \label{spine}
Spine(C) : = \{ x \in supp \Vert C \Vert \, : \, \Theta(C,x) = \Theta(C,0) \}.
\eeq

We clearly have that 

\beq \label{incl-strati}
\Sigma_0 \subset \Sigma_1 \subset .... \subset \Sigma_{N-1} = supp \Vert V \Vert 
\eeq
furthermore (see \cite{Alm2, White}) the Hausdorff dimension of $\Sigma_j$ is at most $j$ (actually $\Sigma_0$
is at most countable).

\medskip

Let us turn to the 3-dimensional case (i.e. $N=3$). 

\medskip

We claim that : $  \Sigma_2 \setminus \Sigma_1 \subset reg(V)$. Let $ x \in  \Sigma_2 \setminus \Sigma_1$ then there  exists a tangent cone $C$ for $V$ at $x$ such that $dim(Spine(C)) =2$. Hence, $C$ must be a two-dimensional plane with integer multiplicity equal to $ \Theta(C,0) =  \Theta(V,x)$ and thus $  \Theta(V,x) =1 $  by \eqref{dens0X}. This proves that $ x \in reg(V)$ by Allard's regularity Theorem. 

We now claim that : $  \Sigma_1 \setminus \Sigma_0 \subset reg(V)$. Let $ x \in \Sigma_1 \setminus \Sigma_0$, then there exists a tangent cone $C$ for $V$ at $x$ such that $dim(Spine(C)) =1$.  Hence, up to a rotation, $ C = C_0 \times \R $, for some 1-stationary integral cone $C_0$ of $ \R^2$. Since a 1-stationary cone is a finite union of rays, 
either $C_0$ is a straight line and so, as before, $C$ must be a two-dimensional plane with integer multiplicity equal to $ \Theta(C,0) =  \Theta(V,x)$ and thus $ \Theta(V,x) =1 $ by \eqref{dens0X}, or $C_0$ is the union of exactly 3 half-lines of multiplicity $1$. The  latter case is impossible since $C_0$ is the boundary of a Borel set $F$ of locally finite perimeter in $\R^2$. 
Therefore we have $ x \in reg(V)$ by Allard's regularity Theorem.

So far we have proved that $supp \Vert V \Vert \setminus \Sigma_0 = \Sigma_2 \setminus \Sigma_0 \subset reg(V)$. Next we prove that $ Sing(V) \subset \{0\}$. If not, by what we have just proven, we would have 
$ \Sigma_0  \cap (Sing(V) \setminus \{0\}) \neq \emptyset $. But $V$ is a cone, and so $ \Sigma_0$ would contain an open ray of singular points, contradicting the fact that the Hausdorff dimension of $ \Sigma_0$ is at most $0$. 

By summarizing, we proved that  $V$ is smooth outside the origin. From this, we see that $V $ is a cone over a smooth stationary 1-varifold $ \gamma $ of $\esse^2$, the unit sphere of $ \R^3$. But then we know from \cite{AlAlm} that $\gamma$ is a geodesic of ${\esse}^2$ and so $ V$ is a plane with multiplicity 1 by \eqref{dens0X}. This and formula \eqref{energiaV}  conclude the proof of item ii). 

\medskip

The proof of item iii) is more involved. Again, we analyse the different strata of $V$. As before, $  \Sigma_3 \setminus \Sigma_2 \subset reg(V)$. Indeed, if $ x \in  \Sigma_3 \setminus \Sigma_2$ then there  exists a tangent cone $C$ for $V$ at $x$ such that $dim(Spine(C)) =3$. Hence, $C$ must be a 3-dimensional plane with integer multiplicity equal to $ \Theta(C,0) =  \Theta(V,x)$ and thus $  \Theta(V,x) =1 $  by \eqref{dens0X}. So $ x \in reg(V)$ by Allard's regularity Theorem. 

Then we prove that $ \Sigma_2 \setminus \Sigma_1 \subset reg(V)$. Let $ x \in \Sigma_2 \setminus \Sigma_1$, then there exists a tangent cone $C$ for $V$ at $x$ such that $dim(Spine(C)) =2$.  Hence, up to a rotation, $ C = C_0 \times \R^2 $, for some 1-stationary integral cone $C_0$ of $ \R^2$, which is the boundary of a Borel set of locally finite perimeter.  At this point, the argument used to exclude the $C_0$ is the union of exactly 3 half-lines of multiplicity $1$ in the 3-dimensional case applies verbatim, and so we get the conclusion. 

Now we prove that also $ \Sigma_1 \setminus \Sigma_0 \subset reg(V)$. This require some new ingredients. Let $ x \in \Sigma_1 \setminus \Sigma_0$, then there exists a tangent cone $C$ for $V$ at $x$ such that $dim(Spine(C)) =1$.  Hence, up to a rotation, $ C = C' \times \R $, for some 2-stationary integral cone $C'$ of $ \R^3$, which is the boundary of a Borel set of locally finite perimeter. 
Then we observe that $ \Theta(C',0') = \Theta(C,0) =  \Theta(V,x) < 2$  by \eqref{dens0X}, so that we can apply (verbatim) to the cone $C'$ the proof of item ii) to deduce that $C'$ is a plane with multiplicity 1, proving thus the claim.

So $supp \Vert V \Vert \setminus \Sigma_0 = \Sigma_3 \setminus \Sigma_0 \subset reg(V)$ and $\Sigma_0 \subset Sing(V)$. On the other hand $ Sing(V) \subset \{0\}$, since otherwise we would have $\Sigma_0  \cap (Sing(V) \setminus \{0\}) \neq \emptyset $. But $V$ is a cone, and so $ \Sigma_0$ would contain an open ray of singular points, contradicting the fact that the Hausdorff dimension of $ \Sigma_0$ is at most $0$. Hence, $V$ is smooth outside the origin. 

To demonstrate the second part of item iii) we proceed as follows. 
Since $V$ is smooth outside the origin, $V$ is a cone over a smooth stationary 2-varifold (a embedded closed smooth minimal surface) $ S $ of $\esse^3$, the unit sphere of $\R^4$. This link $S$ must then be connected by Frankel's theorem. Now, if the genus of $S$ is zero, we know from a result of Almgren \cite{Alm1} (see also \cite{Lawson-S3}) that $S$ is totally geodesic and so its area is $4\pi$.  This yields $ \Theta(V,0) = 1$, so $V$ is a plane with multiplicity 1 by \eqref{dens0X} and $ L_u = {\ell}_{min}$. If the genus of $S$ is positive the deep Theorem B of F.C. Marques and A. Neves \cite{Coda-Neves} tell us that the area($S$) $\ge 2\pi^2$. Thus, $ \Theta(V,0) \ge  \frac{\pi}{2}$, which completes the proof of the Theorem. \end{proof}

\smallskip

\begin{proof}[Proof of Theorem \ref{Teo-valore-secondo-livello} ]
\noindent From Theorem \ref{Teo-esist-secondolivello} we already know that $ \ell_2 \le 2 \ell_{min}$. Also, from Theorem \ref{Teo-Energia-saddle-2D} and Theorem \ref{energia-arbitr-grande} we know that there are entire solutions $u$ of the Allen-Cahn equation with $L_u = 2 \ell_{min}$. To get the desired conclusion is then enough to prove that $ \ell_2 \ge 2 \ell_{min}$.

For $ N=2$ the claim is contained in Theorem \ref{Teo-finite-energy-2D} or in Theorem \ref{Teo-Energia-saddle-2D}. 
This concludes the proof for the two-dimensional case. 

Assume $N=3$. If $ L_u < 2 \ell_{min} $,  item ii) of Theorem \ref{Teo-sotto-2lmin} would give $L_u = \ell_{min}$ a contradiction with the definition of the second energy level. \end{proof}

We can now proceed to the 

\smallskip

\begin{proof} [Proof of Theorem \ref{Secondo-livello4+}]  
From Corollary \ref{CorCabré} we know that there are entire solutions $u$ of the Allen-Cahn equation with $L_u = \frac{\pi}{2} \ell_{min} .$ The latter and item iii) of Theorem \ref{Teo-sotto-2lmin} prove the theorem. 
\end{proof}

\smallskip

\begin{proof} [Proof of Theorem \ref{Secondo-livello5}]  
From Corollary \ref{CorCabré} we know that there are entire solutions $u$ of the Allen-Cahn equation with $L_u = \frac{\pi}{2} \ell_{min}$ for any $ N \geq 4$. This proves the theorem for $N=5$. 
More generally, for any $ m \geq 2$, in view of formula \eqref{val-energ-sad}, 
the saddle-type solution $u=u_m(x,y)$ built in \cite{CT1} for $(x,y) \in \R^m \times \R^m = \R^{2m}$, is an entire solution of \eqref{eqAC} satisfying $L_u(2m) = \Theta(\mathcal{C}_{m,m},0) \ell_{min}(2m)$. Hence, $ \ell_2(2m) \leq \Theta(\mathcal{C}_{m,m},0) \ell_{min}(2m)$. To proceed further, we take $n = \lfloor \frac{N}{2}\rfloor$ and we consider the cylindrical solution $u(x_1,\ldots,x_N) := u_n(x_1,\ldots,x_{2n})$ of \eqref{eqAC} in $\R^N$, where $u_n$ is the $2n-$dimensional saddle-type solution just considered. The desired conclusion then follows by applying \eqref{energiaV} and Proposition \ref{Sim2} about cylindrical solutions. 
\end{proof}

\medskip

Next we carachterize all the possible limit varifolds $V$ arising when $u$ is a solution of \eqref{eqAC} achieving the second level of energy $ \ell_2$, for $ N=2,3,4$. 

\smallskip

Let us first consider the case $ N=2$. 
 
\medskip 
 
\begin{theorem}\label{Teo-carat-varifold-l2-2D}
Assume $N = 2$ and let $u$ be a solution of \eqref{eqAC} such that 
\beq
L_u = {\ell}_{2} = 2 \ell_{min}. 
\eeq

Then, either $V$ is a straight line through the origin with multiplicity 2, or $V$ is the union of exactly two different straight lines through the origin, each of which of multiplicity 1.

\end{theorem}

\begin{proof}

From Theorem \ref{Teo-finite-energy-2D} and Corollary \ref{secondo-livello-2D} we know that either $supp \Vert V \Vert $ is a straight line through the origin, and necessarily of multiplicity 2, or $supp \Vert V \Vert $ is a junction of $ 3 \le k \le 4$ half-lines. 

Let us first prove that the triple junction cannot occur. For contradiction, suppose that $ k=3$ then, by Theorem \ref{Teo-finite-energy-2D}, there are distinct vectors $ v_1, v_2, v_3$ satisfying $ v_1 + v_2 + v_3 = 0 $ and $ \vert v_1 \vert =  \vert v_2 \vert = 1$, $  \vert v_3 \vert =2$. Hence, $ 4 = \vert v_3 \vert^2 =  \vert v_1 + v_2 \vert^2 = 1 + 1 + 2 (v_1, v_2) $, which in turn yields $ (v_1, v_2) =1$. The latter clearly implies $ v_1 = v_2$, a contradiction. 

We are left with the case in which $supp \Vert V \Vert $ is a junction of 4 half-lines. In this case,  from Theorem \ref{Teo-finite-energy-2D} and Corollary \ref{secondo-livello-2D} we know that there are distinct {\it unit} vectors $ v_1, v_2, v_3, v_4$ satisfying the balancing condition $ v_1 + v_2 + v_3 +v_4 = 0 $. We may (and do) suppose that $v_1 = (1,0) $ and that $ v_1, v_2, v_3, v_4$ are ordered (as points on the unit circle when moving counterclockwise).  Since 
$\vert v_1 +v_2 \vert^2 =  \vert v_3 + v_4 \vert^2 $ we get $ 1 + 1 + 2 (v_1, v_2) = 1 + 1 + 2 (v_3, v_4) $ and so 
$ (v_1, v_2) = (v_3,v_4)$. Exchanging the role of $v_2$ and $v_4$ we also obtain $ (v_1, v_4) = (v_3,v_2)$.  These relations imply that the sum of the angle between $v_1$ and $ v_2$ and the angle between $v_2$ and $ v_3$ is $\pi$. This means that $ v_3 = - v_1$ and so $ v_2 = - v_4$ as desired. 

\end{proof}

\medskip

\begin{theorem}\label{Teo-carat-varifold-l2-3D}
Assume $N = 3$ and let $u$ be a solution of \eqref{eqAC} such that 
\beq
L_u = {\ell}_{2} = 2 \ell_{min}. 
\eeq

Then, either $V$ is a plane through the origin with multiplicity 2, or $V$ is the union of exactly 2 distinct, multiplicity one planes through the origin, meeting along a unique common straight line. 

\end{theorem}

\bigskip

\begin{proof}

Since $V$ is stationary and conical about the origin, 
formula \eqref{energiaV} and the assumption $ L_u = 2 \ell_{min}$ imply that 
\beq \label{dens0X2}
0 < \Theta(V,x) \le \Theta(V,0) = 2
\eeq
for every $x \in supp \Vert V \Vert$. 

At this point there are two possibilities : either $ \Theta(V,x) < \Theta(V,0) = 2$ for every $ x \in supp \Vert V \Vert \setminus \{0\}$, or there is $x \in supp \Vert V \Vert \setminus \{0\}$ such that $ \Theta(V,x) = \Theta(V,0) = 2.$

The first case is impossible. Indeed, from the integrality of $V$, we infer that $\Theta(V,x)  =1 $ for ${\cal H}^{2}$ a.e. $x \in supp \Vert V \Vert$ and then we can proceed exactly as is the proof of Theorem \ref{Teo-sotto-2lmin} to get that  $V$ is a plane with multiplicity 1, a contradiction.

In the second case, the $Spine (V)$ is necessarily  a vector space of dimension 1 or 2. If $dim(Spine(V)) =2$, then $V$ is a plane through the origin with multiplicity 2. If  $dim(Spine(V)) =1$ then, up to a rotation, $ V = V_0 \times \R $, for some 1-stationary integral cone $V_0$ of $ \R^2$. Since $L_u = {\ell}_{2} = 2 \ell_{min}$, the geometric classification in the proof of Theorem \ref{Teo-carat-varifold-l2-2D} tell us that $V_0$ must be the union of exactly two different straight lines through the origin, each of which of multiplicity 1. The desired conclusion then follows by recalling that $ V = V_0 \times \R $. 
\end{proof}

\medskip

Finally, we turn to

\medskip
 
\begin{theorem}\label{Teo-carat-varifold-l2-4D}
Assume $N = 4$ and let $u$ be a solution of \eqref{eqAC} such that 
\beq
L_u = {\ell}_{2} = \frac{\pi}{2}\ell_{min}. 
\eeq
Then, up to rotations in $\R^4$, $V$ is the Simons' cone $ \, \, {\cal C}_{2,2}: = \{ (x,y) \in \R^2 \times \R^2 \,\,  : \,\, \vert x \vert^2 = \vert y \vert^2 \}$ (with multiplicty 1). Or, equivalently, $V$ is the cone over the Clifford torus $ S^{1}_{{\frac{1}{\sqrt{2}}}} \times S^{1}_{{\frac{1}{\sqrt{2}}}} \subset \esse^3$, up to isometries of $\esse^3$.

\end{theorem}

\begin{proof} Recall that, from the end of the proof of Theorem \ref{Teo-sotto-2lmin} we already know that $V $ is a cone over a smooth stationary 2-varifold (an embedded closed smooth minimal surface) $ S $ of $\esse^3$ with positive genus.  Another application of the deep Theorem B of F.C. Marques and A. Neves \cite{Coda-Neves}  provides  the desired conclusion.\end{proof}

\medskip

\begin{remark} We stress on the fact the the saddle-type solution \cite{CT1,CT2,CT3} built in $\R^4$ and analyzed in detail in Theorem \ref{energia-saddle-Cabré} has exactly  the Simons' cone $ \, \, {\cal C}_{2,2}$ as (unique) limit varifold $V$ (see  Theorem \ref{energia-saddle-Cabré}). 
\end{remark}

\medskip

\section{Stable solutions} \label{sec: Stabili}

\medskip

A function $u \in \mathcal{C}^2({\R}^N)$ which solves $-\Delta u = u-u^3$ on $\R^N$ is said to be {\it stable} if 

\beq\label{stable}
\forall \psi \in \mathcal{C}^1_c({\R}^N) \qquad {\cal Q}_u (\psi) : = \int_{\R^N}\vert \nabla \psi \vert^2 - (1-3u^2) \psi^2 \ge 0,  
\eeq
which is equivalent to say that $u$ is a {\it stable critical point} of the energy functional $ {\mathcal{E}}(u,\Omega) :=\int_{\Omega} \frac{|\nabla u|^2}{2} +  \frac{(1-u^2)^2}{4}$ for every bounded open set $\Omega\subset\R^N$. 

\medskip

Typical examples of stable solutions are local minimizers and monotone solutions. 

\bigskip 

Let us recall the following well-known result about the classification of stable solutions in dimension $ N=2$ and $N=3$. Its proof is part of the folklore and it can be found and/or infered from \cite{GG, AC, BCN2, FarHab, Far-Dino-Enrico-SNS}. In \cite{FarHab,Far-Dino-Enrico-SNS}, a proof is given under weaker assumptions and for general semilinear equations of the form $ - \Delta u = f(u) $, with $ f \in C^{0,1}_{loc}(\R)$, as well as for quasilinear equations in divergence form. 

\bigskip

\begin{theorem}\label{Teo-stabili-folklore}
Assume $ N=2$ or $N=3$ and let $u$ be a stable solution of \eqref{eqAC} satisfying 
\beq\label{cond-en-cinetica} 
\exists \, x_0 \in \R^N, \quad \exists R_0 >0 \quad : \quad \forall R >R_0 \qquad \quad \int_{B(x_0,R)} \vert \nabla u \vert^2 \le C R^2,
\eeq
for some positive constant $ C$ independent on $R$.

Then $u$ is 1D, that is, either $u$ is constant (necessarily equal to $ \pm1$) or $u(x) = \pm \tanh\big (\frac{(\nu\cdot x) +\alpha}{\sqrt 2} \big)$, for some unit vector $ \nu \in \R^N$ and some real  number $\alpha$. 
\end{theorem}

\bigskip

The main result of this section concerns the classification of stable solutions to \eqref{eqAC} in any dimension $ N\ge 2$. We give a detailed description of the limit interface, the conical varifold $V$ and of its regularity. In particular we prove that, for $ N\le 7$, the varifold $V$ is always a ($N$-1)-dimensional hyperplane through the origin $0 \in \R^N$.
This result seems to corroborate the conjecture about the classification of stable solutions to the Allen-Cahn equation raised on p. 1135 of \cite{PacWei-stab} by F. Pacard and J. Wei.

\bigskip
 
\begin{theorem} {\textbf{(Structure Theorem for stable solutions)}} \label{Teo-stabili}

\noindent Assume $N \ge 2$. Let $u$ be a non-constant, stable, finite energy solution of \eqref{eqAC} and let $(\eps_n)_{n\ge1}$ be a sequence of positive real numbers such that $ \eps_n \to 0$. 

Then there are a subsequence of $(\eps_n)_{n\ge1}$ (not relabeled), a set E of locally finite perimeter in $\R^N$ and a stationary, integral ($N$-1)-varifold $V \neq 0$ on $\R^N$ such that 

\smallskip

\item[i)]  

\beq\label{convL1bis} 
u_{\eps_n} \longrightarrow \, \,  {\bf 1}_E - {\bf 1}_{\R^N \setminus E} \quad  {\text {a.e. \, and \, in}} \quad L^1_{loc}(\R^N),
\eeq

\beq\label{varifold-limite2} 
V_{\eps_n} \longrightarrow \, \,  V \quad  {\text {in the varifolds sense}}. 
\eeq

\smallskip

\item[ii)] 

\beq\label{Vconica-sta}
0 \in supp \Vert V \Vert, \qquad  {\text {$V$ is conical about $0$}},
\eeq

\beq
L_u = \frac{2\sqrt{2}}{3} \omega_{N-1} \Theta(V,0),
\eeq
\noindent  where $ \Theta(V,0) $ denotes the density function of the varifold $V$  at the origin $0$. 

\smallskip

\item[iii)] For $ N \le 7$, $supp \Vert V \Vert$ is a ($N$-1)-dimensional hyperplane through the origin $0 \in \R^N$.

\smallskip

\item[iv)] For $ N = 8 $, $ Sing(V) \subseteq \{ 0 \}$.

\smallskip

\item[v)] For $ N \ge 9 $, the set $ Sing(V) $ has Hausdorff dimension at most $ N-8$. 

\smallskip

\item[vi)] For $ N \ge 8 $, the set $ Reg(V) $ is an embedded smooth stable minimal hypersurface of $ \R^N$. 

\smallskip

\noindent Here, as customary, $Reg(V)$ is the set of regular points of  $supp \Vert V \Vert$ while $Sing(V) := supp \Vert V \Vert \setminus Reg(V)$ denotes the singular set of $V$.

\end{theorem}

\bigskip

\begin{remark}\label{Rk6}
It is interesting to note the strong analogy existing between the case of stable solutions described by the above Theorem \ref{Teo-stabili} and the one of local minimizers described by means of Theorem \ref{Teo-minloc}. 
\end{remark}

\bigskip

\begin{proof}[Proof of Theorem \ref{Teo-stabili}]
Let  $(\eps_n)_{n\ge1}$ be a sequence as in the statement and let $E$ and $V$ be as in Theorem \ref{Teo-finite-energy}. 

Items i) and ii) follow immediately from Theorem \ref{Teo-finite-energy} (since $u$ has finite energy). 
To prove items iii)-vi) we invoke Theorem 2.1. of \cite{TW} (see also \cite{Tone-stable} for $ N=2$) which describes the structure of any limit-interface corresponding to a locally uniformly bounded, locally energy-bounded sequence of stable critical points of the 
Van der Waals-Cahn-Hilliard energy functionals. 

Here we are going to use in a crucial way the notion of stability of $u$, besides the finiteness of its energy, as well as the fact that the limit varifold $V$ is conical about the origin.  

Since $u$ is an entire stable solution of the Allen-Cahn equation, the rescaled functions $u_{\eps} $ are stable solutions to \eqref{eq-riscalata}. In particular, they are 

\beq\begin{split}\label{energia-riscalata3}
{\text {\it stable critical points of the rescaled energy functional }}  \\\\
{\mathcal{E}_{\eps}}(u_{\eps},\Omega):=\int_\Omega\frac{\eps \vert \nabla u_{\eps} \vert^2}{2}
+\frac{(1 - u_\eps^2)^2}{4\eps}, \qquad
\end{split}\eeq
for every bounded open set $\Omega\subset\R^N$. 

Moreover, since $u \in {\cal F}_p$,  for every $R>0$ and every $ n \ge 1$ we have 

\beq \label{e bound-4} 
{\mathcal{E}_{\eps_n}} (u_{\eps_n},B(0,R)) \le L_u R^{N-1}.
\eeq

Also, from item i) of Theroem \ref{Reg-univBound} we get

\beq \label{bound-L-infinito-sta} 
\forall \, n \ge 1 \qquad \Vert u_{\eps_n} \Vert_{L^{\infty}(\R^N)} \le 1.
\eeq

Now, in view of \eqref{energia-riscalata3}-\eqref{bound-L-infinito-sta} we can apply Theorem 2.1. of \cite{TW} to obtain that the limit varifold $V$ (which is not necessarily connected) is such that $ Sing(V) = \emptyset $ if $ N \le 7$, is a discrete set of points if $ N=8$ and has Hausdorff dimension at most $ N-8$ when $ N\ge 9$. Furthermore, $ Reg(V) $ is an embedded smooth stable minimal hypersurface of $ \R^N$.  The latter immediately gives items v) and vi). To prove the remaining items we make use of the crucial new fact (proved before)  that $V$ is also conical. 

Indeed, since $V$ is conical about the origin and $ Sing(V) = \emptyset $ if $ N \le 7$, then $supp \Vert V \Vert$ must be a ($N$-1)-dimensional hyperplane through the origin. 
When $ N=8$ we know that $ Sing(V) $ is a discrete set of points, therefore $ Sing(V) \subseteq \{ 0 \}$  (since $V$ is conical about the origin). 

\end{proof}

\section{Solutions stable outside a compact set} \label{sec: Stabili fuori}

\medskip

In this section we consider solutions to the Allen-Cahn equation \eqref{eqAC} which are {\it stable outside a compact set.} The notion of stability outside a compact set was introduced in \cite{Fa4} (cf. Definition 2.1. therein. See also \cite{Fa2}) as a natural extention of the classical notion of stability considered in Section \ref{sec: Stabili}. It asks stability only "near infinity" and thus it encompasses the important family of solutions with Finite Morse Index. More precisely, a solution $u$ of the 
Allen-Cahn equation \eqref{eqAC} is  {\it stable outside a compact set} ${\cal K}$ of $\R^N$,  if

\beq\label{def-stable-fuori}
\forall \, \psi \in \mathcal{C}^1_c({\R}^N \setminus {\cal K}) \qquad{\cal Q}_u (\psi) : = \int_{\R^N}\vert \nabla \psi \vert^2 - (1-3u^2) \psi^2 \ge 0.
\eeq

In the sequel we shall distinguish between the case $ N=2$ and $ N\ge3$. Surprisingly, (see Remark \ref{Rk7}) the situation in the two-dimensional case is "richer" than the one in dimension $ 3 \le N \le 7$.  Let us start with the case $N \ge 3$.

\medskip

\begin{theorem} {\textbf{(Structure Theorem for solutions stable outside a compact set)}} \label{Teo-stabili-fuori-alta-dim}

\noindent Assume $N \ge 3$ and let $u$ be a non-constant solution of \eqref{eqAC} with finite energy. 
Assume that $u$ is stable outside a compact set and let $(\eps_n)_{n\ge1}$ be a sequence of positive real numbers such that $ \eps_n \to 0$. 

Then there are a subsequence of $(\eps_n)_{n\ge1}$ (not relabeled), a set E of locally finite perimeter in $\R^N$ and a stationary, integral ($N$-1)-varifold $V \neq 0$ on $\R^N$ such that 

\smallskip

\item[i)]  

\beq\label{convL1tris} 
u_{\eps_n} \longrightarrow \, \,  {\bf 1}_E - {\bf 1}_{\R^N \setminus E} \quad  {\text {a.e. \, and \, in}} \quad L^1_{loc}(\R^N),
\eeq

\beq\label{varifold-limite-fuori} 
V_{\eps_n} \longrightarrow \, \,  V \quad  {\text {in the varifold sense}}. 
\eeq

\smallskip

\item[ii)] 

\beq\label{Vconica-stab-fuori}
0 \in supp \Vert V \Vert, \qquad  {\text {$V$ is conical about $0$}},
\eeq

\beq
L_u = \frac{2\sqrt{2}}{3} \omega_{N-1} \Theta(V,0),
\eeq
\noindent  where $ \Theta(V,0) $ denotes the density function of the varifold $V$  at the origin $0$. 

\smallskip

\item[iii)] For $ 3 \le N \le 7$, $supp \Vert V \Vert$ is a ($N$-1)-dimensional hyperplane through the origin $0 \in \R^N$.

\smallskip

\item[iv)] For $ N = 8 $, $ Sing(V) \subseteq \{0 \}$.

\smallskip

\item[v)] For $ N \ge 9 $, the set $ Sing(V) $ has Hausdorff dimension at most $ N-8$. 

\smallskip

\item[vi)] For $ N \ge 8 $, the set $ Reg(V) $ is an embedded smooth stable minimal hypersurface of $ \R^N$. 

\smallskip

\noindent Here, as customary, $Reg(V)$ is the set of regular points of  $supp \Vert V \Vert$ while $Sing(V) := supp \Vert V \Vert \setminus Reg(V)$ denotes the singular set of $V$. 

\end{theorem}

\bigskip

\begin{proof}[Proof of Theorem \ref{Teo-stabili-fuori-alta-dim}]
Let  $(\eps_n)_{n\ge1}$ be a sequence as in the statement and let $E$ and $V$ be as in Theorem \ref{Teo-finite-energy}. 

Since $u$ is a solution of the Allen-Cahn equation and $u$ is stable outside a large ball centered at the origin, we see that, for every $ a>0$  there is $ \eps(a) >0$ such that for every $ \eps_n < \eps(a) $, the rescaled functions $u_{\eps_n} $ are 

\beq\begin{split}\label{energia-riscalata3-fuori}
{\text {\it stable critical points of the rescaled energy functional }}  \\\\
{\mathcal{E}_{\eps}}(u_{\eps},\Omega):=\int_\Omega\frac{\eps \vert \nabla u_{\eps} \vert^2}{2}
+\frac{(1 - u_\eps^2)^2}{4\eps}, \qquad
\end{split}\eeq
for every bounded open set $\Omega\subset \R^N \setminus {\overline {B(0,a)} }$. 

Moreover, since $u \in {\cal F}_p$,  for every $R>0$ and every $ n \ge 1$ we have 

\beq \label{e bound-4-fuori} 
{\mathcal{E}_{\eps_n}} (u_{\eps_n},B(0,R)) \le L_u R^{N-1}.
\eeq

Also, from item i) of Theroem \ref{Reg-univBound} we get

\beq \label{bound-L-infinito-sta-fuori} 
\forall \, n \ge 1 \qquad \Vert u_{\eps_n} \Vert_{L^{\infty}(\R^N)} \le 1.
\eeq

Now, thanks to \eqref{energia-riscalata3-fuori}-\eqref{bound-L-infinito-sta-fuori}, we can apply Theorem 2.1. of \cite{TW} for every $a>0$, to obtain that the limit varifold $V$ satisfies $ Sing(V) \cap (\R^N \setminus \{0\}) = \emptyset $ if $ N \le 7$, has no accumulation points in $\R^N \setminus \{0\}$ if $ N=8$ and has Hausdorff dimension at most $ N-8$ when $ N\ge 9$. Furthermore, $ Reg(V) $ is an embedded smooth {\it stable} minimal hypersurface of $ \R^N$.  Thus, we have proved items v) and vi). 

It remains to prove items iii) and iv).  So far, for $ N \le 7$, we have proved that $ Sing(V) \subseteq \{0\} $ and that $Reg(V)$ is an embedded smooth {\it stable} minimal hypersurface of $ \R^N$. In particular, $supp \Vert V \Vert $ is an 
$(N-1)-$dimensional cone, smooth outside the origin, with zero mean curvature (outside the origin) and satisfying the stability condition

\beq\label{cono-stabile}
\forall \, \varphi \in \mathcal{C}^1_c(M) \qquad \int_{M} \Big ( \vert \nabla^M \varphi \vert^2 - \vert A \vert^2 \varphi^2 \Big ) d {\cal H}^{N-1} \ge 0,
\eeq

\noindent where $ M = supp \Vert V \Vert \setminus \{0\} $, $ A$ is the second fundamental form of $M$ and $ 
\vert A \vert $ is its lenght.  

Since $ 2 \le N-1 \le 6$ we can then invoke a celebrated result of J. Simons \cite{Simons} (see also \cite{SSY, Simon, Giusti}) to conclude that the cone $supp \Vert V \Vert $ must be a hyperplane. 

For $N=8$ we already know that $ Sing(V) \cap (\R^N \setminus \{0\}) $ has no accumulation points in $\R^N \setminus \{0\}$. Let us prove that $ Sing(V) \cap (\R^N \setminus \{0\}) = \emptyset $. If not, since $V$ is a cone, there would exist an open ray of singular points, i.e., an open ray entirely contained in $(supp \Vert V \Vert \setminus \{0\}) \cap Sing(V) $, a clear contradiction.  This concludes the proof.



\end{proof}

\begin{remark}\label{Rk7}
It is interesting to note that, for $ 3 \le N \le 7$, solutions stable outside a compact set and stable solutions share the same properties, that is, the $supp \Vert V \Vert$ of any limit varifold $V$ is a ($N$-1)-dimensional hyperplane through the origin $0 \in \R^N$.  
The analogy is even more surprising if one observes that, in dimension $N=2$, the saddle solution $u$ constructed by Dang, Fife and Peletier in \cite{DFP} is a finite energy solution which is not stable, but which is stable outside a compact set of $\R^2$.  For the finiteness of the energy and its exact value see the Appendix \ref{sec : appendix C}. The instability of $u$ is well-known and can be also directly inferred from Theorem \ref{Teo-stabili}, since the $supp \Vert V \Vert$ coincides with the set of the zeros of $u$ (see the proof of Theorem \ref{Teo-Energia-saddle-2D}), which is not a straight line. The stability outside a compact set is a consequence of Lemma 2.2. in \cite{KL} in which it is proved that the linearized equation around $u$ posseses a positive supersolution outside a large ball (and it is a standard fact that this is equivalent to stabilty of $u$ outside that ball).

\end{remark}


\bigskip

\section{Invariant solutions}\label{sec : InvSol} 

In this section we consider finite energy solutions having some symmetries. For this family of invariant solutions we describe the effect of the presence of these symmetries on the limit varifold $V$. The first result is 

\medskip

\begin{proposition} \label{Sim1}  
Assume $N \ge 2$ and let $ H \in O(N) $ be an orthogonal transformation. Let $u$ be a non-constant solution of \eqref{eqAC} with finite energy and satisfying :

\beq \label{sim1a}
{\it either} \qquad  u(x) = u (H(x)) \qquad \forall \, x \in \R^N,
\eeq
\beq \label{sim1b}
\qquad {\it or} \qquad  u(x) = - u (H(x)) \qquad \forall \, x \in \R^N.  
\eeq

Let $(\eps_n)_{n\ge1}$ be the subsequence and let $V$ be the limit varifold provided by Theorem \ref{Teo-finite-energy} then
\beq \label{dens-sim}
\Theta (V,x) = \Theta (V,(H(x))) \qquad \forall \,\, x \in \R^N
\eeq
and $supp \Vert V \Vert$ is $H$-invariant. 

\end{proposition}

\begin{proof}

For any $ x \in \R^N$ and any $ R>0$ we have 
\beq\begin{split}
& {\mathcal{E}_{\eps_n}} (u_{\eps_n},B(x,R))  = \int_{B(x,R)} \frac{\eps_n \vert \nabla u_{\eps_n} \vert^2}{2}
+\frac{(1 - u_{\eps_n}^2)^2}{4\eps_n} \\
& = \int_{B(x,R)} \frac{1}{\eps_n} \Big[ \frac{\vert \nabla u(\frac{y}{\eps_n})\vert^2}{2}  
+\frac{(1 - u^2(\frac{y}{\eps_n}))^2}{4} \Big]  dy \\
& = \int_{B(x,R)} \frac{1}{\eps_n} \Big[ \frac{\vert \nabla u(H(\frac{y}{\eps_n}))\vert^2}{2}  
+\frac{(1 - u^2(H(\frac{y}{\eps_n})))^2}{4} \Big]  dy \\
& = \int_{B(x,R)} \frac{1}{\eps_n} \Big[ \frac{\vert \nabla u(\frac{H(y)}{\eps_n})\vert^2}{2}  
+\frac{(1 - u^2(\frac{H(y)}{\eps_n}))^2)}{4} \Big]  dy \\
\end{split}
\eeq
thanks to \eqref{sim1a}-\eqref{sim1b} and the linearity of $H$. By the change of variable $ z = H(y) $ we get 
\beq \label{integro-sim}
{\mathcal{E}_{\eps_n}} (u_{\eps_n},B(x,R))  = {\mathcal{E}_{\eps_n}} (u_{\eps_n},H(B(x,R))) =  {\mathcal{E}_{\eps_n}} (u_{\eps_n},B(H(x),R)).
\eeq
For every $x \in \R^N$ we have $\Vert V \Vert (\partial B(x,R)) = \Vert V \Vert (\partial B(H(x),R)) = 0$ for all $R>0$, with an at most countable exception, since $ \Vert V \Vert$ is a Radon measure on $\R^N$. 
Thus, an application of \eqref{integro-sim} immediately yields \eqref{dens-sim}.  

\end{proof}

\medskip

As a first consequence of the previous result we have the following classification result for radially symmetric solutions.

\medskip

\begin{corollary}\label{nonEsist-radiali1}
Assume $N \ge 2$ and let $u$ be a finite energy solution of \eqref{eqAC} which is radially symmetric, i.e.,  such that 
\beq \label{rad}
\qquad  u(x) = u (\vert x \vert ) \qquad \forall \, x \in \R^N.
\eeq
Then either $u \equiv 1 $ or $u \equiv - 1 $.

\end{corollary}

\begin{proof}
Suppose, for contradiction, that $ L_u>0 $. Since $u$ satisfies \eqref{sim1a} for every rotation, the set $supp \Vert V \Vert$ must be a rotationally-invariant stationary, integral ($N$-1)-cone. Thus, $supp \Vert V \Vert = \R^N$, a contradiction. Therefore $L_u =0$ and the conclusion follows from item iv) of Theorem \ref{Energia}.  

\end{proof}

\medskip

In the two-dimensional case the conclusion of the previous corollary can be considerably improved as follows  

\medskip

\begin{corollary}\label{nonEsist-radiali2}
Assume $N = 2$ and let $u$ be a non-constant finite energy solution of \eqref{eqAC} satisfying either \eqref{sim1a} or \eqref{sim1b} for some rotation $H \in SO(2)$.  Then, if {\rm ord(H)} denotes the  order of $H$, we have  

\item i) {\rm ord(H)} is finite,

\item ii) {\rm ord(H)} divides the number $k$ of half-lines composing the cone $V$, 

\item iii) if {\rm ord(H)} is even, then $V$ is the union of a finite number of distinct straight lines through the origin, each of which of constant multiplicity.

\end{corollary}

\medskip

\begin{remark}
VEDI UNICIT\`A DELLA VARIFOLD LIMITE in PRESENZA DI SIMMETRIE.   
\end{remark}

\medskip

\begin{proof} Here we use the notation of Theorem \ref{Teo-finite-energy-2D}. 

\item $i)$ Since  $u$ satisfies either \eqref{sim1a} or \eqref{sim1b}, for every integer $n \ge 1$, the rotation $ H^n $ sends bijectively the half-line $\{ \, t v_1 : t >0 \}$  onto the half-line  $\{ \, t v_j : t >0 \}$, for some $ j=j(n) \in \{1,...,k\}$, where $k$ denotes the number of half-lines composing the cone $V$.  
Therefore it must exist two positive integers $ n_1 < n_2$ such that $H^{n_1} (v_1) = H^{n_2} (v_1) $, which yields 
$ {\rm ord(H)} \le n_2 -n_1$. 

\item $ii)$ As the result is obvious when $H$ is the identity, we may assume (by item  $(i)$) that $H$ is a rotation of angle $\frac{2\pi p}{q}$, where $q>1$ and $p$ are two coprime positive integers. Since no power of $H$ with exponent $ 1, \ldots, q-1$ can fix a half-line composing $V$ (i.e., the action is free) we see that each orbit has exactly $q={\rm ord(H)}$ elements. The total number of half-lines therefore is a multiple of $q$.

\item $iii)$ Since {\rm ord(H)} is even, there is a positive integer $n$ such that $ H^n$ has order two, that is, $H^n$ is the rotation of angle $\pi$ (or $ -\pi$).  By Proposition \ref{Sim1}, for every $ j \in \{1,...,k\}$ the half-line $\{ \, t H^{n}(v_j) : t >0 \} = H^n (\{ \, t v_j : t >0 \})$ and the half-line $\{ \, t v_j : t >0 \}$ have the same multiplicity. This means that $ H^{n}(v_j) = -  v_j $, which concludes the proof. 

\end{proof}

\bigskip

Next we turn to {\it cylidrical solutions}, that is solutions which are independent of one or more Euclidean variables.  In other words, these are solutions $v$ of the Allen-Cahn equation in $\R^N$ of the form 

\beq \label{sim2a}
v(x_1,..., x_N) = u (x_1,...,x_k) \qquad \forall \, x=(x_1,...,x_N) \in \R^N,
\eeq
for some solution $u$ of the Allen-Cahn equation in $ \R^k $ and some integer $k$ satisfying $2 \le k < N$.   

It is immediate to see that $u$ is a finite energy solution in $\R^k$ if and only if $v$ is a finite energy solution in $\R^N$. Thus, we can apply Theorem \ref{Teo-finite-energy} both to $u$ and $v$ to get the existence of a limit stationary, integral ($k$-1)-varifold $V_u$, associated to $u$, and a limit stationary, integral ($N$-1)-varifold $ V_v $, associated to $v$ (and both associated to the same subsequence $(\eps_n)_{n\ge1}$). 

In the next result we shall describe the varifold $V_v$ and its relation with the varifold $V_u$. 

\medskip

\begin{proposition} \label{Sim2}  
Assume $2 \le k < N$ and let $v$ be a non-constant solution of \eqref{eqAC} in $ \R^N$, with finite energy and satisfying \eqref{sim2a} for some solution $u$ of \eqref{eqAC} in $ \R^k $. 

Let $(\eps_n)_{n\ge1}$ be the subsequence and let $V_u, V_v$ be the limit varifolds provided by Theorem \ref{Teo-finite-energy} then the limit varifold $V_v$ is a cylinder with

\beq\label{supp-cil}
supp (V_v) = supp (V_u) \times \R^{N-k}
\eeq
and satisfying 

\beq \label{dens-sim2a}
\Theta (V_v,x) = \Theta (V_v, \pi_1(x)) = \Theta (V_u, x^{'}) \qquad \forall \,\, x \in \R^N,
\eeq
where $\pi_1(x)$ denotes the orthogonal projection of $x$ onto the $k$-dimensional subspace $\R^{k} \times \{0\} $, that is, $ \pi_1(x) = (x_1,..., x_k,0,...,0) : = (x^{'}, 0,...,0) \in \R^k \times \R^{N-k} : = \R^N$.   

In particular, we have 
\beq \label{dens-sim2b}
\Theta (V_v,0) = \Theta (V_u, 0).
\eeq

\end{proposition}

\begin{proof}

For any $ x \in \R^N$ we set $ \pi_1(x) = (x_1,...,x_k,0,...,0) \in \R^N$ and $ \pi_2(x) = (0,...,0,x_{k+1},...,x_N) \in \R^N$. Then, for every bounded Borel set $ \Omega \subset \R^N$ we have 
\beq\begin{split} \label{formula-gen}
& {\mathcal{E}_{\eps_n}} (v_{\eps_n},\Omega)  = \int_{\Omega} \frac{\eps_n \vert \nabla v_{\eps_n} \vert^2}{2}
+\frac{(1 - v_{\eps_n}^2)^2}{4\eps_n} \\
& = \int_{\Omega} \frac{1}{\eps_n} \Big[ \frac{\vert \nabla v(\frac{y}{\eps_n})\vert^2}{2}  
+\frac{(1 - v^2(\frac{y}{\eps_n}))^2}{4} \Big]  dy \\
& = \int_{\Omega} \frac{1}{\eps_n} \Big[ \frac{\vert \nabla v(\frac{y-\pi_2(x)}{\eps_n})\vert^2}{2}  
+\frac{(1 - v^2(\frac{y-\pi_2(x)}{\eps_n}))^2}{4} \Big]  dy \\
\end{split}
\eeq
where in the latter we have used \eqref{sim2a}.

Taking $ \Omega = B^N(x,R)$, $R>0$ in \eqref{formula-gen}, the change of variable $ z = y-\pi_2(x)$ leads to 

\beq\begin{split}\label{formula-energia-generale}
& {\mathcal{E}_{\eps_n}} (v_{\eps_n},B^N(x,R))  =  \int_{B^N(x -\pi_2(x) ,R)} \frac{1}{\eps_n} \Big[ \frac{\vert \nabla v(\frac{z}{\eps_n})\vert^2}{2}  + \frac{(1 - v^2(\frac{z}{\eps_n}))^2}{4} \Big]  dz \\
& = \int_{B^N(\pi_1(x) ,R)} \frac{1}{\eps_n} \Big[ \frac{\vert \nabla v(\frac{z}{\eps_n})\vert^2}{2}  
+\frac{(1 - v^2(\frac{z}{\eps_n}))^2}{4} \Big]  dz \\
& =  \int_{B^{N-k}(x^{''}, R)} \Big[ \int_{B^k(x^{'}, \sqrt{R^2 - \vert \xi - x^{''} \vert^2})} \frac{1}{\eps_n} \Big[ \frac{\vert \nabla u(\frac{\eta}{\eps_n})\vert^2}{2}  + \frac{(1 - u^2(\frac{\eta}{\eps_n}))^2}{4} \Big]  d\eta \, \Big] d\xi \\
&  =  \int_{B^{N-k}(x^{''}, R)}  {\mathcal{E}_{\eps_n}} (u_{\eps_n},B^k(x^{'}, \sqrt{R^2 - \vert \xi - x^{''} \vert^2})) d\xi, \\
\end{split}
\eeq
where $ x^{'} = (x_1,...,x_k) \in \R^k$ and $x^{''} = (x_{k+1},...,x_N) \in \R^{N-k}$.

Now pick $x \in \R^N$. Then, as before we have $\Vert V_v \Vert (\partial B(x,R)) = 0$ for all $R>0$, with an at most countable exception. For each of these $R>0$ and for every $ t \in (0,R)$ there is a point $ \xi_t \in \partial B^{N-k}(x^{''}, \sqrt{R^2 - t^2})$ such that $ \Big \{ \partial B^{k}(x^{'}, \sqrt{R^2 - \vert \xi_t - x^{''}\vert^2}) \Big \}_{t \in (0,R)}$ is a disjoint family of Borel sets in $ \R^k$. Therefore, since $ \Vert V_u \Vert $ is a Radon measure on $ \R^k$ by construction, we obtain that  $\Vert V_u \Vert ( \partial B^{k}(x^{'}, \sqrt{R^2 - \vert \xi_t - x^{''}\vert^2})  )  = 0$ for every $t \in (0,R) \setminus N$, where $ N $ is at most countable. Furthermore, ${\cal L}_{N-k} \Big ( \bigcup_{t \in (0,R) \setminus N}  \partial B^{N-k}(x^{''}, \sqrt{R^2 - t^2}) \Big ) = {\cal L}_{N-k} (B^{N-k}(x^{''}, R ))$, where $ {\cal L}_{N-k}$ denotes the $(N-k)-$dimensional Lebesgue measure. From the previous analysis we infer that, for almost every $\xi \in  B^{N-k}(x^{''}, R )$ (w.r.t the Lebsegue measure ${\cal L}_{N-k}$) 

\beq \label{Leb1}
{\mathcal{E}_{\eps_n}} (u_{\eps_n},B^k(x^{'}, \sqrt{R^2 - \vert \xi - x^{''} \vert^2})) \longrightarrow  c_{\star} \Vert V_u \Vert ( B^{k}(x^{'}, \sqrt{R^2 - \vert \xi - x^{''}\vert^2})  ) 
\eeq
as $ n \to +\infty$. On the other hand, $u$ has finite energy so  

\beq\label{Leb2}
 {\mathcal{E}_{\eps_n}} (u_{\eps_n},B^k(x^{'}, \sqrt{R^2 - \vert \xi - x^{''} \vert^2}))  \le L_u \Big ( \sqrt{R^2 - \vert \xi - x^{''} \vert^2} \Big )^{k-1} \le L_u R^{k-1}
\eeq
for every $\xi \in  B^{N-k}(x^{''}, R )$. \eqref{Leb1} and \eqref{Leb2} enable us to use Lebesgue dominated convergence theorem to pass to the limit in \eqref{formula-energia-generale} to get

\beq\begin{split}\label{formula-densità-generale}
&\Vert V_v \Vert (B^N(x,R))  =  \int_{B^{N-k}(x^{''}, R)}  \Vert V_u \Vert ( B^{k}(x^{'}, \sqrt{R^2 - \vert \xi - x^{''}\vert^2}))  d\xi \\
& = \int_{B^{N-k}(0, R)}  \Vert V_u \Vert ( B^{k}(x^{'}, \sqrt{R^2 - \vert z \vert^2}))  dz
\end{split}
\eeq
for all $R>0$, with an at most countable exception. 
The latter tell us that $x \in supp \Vert V_v \Vert$ if and only if $ x^{'} \in supp \Vert V_u \Vert$ and thus $ V_v$ is a cylinder of the form \eqref{supp-cil}.

To prove \eqref{dens-sim2a} and \eqref{dens-sim2b} we observe that by definition of the density function we have 

\beq
\forall \, n \ge 1, \, \exists \, 0 < R_n < \frac{1}{n} \quad : \quad \forall \, 0 < R \le R_n\qquad 
\Big \vert \frac{\Vert V_u \Vert (B^k(x^{'},R)) }{\omega_{k-1} R^{k-1}} - \Theta (V_u, x^{'}) \Big \vert < \frac{1}{n}
\eeq
thus
\beq\label{infinitesimo}
\forall \, n \ge 1, \, \forall \, z \in  B^{N-k}(0, R_n) \qquad 
\Big \vert \frac{\Vert V_u \Vert (B^k(x^{'},  \sqrt{R_n^2 - \vert z \vert^2})) }{\omega_{k-1} (\sqrt{R_n^2 - \vert z \vert^2}))^{k-1}} - \Theta (V_u, x^{'}) \Big \vert < \frac{1}{n}
\eeq
From \eqref{formula-densità-generale} we then get 
\beq\begin{split}\label{quasi-fine}
&\Vert V_v \Vert (B^N(x,R_n)) = \Theta (V_u, x^{'})  \int_{B^{N-k}(0, R_n)} \omega_{k-1} 
(\sqrt{R_n^2 - \vert z \vert^2}))^{k-1}  d \, z \\
& + \int_{B^{N-k}(0, R_n)} (\Vert V_u \Vert (B^k(x^{'}, \sqrt{R_n^2 - \vert z \vert^2})) - \Theta (V_u, x^{'})  \omega_{k-1} (\sqrt{R_n^2 - \vert z \vert^2}))^{k-1} ) d\, z  \\
&= I_1(n) + I_2(n). 
\end{split}
\eeq
Now we observe that 

\beq\label{quasi-fine-palla}
I_1(n) =  \Theta (V_u, x^{'}) \omega_{N-1} R_n^{N-1}
\eeq
and
\beq\label{quasi-fine-resto}
\vert I_2(n) \vert  \le  \frac{1}{n} \int_{B^{N-k}(0, R_n)} \omega_{k-1} (\sqrt{R_n^2 - \vert z \vert^2}))^{k-1} ) d\, z =  \frac{1}{n} \omega_{N-1} R_n^{N-1}, 
\eeq
where in the latter we have used \eqref{infinitesimo}. Hence

\beq\label{}
\Theta( V_v, x) = \lim_{n \to +\infty} \frac{\Vert V_v \Vert (B^N(x,R_n))}{\omega_{N-1} R_n^{N-1}} = \lim_{n \to +\infty} \frac{ I_1(n) + I_2(n)}{\omega_{N-1} R_n^{N-1}} =   \Theta (V_u, x^{'}) 
\eeq
thanks to \eqref{quasi-fine-palla} and \eqref{quasi-fine-resto}.





\end{proof}

\section{Ordered solutions}\label{sec : Sol-Ordinate}

In this section we deal with ordered solutions to the Allen-Cahn equation, i.e., with solutions $u$ and $v$ of \eqref{eqAC}, such that $ u\le v$ on $ \R^N$. We shall focus either on pair of ordered local minimizers or on pair of ordered finite energy solutions. For local minimizers we are interested in the comparison of the limit minimal cones $C_u$ and $C_v$ provided by Theorem \ref{Teo-minloc} while, for general finite energy solutions we would like to compare the limit varifold $V_u$ and $V_v$ provided by Theorem \ref{Teo-finite-energy}. In both cases, the limit objects are always assumed to be associated to the same subsequence $(\eps_n)_{n\ge1}$ (note that this is always possible, after passing, if necessary, to another subsequence).


To this end, for any finite energy solution $w$, we consider the open set $ X_{w} := \R^N \setminus supp \Vert V_{w} \Vert $ and we recall (see Theorem 1 and Proposition 4.2 of \cite{HT}) that  $\eps_n \frac{\vert \nabla w_{\eps_n} \vert^2}{2}
+ \frac{(1 - w_{\eps_n}^2)^2}{4\eps_n} $ converges locally uniformly to {\it zero} on $X_w$. In particular, $w_{\eps_n} $ converges locally uniformly to $\pm 1$ on $ X_{w}. $ We can therefore define the sets $ X^{+}_{w} := \{ \, x \in X_{w} \, : \, w_{\eps_n}(x) \to 1 \}$ and $ X^{-}_{w} := \{ \, x \in X_{w} \, : \, w_{\eps_n}(x) \to -1 \}$ (which are also open due to the above mentioned locally uniform convergence property). With this notations we have 

\medskip

\begin{lemma} 
Assume $ N \ge 2$ and let $w$ be a finite energy solution of the Allen-Cahn equation \eqref{eqAC}. Then, $X_{w} , X^{+}_{w} $ and each connected component of $ X^{+}_{w}$ are open cones of $\R^N$. The same conclusion holds true for $X^{-}_{w} $ and each connected component of $ X^{-}_{w}$. 
\end{lemma}

\medskip

\begin{proof}
We need to prove only the first part of the claim. The second part follows by applying the first one to the solution $-w$. 

Since $X_{w}$ and $X^{+}_{w} $ are open sets of $\R^N$, also each connected component of $ X^{+}_{w}$ is an open set of $\R^N$. Moreover, $X_{w}$ is a cone since it is complement,  $ supp \Vert V_{w} \Vert$, is a cone. 

Pick $x \in X^{+}_{w} $ and suppose there is $t>0$ such that $tx \in  supp \Vert V_{w} \Vert \cup  X^{-}_{w} $.  The point $tx$ cannot belong to $ supp \Vert V_{w} \Vert$, since this set is a cone. Thus $ tx \in  X^{-}_{w} $. But this implies the existence of a point $Tx \in supp \Vert V_{w} \Vert$, for some $T>0$. As before, this is impossible, since $ supp \Vert V_{w} \Vert$ is a cone.
 
Let $ \omega$ be a connected component of $X^{+}_{w}$ and fix $x \in \omega$. Let us prove that the open line segment $ \{ \, tx : t>1 \} $ belongs to $\omega$. Since $ \omega$ is open, the set $ L_x := \{ \, s >1 \, \, : \, \, \forall \,  t \in (1, s), \,\,  tx \in \omega \} $ is not empty. To conclude we shall prove that $ S : = \sup L_x = + \infty$. If not, then $1 < S < +\infty$, $Sx \in  X^{+}_{w}$ and there exists $R(Sx) >0$ such that $B(Sx, R(Sx)) \subset X_w^{+}$. Thus, by the definiton of $S$, the point $Sx$ must belong to another connected component of $X_u^{+}$, say $ \omega^{'} (\not \equiv \omega) $. Since $ \omega^{'}$ is open, there would exist $R^{'}(Sx) >0$ such that $B(Sx, R^{'}(Sx)) \subset \omega^{'}$,which entails $ \omega \cap \omega^{'} \neq  \emptyset$, a contradiction.  A similar argument proves that the open line segment $ \{ \, tx : 0 < t <1\} $ belongs to $\omega$.
\end{proof}

\medskip

\begin{proposition} \label{Teo-Ord}  
Assume $N \ge 2$ and let $u,v$ be two non-constant solution of \eqref{eqAC} satisfying :

\beq \label{ordinate}
u(x) \le v(x) \qquad \forall \, x \in \R^N,
\eeq

\item[i)] if $u$ and $v$ are local minimizers, then 

\beq \label{ordminloc}
C_u = C_v 
\eeq
where $C_u$ and $C_v $ denote the limit minimal cones (associated to $u$ and $v$ resp.) provided by Theorem \ref{Teo-minloc}. 

\item[ii)] if $u$ and $v$ have both finite energy, then 

\beq \label{ordfinite}
X^{+}_{u} \subseteq X^{+}_{v}, \qquad \quad X^{-}_{v} \subseteq X^{-}_{u}
\eeq
\end{proposition}

\medskip

\begin{proof} Since $u$ and $v$ are ordered, we have that 

\beq\label{ordeps}
 u_{\eps_n}(x) \le v_{\eps_n}(x) \qquad \forall \, x \in \R^N, \qquad \forall n \ge 1.
\eeq

i) If $ x \in C_u$, then $ u_{\eps_n}(x) \rightarrow 1 $ and, by \eqref{ordeps}, we also have that $ v_{\eps_n}(x) \rightarrow 1$. Thus $ C_u \subseteq C_v$, but then, the minimality of both cones implies that they must coincide. 

\smallskip

ii) Consider $ x \in X_u^{+}$. By the above mentioned locally uniform convergence property of $u_{\eps_n}$ there exists $R(x) >0$ such that $B(x, R(x)) \subset X_u^{+}$, $\frac{(1 - u_{\eps_n}^2)^2}{4\eps_n} \longrightarrow 0$ uniformly on $B(x, R(x))$  and also $u_{\eps_n} \longrightarrow 1$ uniformly on $B(x, R(x))$.  
Thus \eqref{ordeps} yields that $ \frac{(1-u_{\eps_n}^2)^2}{4 \eps_{n}} \ge  \frac{(1-v_{\eps_n}^2)^2}{4 \eps_{n}} $ on  $B(x, R(x))$ for $n$ large enough. On the other hand, we also have that $ \frac{(1-v_{\eps_n}^2)^2}{4 \eps_{n}} \ge {\eps_n} \frac{\vert \nabla v_{\eps_n} \vert^2}{2} $ on $\R^N $ for every $n \ge 1$,  by Modica's gradient bound. To summarize we have 

\beq\label{energia-ord}
\int_{B(x, R(x))} \frac{\eps_n \vert \nabla v_{\eps_n} \vert^2}{2}
+\frac{(1 - v_{\eps_n}^2)^2}{4\eps_n} \le 2 \int_{B(x, R(x))} \frac{(1 - u_{\eps_n}^2)^2}{4\eps_n} \longrightarrow 0
\eeq
The latter and the definition of $\Vert V_v \Vert$ immediately imply that $ x \not\in supp \Vert V_v \Vert $. 
Since also $ v_{\eps_n} \longrightarrow 1$ uniformly on $B(x, R(x))$ we see that $ x \in X_v^{+}$, which proves the first inclusion  in \eqref{ordfinite}. The second one follows by a similar argument applied to $ X^{-}_{v} $. 

\end{proof}

\section{Monotone solutions}\label{sec : Soluzioni-Monotone}

This section is motivated by the following celebrated question raised by E. De Giorgi in 1978 (\cite{CongDeG}) :

\bigskip

\noindent {\bf De Giorgi's conjecture. } {\it Let $u\in \mathcal{C}^2(\R^N, [-1,1])$ satisfy

\beq\label{congDeG}
-\Delta u= u-u^3\qquad and \qquad\frac{\partial u}{\partial x_N}>0\qquad on \quad \R^N.
\eeq

\noindent Is it true that  all the level sets of~$u$ are hyperplanes, at least if $N\le 8$? }

\bigskip

It is easily seen that the claim of De Giorgi's conjecture is equivalent to $u$ being 1D and monotone, i.e., $u$ being the one-dimensional monotone solution 

$$u=u(x) = \tanh\Big (\frac{(\nu\cdot x) +\alpha}{\sqrt 2} \Big) \quad \qquad \forall \, x \in \R^N,$$
for some $\alpha \in \R$ and some unit vector $ \nu \in \R^N$ such that $ \nu_N := (\nu \cdot e_N) >0$.  

\bigskip

De Giorgi's conjecture has been proved in dimension $ N=2$ by N. Ghoussoub and C. Gui \cite{GG} (see also \cite{BCN2}) and in dimension $ N=3$ by L. Ambrosio and X. Cabr\'e \cite{AC}. {\it The original conjecture is still open when $ 4 \le N \le 8$. }

\smallskip

O. Savin \cite{S} proved the validity of the conjecture for $4 \le N \le 8$ under the following additional assumptions :   
\beq\label{limiti-cost}
\lim_{x_N\rightarrow-\infty} u(x',x_N)=-1
\quad {\mbox{ and }} \quad
\lim_{x_N\rightarrow+\infty} u(x',x_N)=1.
\eeq

Later, in \cite{FVtams}, the present author and E. Valdinoci proved De Giorgi's conjecture (and other symmetry results) under conditions at infinity which are weaker than \eqref{limiti-cost}. In particular, it is established in \cite{FVtams} that condition \eqref{limiti-cost} can be relaxed to suitable symmetry assumptions on the asymptotic profiles
\begin{equation}\label{N-profiles-conj}\begin{split}
&\underline u(x'):=\lim_{x_N\rightarrow-\infty} u(x',x_N), 
\\ 
&\overline u(x'):=\lim_{x_N\rightarrow+\infty} u(x',x_N).
\end{split}\end{equation}

More precisely, in \cite{FVtams} it is proved that solutions of \eqref{congDeG} 
are 1D if :

\begin{itemize}
\item $2\le N\le 4$ and at least one between~$\overline u$
and~$\underline u$ are $2$D,
\item $2\le N\le 8$ and both~$\overline u$
and~$\underline u$ are $2$D.
\end{itemize}

\medskip

On the other hand, for $N \ge 9$, M. del Pino, M. Kowalczyk and J. Wei have constructed in \cite{DKW} a monotone solution to \eqref{eqAC} which is not 1D, showing that the dimensional constraint in De Giorgi's conjecture cannot be removed.

\smallskip

We also refer to \cite{FVstate} for further results and a review on De Giorgi's conjecture.

\medskip

In \cite{FVtams} (see Theorem 1.3 therein) we also proved the following result for monotone solutions 

\begin{theorem}{\bf (Symmetry from a level set being a graph)} \label{livello-grafico-intero} 
Assume $2 \le N \le 8 $ and let $u$ be a non-constant, monotone solution of \eqref{eqAC}. 

\noindent Suppose that, at least one level set $u$ is an entire graph over $\R^{N-1}$. 

\smallskip

\noindent Then $u$ is 1D, that is, $u$ is of the form 

$$u=u(x) = \tanh\Big (\frac{(\nu\cdot x) +\alpha}{\sqrt 2} \Big) \quad \qquad \forall \, x \in \R^N,$$
for some $\alpha \in \R$ and some unit vector $ \nu \in \R^N$ such that $ \nu_N := (\nu \cdot e_N) >0$.  
\end{theorem}

\bigskip

The main results of this section concern the study and the classification of {\it monotone} solutions {\it with finite energy.}  Without loss of generality, by {\it monotone solution} we always denote a solution $u$ of the Allen-Cahn equation \eqref{eqAC} satisfying 

\beq 
\frac{\partial u}{\partial x_N}(x)>0\qquad \forall \, x \in \R^N,
\eeq
as in De Giorgi's conjecture.

\medskip

Let us start with the following auxiliary result. It makes the relation between the energy of a monotone solution $u$ and  the energies of its asymptotic profiles $\overline u$ and $\underline u$, defined by \eqref{N-profiles-conj}. It will be ancillary to the main results of the present section. 

\bigskip

\begin{proposition}\label{Energie-u-limiti}
Assume $N \ge 2$. Let $u$ be a non-constant, monotone solution of \eqref{eqAC}. Then, the following assertions are equivalent :

\smallskip

\item[i)] $u$ has finite energy,

\item[ii)] both the profiles $\overline u$ and $\underline u$ have finite energy,

\item[iii)] at least one between $\overline u$ and $\underline u$ has finite energy.

\medskip

\noindent Furthermore, when $u$ has finite energy we have the following relations 

\beq\label{relEn-u-profili}
\vert L_u - L_{\overline u}(N) \vert \le \sqrt{2} s_{N-1}, \quad  \vert L_u - 
L_{\underline u}(N) \vert \le \sqrt{2} s_{N-1}, \quad \vert L_{\underline u}(N) - L_{\overline u}(N) \vert \le \sqrt{2} s_{N-1}, 
\eeq
where $s_{N-1}$ is the surface area of the unit sphere of $\R^N$. \\
Where, in \eqref{relEn-u-profili} above, with a slight abuse of notation, $L_{\overline u}(N)$ (resp. $L_{\underline u}(N)$) denotes the energy of ${\overline u}$ (resp. $\underline u$) considered as a function of $N$ variables. 
\end{proposition}

\medskip

\begin{proof}
We use a dynamical system-type approach as in \cite{AC}. For $ t \in \R $ and $ x \in \R^N$ we set $ u^{t} (x) := 
u(x^{'}, x_N + t)$. For every $ R>0$ we then obtain

\beq\begin{split}
& \frac{d}{dt} {\mathcal{E}}(u^{t}, B(0,R))  = \int_{B(0,R)} \nabla u^{t} \cdot \nabla (\partial_t u^{t}) 
+  [(u^{t})^3 - u^t]\partial_t u^{t} \\
& = \int_{\partial B(0,R)} \partial_t u^{t} (\nabla u^{t} \cdot \nu)  d\sigma \ge - \Big (\sup_{x \in \R^N} \vert \nabla u(x) \vert \Big) \int_{\partial B(0,R)} \partial_t u^{t}  d\sigma,
\end{split}
\eeq
where in the latter we have used the monotonicity assumption $ \frac{\partial u}{\partial x_N}  > 0$ on $ \R^N$. Integrating the above inequality over $ [S,T]$ and using Tonelli's Theorem we get
\beq\begin{split}
& {\mathcal{E}}(u^{T}, B(0,R))-  {\mathcal{E}}(u^{S}, B(0,R)) = \int_S^T \frac{d}{dt} {\mathcal{E}}(u^{t}, B(0,R))  dt \\
& \ge -  \Big(\sup_{x \in \R^N} \vert \nabla u(x) \vert \Big)  \int_S^T \int_{\partial B(0,R)} \partial_t u^{t}  d\sigma dt = - \Big(\sup_{x \in \R^N} \vert \nabla u(x) \vert \Big) \int_{\partial B(0,R)}  \int_S^T  \partial_t u^{t}  dt d\sigma \\
& = - \Big(\sup_{x \in \R^N} \vert \nabla u(x) \vert \Big) \int_{\partial B(0,R)} (u^{T}-u^{S})  d\sigma.
\end{split}
\eeq
Therefore, denoting by $s_{N-1}$ the surface area of the unit sphere of $\R^N$,  we have that for every $ S < T $ and every $R>0$

\beq\label{En-u-e-profili}\begin{split}
&\vert  {\mathcal{E}}(u^{T}, B(0,R)) -  {\mathcal{E}}(u^{S}, B(0,R)) \vert \le \Big(\sup_{x \in \R^N} \vert \nabla u(x) \vert \Big) \int_{\partial B(0,R)} (u^{T}-u^{S})  d\sigma \\
& \le \Big(\sup_{x \in \R^N} \vert \nabla u(x) \vert \Big) 2 s_{N-1} R^{N-1}\le \sqrt{2} s_{N-1} R^{N-1},
\end{split}
\eeq
where in the latter we have used the universal bound $ \Vert u \Vert_{\infty} \le 1$ (cf. item i) of Theorem \ref{Reg-univBound}) and the Modica's gradient bound 
\beq 
\frac{|\nabla u|^2}{2}  \le \frac{(1-u^2)^2}{4} \qquad {\rm on} \qquad \R^N. 
\eeq
Taking $ S=0$ and letting $T \to+\infty$ in \eqref{En-u-e-profili} we get

\beq\label{En-u-profili1}
\forall \, R >0 \qquad \vert  {\mathcal{E}}(\overline u, B(0,R)) -  {\mathcal{E}}(u, B(0,R)) \vert \le \sqrt{2} s_{N-1} R^{N-1}.
\eeq
Taking $ T=0$ and letting $S \to-\infty$ in \eqref{En-u-e-profili} we get

\beq\label{En-u-profili2}
\forall \, R >0 \qquad \vert  {\mathcal{E}}(\underline u, B(0,R)) -  {\mathcal{E}}(u, B(0,R)) \vert \le \sqrt{2} s_{N-1} R^{N-1}.
\eeq
Finally, letting $S \to -\infty, T \to+\infty$ in \eqref{En-u-e-profili} we get

\beq\label{En-u-profili3}
\forall \, R >0 \qquad \vert  {\mathcal{E}}(\overline u, B(0,R)) -  {\mathcal{E}}(\underline u, B(0,R)) \vert \le \sqrt{2} s_{N-1} R^{N-1}.
\eeq
The desired conclusion then follows immediately from \eqref{En-u-profili1}-\eqref{En-u-profili3} and the definition of energy. 
\end{proof}

\medskip

Now we can prove De Giorgi's conjecture in dimension $N=4$ for monotone, finite energy solutions.

\medskip

\begin{theorem}\label{DeG4}
Let $u\in \mathcal{C}^2(\R^4)$ satisfy

\beq
-\Delta u= u-u^3\qquad and \qquad\frac{\partial u}{\partial x_4}>0\qquad on \quad \R^4.
\eeq

\noindent If $u$ has finite energy, then $u$ is 1D, that is, $u$ is of the form 

$$u=u(x) = \tanh\Big (\frac{(\nu\cdot x) +\alpha}{\sqrt 2} \Big) \quad \qquad \forall \, x \in \R^4,$$
for some $\alpha \in \R$ and some unit vector $ \nu \in \R^4$ such that $ \nu_4 := (\nu \cdot e_4) >0$.  

\end{theorem} 

\medskip

\begin{proof}
By Proposition \ref{Energie-u-limiti} we know that the limit profiles $\overline u$ and $\underline u$ have finite energy, in particular they satisfy 

\beq
\forall \, R >0 \qquad \int_{B^4(0,R)} \vert \nabla u_{\overline u} \vert^2 \le 2L_{\overline u} R^3, \qquad
\int_{B^4(0,R)} \vert \nabla u_{\overline u} \vert^2 \le 2L_{\overline u} R^3. 
\eeq
On the other hand, both $\overline u$ and $\underline u$ depend only on the first 3 Euclidean variables and we can find a constant $ C>0 $, independent of $R$ (depending only on $ L_{\overline u}, L_{\underline u}$ and the dimension $N=4$),  such that

\beq\label{stima-ener-3D}
\forall \, R >0 \qquad \int_{B^3(0,R)} \vert \nabla^{'} u_{\overline u} \vert^2 \le C R^2, \qquad
\int_{B^3(0,R)} \vert \nabla^{'} u_{\overline u} \vert^2 \le C R^2. 
\eeq

Also, both $\overline u$ and $\underline u$ are stable solutions of the Allen-Cahn equation \eqref{eqAC} in $\R^3$, because they are the limit profiles of a monotone (and hence stable) solution. Since $\overline u$ and $\underline u$ are stable solutions of the Allen-Cahn equation \eqref{eqAC} in $\R^3$ and \eqref{stima-ener-3D} is in force, an application of Theorem \ref{Teo-stabili-folklore} yields that both $\overline u$ and $\underline u$ are 1D.  But now, as recalled in the first part of this section, we can apply Theorem 1.1. (or Theorem 1.2) in our previous work \cite{FVtams} to conclude that also $u$ is 1D. 

\end{proof}

\medskip

\begin{remark}\label{Rk-monot-finite-energy}
Note that in $\R^2$ and $\R^3$ monotone solutions are automatically of {\it finite energy} (and thus 1D). Indeed, the limit profiles $\overline u$ and $\underline u$ are 1D (resp. 2D) stable solutions of the Allen-Cahn equation in $\R$ (resp. $\R^2$).  
Therefore, both $\overline u$ and $\underline u$ have finite energy. This follows by applying Theorem \ref{stab-1D}  and Theorem \ref{Teo-clas1D} when $ N=2$ and by Theorem \ref{Teo-stabili-folklore} for $N=3$, since every solution of the Allen-Cahn equation has bounded gradient in any dimension (either by item i) of Theorem \ref{Reg-univBound} and standard elliptic estimates (cf. \eqref{stime-derivate}) or by item i) of Theorem \ref{Reg-univBound} and Modica's gradient bound). 

This remark and the above Theroem \ref{DeG4} strongly suggest the validity of the original De Giorgi's conjecture in $\R^4$. 

\end{remark}

\medskip

In the introductory part of the present section we recalled that in \cite{FVtams} we proved the validity of De Giorgi's conjecture for every $ N \le 8$ if both the profiles at infinty $\overline u$ and $\underline u$ are 2D. In the next Theorem we improve this result by further relaxing the 2D condition at infinity. We shall prove that De Giorgi's conjecture holds true for $N \le 8$, if  $u$ has finite energy and both the profiles at infinty $\overline u$ and $\underline u$ are 3D. 

\medskip

\begin{theorem}\label{DeG3-prof3D}
Assume $ 3 \le N \le 8$. Let $u\in \mathcal{C}^2(\R^N)$ satisfy

\beq
-\Delta u= u-u^3\qquad and \qquad\frac{\partial u}{\partial x_N}>0\qquad on \quad \R^N.
\eeq

\noindent Suppose that $u$ has finite energy and that both $\overline u$ and $\underline u$ are 3D. 

\medskip

\noindent Then $u$ is 1D, that is, $u$ is of the form 

$$u=u(x) = \tanh\Big (\frac{(\nu\cdot x) +\alpha}{\sqrt 2} \Big) \quad \qquad \forall \, x \in \R^N,$$
for some $\alpha \in \R$ and some unit vector $ \nu \in \R^N$ such that $ \nu_N := (\nu \cdot e_N) >0$.  

\end{theorem} 

\medskip

\begin{proof}
By Proposition \ref{Energie-u-limiti} we see that both $\overline u$ and $\underline u$ have finite energy, and so they satisfy 

\beq
\forall \, R >0 \qquad \int_{B^N(0,R)} \vert \nabla u_{\overline u} \vert^2 \le 2L_{\overline u} R^{N-1}, \qquad
\int_{B^N(0,R)} \vert \nabla u_{\overline u} \vert^2 \le 2L_{\overline u} R^{N-1}. 
\eeq
On the other hand, both $\overline u$ and $\underline u$ depend only on 3 variables, by the 3D assumption. Therefore we can find a constant $ C>0 $, independent of $R$ (depending only on $ L_{\overline u}, L_{\underline u}$ and the dimension $N$),  such that

\beq\label{stima-ener-(N-1)D}
\forall \, R >0 \qquad \int_{B^3(0,R)} \vert \nabla^{'} u_{\overline u} \vert^2 \le C R^2, \qquad
\int_{B^3(0,R)} \vert \nabla^{'} u_{\overline u} \vert^2 \le C R^2. 
\eeq

Since $u$ is monotone, both $\overline u$ and $\underline u$ are stable solutions of the Allen-Cahn equation \eqref{eqAC} in $\R^3$ satisfying \eqref{stima-ener-(N-1)D}. By applying Theorem \ref{Teo-stabili-folklore} we see that both $\overline u$ and $\underline u$ are 1D.  To conclude that also $u$ is 1D, we invoke once again Theorem 1.1. (or Theorem 1.2) of \cite{FVtams}.

\end{proof}

\medskip

The following result is a variation on the same theme, but differently from the previous results, it does not require the finiteness of the energy. 

\medskip

\begin{theorem}
Assume $ 3 \le N \le 8$. Let $u\in \mathcal{C}^2(\R^N)$ satisfy

\beq
-\Delta u= u-u^3\qquad and \qquad\frac{\partial u}{\partial x_N}>0\qquad on \quad \R^N.
\eeq

\noindent Suppose that $\overline u$ is 2D and that $\underline u$ is 3D. 

\medskip

\noindent Then $u$ is 1D, that is, $u$ is of the form 

$$u=u(x) = \tanh\Big (\frac{(\nu\cdot x) +\alpha}{\sqrt 2} \Big) \quad \qquad \forall \, x \in \R^N,$$
for some $\alpha \in \R$ and some unit vector $ \nu \in \R^N$ such that $ \nu_N := (\nu \cdot e_N) >0$.  

\end{theorem} 

\medskip

\begin{remark}
Note that we {\it do not} assume that $u$ has finite energy. This is a consequence of the others hypotheses. Also, the same conclusion holds if the assumptions on $\overline u$ and $\underline u$ are exchanged. 
\end{remark}

\medskip

\begin{proof}
Since $\overline u$ is 2D, an application of Theorem \ref{Teo-stabili-folklore} with $ N=2$ tell us that $\overline u$ is 1D and with finite energy. We can therefore invoke Proposition \ref{Energie-u-limiti} to get that also $u$ has finite energy. The claim then follows by applying Theorem \ref{DeG3-prof3D}. 
 
\end{proof}

\medskip

Then we provide a one-dimensional symmetry criterium for monotone solutions in any dimension $ N \ge 2$.  It involves, in a natural way,  the limit (interface) conical varifold $V$ provided by the Structure Theorem for finite energy solutions, i.e., Thereom \ref{Teo-finite-energy}. This result will be useful also in the sequel. 

\medskip

\begin{theorem}\label{non-vertical} Assume $N \ge 2$ and let $u$ be a non-constant, monotone solution to \eqref{eqAC} with finite-energy. Then u is 1D if and only if a limit conical varifold $V$ is a non-vertical hyperplane. 

\end{theorem}

\medskip

\begin{remark}\label{RkDeG}

Note that, if $ N \le 7$, the varifold $V$ is {\underline {always}} a hyperplane through the origin by Theorem \ref{Teo-stabili}. This strongly suggests that De Giorgi's conjecture is true for finite energy solutions, at least up to dimension 7, as already fully proved in the above Theorem \ref{DeG4} for $N=4$. 
 
\end{remark}

\medskip

We also refer to Corollary 1.10 of our recent work \cite{FVcomPDE}, where we proved a similar result but under more restrictive assumptions on the limit interface. 

\medskip

\begin{proof}[Proof of Theorem \ref{non-vertical}]
If $u$ is 1D, i.e., of the form 

$$u=u(x) = \tanh\Big (\frac{(\nu\cdot x) +\alpha}{\sqrt 2} \Big) \quad \qquad \forall \, x \in \R^N,$$
for some $\alpha \in \R$ and some unit vector $ \nu \in \R^N$ such that $ \nu_N := (\nu \cdot e_N) \neq 0$, then $supp \Vert V \Vert$ is the non-vertical ($N$-1)-dimensional hyperplane $\{\, \nu \cdot x =0 \,\} $. 

Coversely, if $V$ is a non-vertical hyperplane through the origin, the half-line $(0',-t)$, $t>0$, is not in $supp \Vert V \Vert$ and so

\begin{equation}
u(0', -\frac{1}{\eps_n}) = u_{\eps_n}(0',-1)  \to \pm 1.
\end{equation}
On the other hand, the monotonicity of $u$ also provides  
\begin{equation}
{\underline u} (0') = \lim_{n\to +\infty} u(0', -\frac{1}{\eps_n}),
\end{equation}
where $ {\underline u}$ denotes, as usual, the profile of $u$ at $- \infty$.  Therefore,
\begin{equation}
{\underline u} (0') = \pm 1,
\end{equation}
and the strong maximum principle implies that ${\underline u} \equiv \pm 1$. The monotonicity assumption on $u$ then forces ${\underline u} \equiv -1$. The same argument applies to the half-line $(0',t)$, $t>0$, and provides ${\overline u} \equiv 1$.

We have demostrated that
\begin{equation}
\lim_{x_N\rightarrow-\infty} u(x',x_N)=-1
\quad {\mbox{ and }} \quad
\lim_{x_N\rightarrow+\infty} u(x',x_N)=1,
\end{equation}
and so $u$ is a non-constant local minimizer by \cite{JM}.

We conclude by providing two different proofs.
 
\medskip 
 
I) Since $u$ is a local minimizer, $V$ is a multiplicity 1 hyperplane. Therefore, $u$ attains the {\it minimum positive energy level} and an application of Theorem \ref{1D-minloc} immediately gives the desired conclusion.

\smallskip

II) The aimed conclusion follows from the last part of Theorem \ref{1D-minloc}, when $N \le  7$. From Theorem \ref{monot-1D-minloc}, when $N=8.$ From Theorem \ref{monot-cono-cilindro} (or from Theorem \ref{monot-con-limiti-cost}), when $ N\ge 9$.   
 
\medskip

\end{proof}

\medskip

We also have the following

\medskip

\begin{theorem} \label{monot-rad-sim-N-1}  Assume $N \ge 2$ and let $u$ be a non-constant monotone solution to \eqref{eqAC} such that  $ u =u(\vert x^{'} \vert, x_N) $. Then  

$$u=u(x) = \tanh\Big (\frac{x_N +\alpha}{\sqrt 2} \Big) \quad \qquad \forall \, x \in \R^N,$$
for some $\alpha \in \R$.   
\end{theorem}

\begin{proof}
By the monotonicity assumption, the two limit profiles $\overline u$ and $\underline u$ are radially symmetric stable solutions of the Allen-Cahn equation on $\R^{N-1}$, i.e., of the form $\overline u = \overline u(x^{'})$ and $\underline u = \underline u(x^{'})$. We can then apply item ii) of Theorem \ref{Solradiali} to infer that both $\overline u$ and $\underline u$ are constant solutions and thus $\overline u = 1$ and $\underline u = - 1$ (again by stability). Therefore we have 

\begin{equation}
\lim_{x_N\rightarrow-\infty} u(x',x_n)=-1
\quad {\mbox{ and }} \quad
\lim_{x_N\rightarrow+\infty} u(x',x_n)=1.
\end{equation}

The latter immediately imply that $u$ is a local minimizer and to conclude we can proceed as in the proof of Theorem \ref{monot-cono-cilindro}. By the monotonicity assumption and   the above asymptotics for $u$, we see that the zero level set of $u$ is a smooth function $ \varphi : \R^{N-1} \rightarrow \R$ such that 

\beq
\{ \,  u < 0 \, \} = \{ (x^{'},x_N) \in \R^{N-1} \times \R \, : \, x_N < \varphi(x^{'}) \, \} 
\eeq
and
\beq\label{ins-livello-sim}
\forall \, x^{'} \, \in \R^{N-1} \qquad \varphi (x^{'})= \varphi (\vert x^{'} \vert)
\eeq

\smallskip

\noindent since $u(x) =u(\vert x^{'} \vert, x_N) $ for every $ x \in \R^N$. 

\smallskip

\noindent Therefore
\beq
\{ \,  u_{\eps_n} < 0 \, \} = \{ (x^{'},x_N) \in \R^{N-1} \times \R \, : \, x_N < \varphi_n(x^{'}) \, \} : = U_n
\eeq
where 
\beq
\varphi_n (x^{'}) : = \eps_n  \varphi \Big (\frac{x^{'}}{\eps_n} \Big)
\eeq
and
\beq\label{ins-livello-risc-sim}
\forall x^{'} \, \in \R^{N-1} \qquad \varphi_n (x^{'})= \varphi_n (\vert x^{'} \vert)
\eeq

By item i) of Theorem \ref{Teo-minloc} and Lebesgue dominated convergence theorem we immediately have that

\beq
 {\bf 1}_{U_n} \longrightarrow {\bf 1}_{\R^N \setminus C} \quad  {\text {in}} \quad L^1_{loc}(\R^N).
\eeq
The latter and Lemma 16.3 of \cite{Giusti} tell us that $\R^N \setminus C$ is itself a subgraph of a measurable function 
$ \varphi_{\infty} : \R^{N-1} \rightarrow [- \infty, + \infty]$ which is the a.e.-limit of $\varphi_n$ (up to a subsequence) . In particular, according to \cite{Giusti, GMM}, we have that $ \varphi_{\infty} $ is a quasi-solution of the minimal surface equation, since $C$ is a minimal cone.  
Therefore the sets 
\beq\label{P-N-non}
P : = \{ x^{'} \in \R^{N-1} \, : \, \varphi_{\infty}(x^{'}) = + \infty  \}, \qquad N : = \{ x^{'} \in \R^{N-1} \, : \, \varphi_{\infty}(x^{'}) = - \infty  \}
\eeq
are both minimal cones in $\R^{N-1}$ with vertex at the origin. 

\smallskip

We also recall that $0 \in \partial C,$ by item i) of Theorem \ref{Teo-minloc}.

\smallskip

Let us prove that $ P = N = \emptyset $. If not, up to replace $u=u(x)$ by $ -u(x^{'},- x_N)$ if necessary, we can suppose that $ N \not\equiv \emptyset$ (i.e., $N$ is of positive measure). In this case $N$ must be rotationally-invariant  as a set of $\R^{N-1}$, indeed by \eqref{ins-livello-risc-sim} and by the fact that $ \varphi_{\infty} $ is the a.e.-limit of $\varphi_n$ (up to a subsequence), we have that 
\beq
 x^{'} \in N \quad \Longrightarrow  \quad \varphi_{\infty}(y^{'}) = -\infty, \,\, \quad \forall \, y^{'} \in \R^{N-1} \quad : \quad \vert y^{'} \vert = \vert x^{'} \vert
\eeq  

Therefore, $N$ is a rotationally-invariant ($N$-1)-dimensional minimal cone and so 
$N = \R^{N-1}$. But this is a contradiction, since $ \R^{N-1} \times \R = N \times \R \subset C $ and $ C \not \equiv \R^N$ (in view of the fact that $0 \in \partial C$).  

\medskip

So far we have proved that $ P = N = \emptyset$, hence $\varphi_{\infty}$ is locally  bounded and therefore it is a classical solution of the minimal surface equation. Then, the limit minimal cone $C$ satisfies 
\beq
C + e_N \subset C \quad {\text and } \quad  C + e_N \not \equiv C 
\eeq
and it must be a half-space, thanks to Lemma 2.3 of \cite{GMM}. This information enables us to use Theorem \ref{Teo-minloc} (item ii)) and Theorem \ref{1D-minloc} to deduce that $u$ is necessarily a 1D monotone profile. Finally, $u$ depends only on the variable $x_N$, due to the fact that it is 1D and it also satisfies $u(x) =u(\vert x^{'} \vert, x_N) $ for every $ x \in \R^N$. 

\end{proof}

\section{Solutions with controlled level sets}\label{sec : Soluzioni-controllate}

This section is devoted to prove some rigidity and symmetry results for entire solutions of the Allen-Cahn equation having controlled level sets. 

The point of view that we adopt here is the following : {\it rigidity and symmetry properties of entire solutions $u$ to the Allen-Cahn equation are a consequence of the "shape" of their set of zeros.} 

In this respect, let us recall that in our previous works \cite{FarGibbons, Fa} we proved the following 

\medskip

\begin{theorem}\label{Teo-superGibbons}
Assume $ N\ge 2$. Let $u$ be an entire solution of the Allen-Cahn equation \eqref{eqAC} and suppose that the set of zeros of $u$ is bounded with respect to some direction $\nu \in S^{N-1}$.  Then, either $u$ is constant ($u \equiv 1$ or $u \equiv -1$) or $u$ is 1D, that is, $u$ is of the form 

$$u(x) = \pm \tanh\Big (\frac{(\nu\cdot x) +\alpha}{\sqrt 2} \Big) \quad \qquad \forall \, x \in \R^N,$$
for some $\alpha \in \R$. 

\end{theorem}

\medskip

The proof of the above result is based on a {\it modified double moving plane procedure} which can be found in \cite{FarGibbons, FarLincei, Fa}.

\medskip

In particular, Theorem \ref{Teo-superGibbons}  proves the following conjecture proposed by G.W. Gibbons \cite{Car, GT} :

\bigskip

\noindent {\bf Gibbons' conjecture.} {\it Assume $ N \ge 2$. Let $u \in C^2(\R^N)$ be a bounded solution of then Allen-Cahn equation such that 
\begin{equation}
\lim_{x_n\rightarrow-\infty} u(x',x_N)=-1
\qquad {\mbox{ and }} \qquad
\lim_{x_n\rightarrow+\infty} u(x',x_N)=1\end{equation}
uniformly w.r.t. $ x' \in \R^{N-1}$. 

\smallskip
\noindent Then 
\[
\quad u(x) = \tanh \Big ( \frac{x_N + \alpha }{\sqrt 2}  \Big).
\]
for some $\alpha \in \R$. 
}

\medskip

We also refer to \cite{BBG, BHM} for a proof of Gibbons' conjecture with different methods.

\bigskip

In the first result of the present section we continue to adopt the above mentioned point of view. Here we are concerned with solutions whose zero level set is given by the graph of a uniformly continuous function $\varphi \in \mathcal{C}^0(\R^{N-1},\R)$, i.e., $ \, \{ u = 0 \} \, \equiv \{(x',x_N) \in \R^{N-1} \times \R \, : \, x_N =\varphi (x') \}$). 

{\it We stress on the fact that no other assumptions are imposed on the solution,  neither monotonicity, nor the regularity of the set of zeros, nor the finiteness of the energy of the considered solution. }

In order to state this result, let us recall that a domain $\Omega \subset \R^N$, $N \ge 2$, is a \emph{continuous epigraph} if, up to a rotation, there exists a function $\varphi \in \mathcal{C}^0(\R^{N-1},\R)$ such that $ \Omega = \{(x',x_N) \in \R^{N-1} \times \R \, : \, x_N > \varphi (x') \}$.

We also recall the notion of modulus of uniform continuity. A function $ \omega : [0, +\infty) \longrightarrow [0, +\infty) $ is a {\it modulus of uniform continuity} if $\omega(0) =0$, $\omega$ is continuous, non-decreasing and subadditive.  

We shall say that a  function $\varphi : \R^{N-1} \longrightarrow \R$ has a modulus of uniform continuity if 
\beq
\forall \, x', y' \in \R^{N-1} \qquad \vert \varphi(x') - \varphi(y') \vert \le \omega(\vert x'-y'\vert) 
\eeq
for some modulus of uniform continuity $\omega$. 

In particular $\varphi$ is a uniformly continuous function on $\R^{N-1}$ with at most linear growth at infinity.

Typical examples of modulus of uniform continuity are provided by the functions $ \omega(t) = t^{\alpha},$ where $\alpha \in (0,1]$.

We also recall that, if $w$ is an entire solution of the Allen-Cahn equation \eqref{eqAC} whose set of zeros is given by the graph of a function  $\varphi \in \mathcal{C}^0(\R^{N-1},\R)$, i.e., $ \, \{ w = 0 \} \, \equiv \{(x',x_N) \in \R^{N-1} \times \R \, : \, x_N =\varphi (x') \}$,  then $w$ must change sign (just apply item iii) of Theorem \ref{Reg-univBound}) and so both $\Omega^{+}_w : = \{ w > 0 \} $ and $\Omega^{-}_w : = \{ w < 0 \} $ are non-empty open connected sets. In particular, up to replace $w$ by $-w$ if necessary, we can always assume that $\Omega^{+}_w = \{ w > 0 \} =  \{(x',x_N) \in \R^{N-1} \times \R \, : \, x_N > \varphi (x') \}$. 

\medskip

We can now state the following

\medskip

\begin{theorem}\label{Teo-grafico=sim}
Assume $ N\ge 2$. Let $u$ be an entire solution of the Allen-Cahn equation \eqref{eqAC} and suppose that the set of zeros of $u$ is given by the graph of some $\varphi \in \mathcal{C}^0(\R^{N-1},\R)$ having a  modulus of uniform continuity.

\noindent Then $u$ is 1D, that is, $u$ is of the form 

$$u=u(x) = \pm \tanh\Big (\frac{(\nu\cdot x) +\alpha}{\sqrt 2} \Big) \quad \qquad \forall \, x \in \R^N,$$
for some $\alpha \in \R$ and some unit vector $ \nu \in \R^N$, with $ \nu_N >0$. 

\end{theorem}


\medskip

The above Theorem recovers and improves upon (with a different method) the following results of

\begin{itemize}
\item L. Modica and S. Mortola \cite{ModMort-DeG} for $N=2$ and where $u$ is a monotone increasing solution whose level lines are the graphs of an equi-lipschitzian family of functions and also satisfying 

\beq
\lim_{x_2\rightarrow-\infty} u(x_1,x_2)=-1
\quad {\mbox{ and }} \quad
\lim_{x_2\rightarrow+\infty} u(x_1,x_2)=1.
\eeq

\item  M.T. Barlow, R.F. Bass and C. Gui \cite{BBG}, for $N\ge 2$ and where $u$ is a monotone increasing solution whose level sets are the graphs of an equi-lipschitzian family of functions and also satisfying 

\beq
\lim_{x_N\rightarrow-\infty} u(x',x_N)=-1
\quad {\mbox{ and }} \quad
\lim_{x_N\rightarrow+\infty} u(x',x_N)=1.
\eeq

\item L. Caffarelli and A. Cordoba \cite{CafCord}, for $N \ge 2$ and where $u$ is a local minimizer  whose zero level set is a globally Lipschitz graph. 

\item O. Savin \cite{S}, for $N \ge 2$ and where $u$ is a monotone increasing solution whose zero level set is a graph with at most linear growth at infinity and also satisfying 

\beq
\lim_{x_N\rightarrow-\infty} u(x',x_N)=-1
\quad {\mbox{ and }} \quad
\lim_{x_N\rightarrow+\infty} u(x',x_N)=1.
\eeq
\end{itemize}

\medskip

To prove the main results of this section we need some preliminary results, which we believe are of independent interest. 


\bigskip

\begin{lemma} \label{comp-infinito} Assume $ N\ge 2$ and let $\Omega \subset \R^N$ be a continuous epigraph.  Let $v \in \mathcal{C}^2(\Omega) \cap \mathcal{C}^0({\overline{\Omega}})$ be a solution of 

\beq\label{pbl-in-Omega-controllato}
\begin{cases}
-\Delta v = v-v^3& \text{in \quad $\Omega$},\\
0 < v \le 1 & \text{in \quad $\Omega$}.
\end{cases}
\eeq
Then 
\beq \label{conva1}
v(y) \longrightarrow 1 \quad uniformly \quad in \quad \Omega \quad whenever \quad dist(y, \partial \Omega) \to +\infty.
\eeq
\end{lemma}

\begin{proof}
This result is essentially contained in the work \cite{BCN3} (see Lemma 3.2 and Lemma 3.3 therein, where only the fact that $\Omega$ is a continuous epigraph is used). For this reason, we regard it as proved.  We also refer to Proposition 7 in \cite{FMV} for a (more) general version of this result on Riemannian manifolds. 
\end{proof}

\bigskip

\begin{lemma} \label{lontano-dagli-estremi} Let $u$ be an entire solution of the Allen-Cahn equation \eqref{eqAC} for $ N \ge 2$ and suppose that the set of zeros of $u$  is given by the graph of some $\varphi \in \mathcal{C}^0(\R^{N-1},\R)$.

\noindent Then, for every $(a,b) \in \R^2$, $a < b$, there exists $ \varepsilon = \varepsilon(a,b) \in (0,1) $ such that

\beq\label{stima-lontano-estremi}
 - 1 + \varepsilon \le u \le 1 - \varepsilon
\eeq
on the set

\beq
\Omega_{a,b} : = \{ \, x \in \R^N \, : \, \varphi(x') + a < x_N < \varphi (x') + b  \, \}.
\eeq
\end{lemma}

\begin{proof} By the assumptions, $u$ is not constant and sign-changing. 
Let $(x^n)_{n \ge1}$ be a sequence of points of $\Omega_{a,b}$ such that $u(x^n) \longrightarrow \sup_{\Omega_{a,b}} u \in [-1 ,1]$ and set $ u_n(x) : = u(x + x^n) $, for every $ x \in \R^N, \, n \ge 1$. By Theorem \ref{Compattezza}, up to a subsequence, we have $ u_n \rightarrow u_{\infty} $ in $ \mathcal{C}^{3}_{loc}(\R^N)$, where $u_{\infty} $ is still a solution of \eqref{eqAC}. Also, for any $ n \ge1$ we consider the point $y^n : = (0^{'}, y_N^n) : = (0^{'}, \varphi((x^n)^{'}) - x_N^n) \in \R^N$. Since $x^n \in \Omega_{a,b}$ we have that $ - b < y_N^n < -a $ and so, up to a subsequence, we may suppose that  $ y_N^n \rightarrow y_N^{\infty} $ in $\R$. 

\noindent We also note that $ u_n(y^n) = u(y^n + x^n) = u((x^n)^{'}, \varphi((x^n)^{'})) = 0 $ for every $ n\ge 1$. 
Therefore, by the above mentioned convergences we obtain that 

\beq
u_n(0) \longrightarrow \sup_{\Omega_{a,b}} u,  
\eeq
\beq
u_n(0) \longrightarrow u_{\infty}(0), 
\eeq
\beq
0= u_n(y^n) = u_n(0^{'}, y_N^n) \longrightarrow u_{\infty}(0^{'}, y_N^{\infty}).
\eeq

Then, by summarizing we have 

\beq
\begin{cases}

-\Delta u_{\infty} =  u_{\infty} - {(u_{\infty})}^3 & \text{in \quad $\R^N$} \\
u_{\infty}(0) = \sup_{\Omega_{a,b}} u \in [-1, 1]   & \\
u_{\infty}(0^{'}, y_N^{\infty}) = 0. &\\
\end{cases}
\eeq
The latter and the strong maximum imply that $\sup_{\Omega_{a,b}} u < 1$. Which provides the inequality on the right-hand side of \eqref{stima-lontano-estremi}. The same argument, applied to a minimizing sequence, yields the other inequality in \eqref{stima-lontano-estremi}. This completes the proof.

\end{proof}

\medskip

We can then prove the following general comparison principle 

\medskip

\begin{proposition} \label{confronto-w-v}
Assume $ N\ge 2$ and let $ \varphi : \R^{N-1} \longrightarrow \R $ be a uniformly continuous function.  Let $\Omega$ be the continuous epigraph $ \{(x',x_N) \in \R^{N-1} \times \R \, : \, x_N > \varphi (x') \}$ bounded by the graph of $ \varphi$. 

Let $w \in \mathcal{C}^2(\Omega) \cap \mathcal{C}^0({\overline{\Omega}})$ such that 

\beq\label{w-sottosol}
\begin{cases}
-\Delta w \le  w-w^3& \text{in \quad $\Omega$},\\
  \quad w \le C_w  & \text{in \quad $\Omega$},
\end{cases}
\eeq
for some constant $ C_w>0$. 

\noindent Let $ v \in \mathcal{C}^2(\Omega) \cap \mathcal{C}^0({\overline {\Omega}})$ satisfying 

\beq\label{soprasol-su-altro-epigr}
-\Delta v \ge  v-v^3 \qquad \text{in \quad $\Omega$}
\eeq
and 
\beq\label{sopra-s_1}
\exists \, b>0 \quad : \quad v \ge s_1 : = \frac{1}{\sqrt 3}\quad on \quad  \{ \, x \in \R^N \, : \,  x_N \ge \varphi (x') + b  \, \}.
\eeq
Then, 

\beq\label{v-domina-w}
v \ge w \quad on \quad  \{ \, x  \, : \,  \varphi (x') \le  x_N \le \varphi (x') + b  \, \} \quad \Longrightarrow  \quad v \ge w \quad  on \quad \Omega
\eeq
 
\end{proposition}

\medskip

\begin{remark}\label{Rk8}

In dimension $N=2$ we do not need to require that $\varphi$ is a uniformly continuous function. The mere continuity of $\varphi$ is enough.  Indeed, in dimension $N=2$, the maximum principle holds for \eqref{PMaxUnb}, as soon as $A$ is a non empty, non dense domain of $\R^2$. 
 
\end{remark}

\medskip

\begin{proof}
Set $ z = w -v \in \mathcal{C}^2(\Omega) \cap \mathcal{C}^0({\overline{\Omega}})$.Then 
$ - \Delta z = - \Delta w + \Delta v \le (w-w^3) - (v-v^3) = c z $ on $ \Omega$, where 
\beq
c(x) = \begin{cases}
\frac{(w-w^3) - (v-v^3) }{w-v} & \text{if \quad $ w(x) \neq v(x)$},\\
0 & \text{otherwise},
\end{cases}
\eeq
is well-defined and belongs to $L^{\infty}_{loc}(\Omega) $. 

For contradiction, suppose that $z>0$ somewhere in $\Omega$ , then the open set $X : = \{ \, x \in \Omega \, : \,  x_N > \varphi (x') + b  \, \} \cap \{ \, z>0 \, \}$ is non-empty and thus we can consider a connected component $A$ of $X$ (which is also open in $\R^N$).
If $ x \in \partial A$ then either $x \in \Omega$ and so $ z(x) =0$, or $x \in  \partial  \{ \, x \in \Omega \, : \,  x_N > \varphi (x') + b  \, \} $ and so $x \in \{ \, x  \, : \,  \varphi (x') \le  x_N \le \varphi (x') + b  \, \} $, which in turn implies that $ z(x) \le 0$.  Thus, $z$ satisfies
\beq
\begin{cases}
- \Delta z \le cz  & \text{in \quad $A$},\\
z \le 0 & \text{on \quad $\partial A $},
\end{cases}
\eeq

Observe that $A \subset \Omega$ yields that $\overline A$ is disjoint from the closure of an infinite open connected cone (since $\Omega$ satisfies this property, in view of the uniform continuity of $\varphi$). 
So, in order to apply the maximum principle in unbounded domains (see for instance Lemma 2.1 in \cite{BCN3}) whose closure is disjoint from the closure of an infinite open connected cone, we need only to check that $z$ is bounded from above on $A$ and $ c \le 0$ on $A$.  

To this end, we observe that $ x \in A $ implies $ w(x) > v(x) \ge s_1 $ (the first inequality being true by definition of $A$ and the second one by \eqref{sopra-s_1}, since $ A\subset X$) and then we deduce that $0 < z = w-v \le C_w + s_1 < +\infty$ on $A$.  Therefore $z$ is bounded on $A$. Moreover, the nonlinear function $ t \longrightarrow t - t^3$ is non-increasing on the interval $[s_1, +\infty)$ and thus $c\le 0$ on $A$. 

By summarizing, $z$ satisfies 
\beq
\begin{cases}\label{PMaxUnb}
z \in \mathcal{C}^2(A) \cap \mathcal{C}^0({\overline{A}}),\\
- \Delta z \le 0 & \text{in \quad $A$},\\
z \le 0 & \text{on \quad $\partial A $},
\end{cases}
\eeq
and the maximum principle then provides $ z \le 0$ on $A$, a contradiction. Therefore $X =\emptyset$ which, together with \eqref{v-domina-w}, immediately gives the desired conclusion.
 
\end{proof}


\bigskip

\begin{theorem} \label{Teo-v-sopra-w}
Assume $ N\ge 2$. Let $w$ be an entire solution of the Allen-Cahn equation \eqref{eqAC} and suppose that the set of zeros of $w$ is given by a graph of some $\varphi \in \mathcal{C}^0(\R^{N-1},\R)$ having a  modulus of uniform continuity. 

Let $ v $  be another solution of the Allen-Cahn equation \eqref{eqAC} satisfying 
\beq\label{positiva-su-altro-epigraf}
v > 0  \qquad \text{in  \quad $\Omega^{+}_w = \{ w > 0 \} $}.
\eeq
Then, 

\beq\label{v-sopra-w}
 v \ge w \qquad \text{in  \quad $\Omega^{+}_w = \{ w > 0 \} $}.
\eeq
\end{theorem}

\begin{proof}
By applying Lemma \ref{comp-infinito} to both $w$ and $v$ on the continuous epigraph $\Omega^{+}_w = \{ w > 0 \} $, we get 
\beq \label{conva1-bis}
w(y), v(y) \longrightarrow 1 \quad  uniformly \quad in \quad \Omega^{+}_w \quad whenever \quad dist(y, \partial \Omega^{+}_w) \to +\infty.
\eeq
Since $\varphi$ has a modulus of uniform continuity, the latter implies the existence of $b=b(w,v) >0$ such that

\beq\label{w-v-sopra-s_1}
v,\, w \, \ge s_1 : = \frac{1}{\sqrt 3}\qquad on \qquad  \{ \, x \in \R^N \, : \,  x_N \ge \varphi (x') + b  \, \}
\eeq
and so in view of Proposition \ref{confronto-w-v}, to demonstrate the claim it is enough to prove that $ v \ge w$ on the set $ \{ \, x  \, : \,  \varphi (x') \le  x_N \le \varphi (x') + b  \, \} = {\overline { \Omega_{0,b}}}$. 

\medskip

For all $t >0$ we consider the translated function $ v_t (x) := v(x + t e_N)$, where $e_N$ is the last vector of the canonical basis of $\R^N$. Then $v_t$ satisfies 

\beq
\begin{cases}
-\Delta v_t = v_t-(v_t)^3& \text{in \quad $\R^N$},\\
 - 1 \le  v_t \le  1 & \text{in \quad $\R^N$},\\
0 < v_t \le 1 & \text{in \quad $\Omega^{+}_w$}
\end{cases}
\eeq
and 
\beq
v_t \ge s_1 : = \frac{1}{\sqrt 3}\qquad on \qquad  \{ \, x \in \R^N \, : \,  x_N \ge \varphi (x') + b  \, \},
\eeq
 since $v$ is an entire solution of the Allen-Cahn equation satisfying \eqref{positiva-su-altro-epigraf} ,\eqref{w-v-sopra-s_1} and $t>0$. 

By Lemma \ref{lontano-dagli-estremi} there exists $\varepsilon = \varepsilon(b) \in (0,1) $ such that $ 0 \le  w \le 1- \varepsilon $ on the set $ \{ \, x  \, : \,  \varphi (x') \le  x_N \le \varphi (x') + b  \, \} = {\overline { \Omega_{0,b}}}$ therefore, by making use of the latter and \eqref{conva1-bis}  we can find $ t_0 > b $ such that

\beq\label{t-grandi}
\forall \, t \ge t_0 \qquad v_t \ge w \qquad \text{on \quad $ {\overline { \Omega_{0,b}}}$}.
\eeq
Let 
\beq\label{defT}
T : = \inf \, \{ \, t_0 > 0 \, : \, \forall \, t \ge t_0 \quad v_t \ge w \quad \text{on \quad $ {\overline { \Omega_{0,b}}}$} \, \} 
\eeq
then $ T \in [0, + \infty)$ by \eqref{t-grandi} and to conclude the proof it is enough to prove that $ T=0$.  To this end we first observe that, by Proposition \ref{confronto-w-v} applied to $w$ and $v_t$, we have 

\beq\label{vT-sopra-w}
\forall \, t \ge T \qquad v_t \ge w \qquad \text{on \quad $\Omega^{+}_w$}
\eeq
and then we argue by contradiction. Suppose $ T>0 $, then there are a sequence $(x^n)_{n \ge1}$ of points of $ {\overline { \Omega_{0,b}}}$ and an increasing sequence of real numbers $ (t_n)_{n \ge 1} $ such that 

\beq\label{cons-ipot-contrad}
\forall \, n \ge 1 \qquad  v_{t_n} (x^n) < w (x^n), \qquad 0 < t_n < T \quad \text {and} \quad t_n \longrightarrow T. 
\eeq
For every $n \ge 1$ we set $ w_n(x) : = w(x + x^n) $, $ v_n(x) : = v(x + x^n) $ and $ \varphi_n(x^{'}) := \varphi(x^{'} + (x^n)^{'}) - x_N^n $. Then by Theorem \ref{Compattezza}, up to subsequences, we have $ w_n \rightarrow w_{\infty} $ and $ v_n \rightarrow v_{\infty} $ in $ \mathcal{C}^{3}_{loc}(\R^N)$, where $w_{\infty},v_{\infty} $ are still entire solutions of the Allen-Cahn equation \eqref{eqAC}. Moreover, the sequence $(\varphi_n)_{n \ge1}$ satisfies

\beq\label{stima-punt}
\forall \, n \ge 1 \qquad -b \le \varphi_n(0^{'}) \le 0,  
\eeq
and for every $ n \ge 1$ and for every $ x^{'}, \, y^{'} \in \R^{N-1}$
\beq\begin{split}
\vert \varphi_n (x^{'}) - \varphi_n (y^{'}) \vert = \vert \varphi(x^{'} + (x^n)^{'}) - x_N^n  - \varphi(y^{'} + (x^n)^{'}) + x_N^n\vert = \\ \vert \varphi(x^{'} + (x^n)^{'})  - \varphi(y^{'} + (x^n)^{'}) \vert \le \omega( \vert x^{'} + (x^n)^{'} - y^{'} - (x^n)^{'}\vert) = &  \, \omega( \vert x^{'} - y^{'} \vert). 
\end{split}\eeq
Therefore, by Ascoli-Arzel\`a Theorem there exists $\varphi_{\infty} \in \mathcal{C}^0(\R^{N-1},\R)$ such that, up to a subsequence, $ \varphi_n \rightarrow \varphi_{\infty} $ in $ \mathcal{C}^{0}_{loc}(\R^{N-1})$, also $\varphi_{\infty}$ has modulus of uniform continuity given by $\omega$. Those properties also imply that $w_{\infty}$ vanishes on the graph of $\varphi_{\infty}$ since 
\beq\begin{split}
\vert w_n(x^{'}, \varphi_{\infty}(x^{'} ) \vert = \vert w_n(x^{'}, \varphi_{\infty}(x^{'} ) - w_n(x^{'}, \varphi_{n}(x^{'} ) + w_n(x^{'}, \varphi_{n}(x^{'} ) \vert = \\  \vert w_n(x^{'}, \varphi_{\infty}(x^{'} ) - w_n(x^{'}, \varphi_{n}(x^{'} ) + w(x^{'} + (x^n)^{'}, \varphi(x^{'} + (x^n)^{'}) \vert = \\ \vert w_n(x^{'}, \varphi_{\infty}(x^{'} ) - w_n(x^{'}, \varphi_{n}(x^{'} ) + 0 \vert \le \Vert \nabla w_n \Vert_{\infty} \Vert \big( x^{'}, \varphi_{\infty}(x^{'} ) \big) - \big( x^{'}, \varphi_{n}(x^{'} ) \big) \Vert_{\R^N} =\\ \Vert \nabla w \Vert_{\infty} \vert \varphi_{\infty}(x^{'})  - \varphi_{n}(x^{'}) \vert \longrightarrow 0 \quad \text{as} \quad n \to \infty \hskip3 true cm
\end{split}\eeq
and thus
\beq
\forall \, x^{'} \in \R^{N-1} \qquad \vert w_{\infty}(x^{'}, \varphi_{\infty}(x^{'} ) \vert =\lim_{n\to +\infty} \vert w_n(x^{'}, \varphi_{\infty}(x^{'} ) \vert = 0. 
\eeq
Also from \eqref{stima-punt} we get 
\beq\label{stima-punt2}
-b \le \varphi_{\infty}(0^{'}) \le 0.   
\eeq
Next we prove that

\beq\label{pos-prof-infinito}
w_{\infty} > 0 \quad \text {on the continuous epigraph} \quad   \{ \, x \in \R^N \, : \,  x_N > \varphi_{\infty} (x')  \, \} : = \Omega_{\infty}.
\eeq
Indeed, by $ \varphi_n \rightarrow \varphi_{\infty} $ in $ \mathcal{C}^{0}_{loc}(\R^{N-1})$, we have that for any compact set $K \subset \Omega_{\infty}$, there is an integer $n(K) \ge 1$ such that $K \subset \Omega^{+}_{w_n} = \{ w_n > 0 \} $, for every $ n \ge n(K) $.  Then, $ w_n \rightarrow w_{\infty} $ in $ \mathcal{C}^{3}_{loc}(\R^N)$ implies that $w_{\infty} \ge 0$ on $\Omega_{\infty}$. On the other hand, by definition of $ w_n, \varphi_n$ of  $ w_{\infty}, \varphi_{\infty}$ and by \eqref{w-v-sopra-s_1} we also have 

\beq\label{winfinto-sopra-s_1}
w_{\infty} \, \ge s_1 : = \frac{1}{\sqrt 3}\qquad on \qquad  \{ \, x \in \R^N \, : \,  x_N \ge \varphi_{\infty} (x') + b  \, \}
\eeq
and by summarizing we see that $w_{\infty}$ satisfies
\beq\label{prop-w-infinito}
\begin{cases}
-\Delta w_{\infty} = w_{\infty}-(w_{\infty})^3& \text{in \quad $\R^N$},\\
w_{\infty} \ge 0 & \text{on \quad $\Omega_{\infty}$},\\
w_{\infty} = 0 & \text{on \quad $\partial \Omega_{\infty}$}, \\
w_{\infty} \, \ge s_1 : = \frac{1}{\sqrt 3}& \text{on $ \{ \, x \in \R^N \, : \,  x_N \ge \varphi_{\infty} (x') + b  \, \}  \subset \Omega_{\infty} $}.
\end{cases}
\eeq
The strong maximum principle then implies \eqref{pos-prof-infinito}. 

Analogously we find that
\beq\label{prop-v-infinito}
\begin{cases}
-\Delta v_{\infty} = v_{\infty}-(v_{\infty})^3& \text{in \quad $\R^N$},\\
v_{\infty} \ge 0 & \text{on \quad $\Omega_{\infty}$},\\
v_{\infty} \, \ge s_1 : = \frac{1}{\sqrt 3}& \text{on $ \{ \, x \in \R^N \, : \,  x_N \ge \varphi_{\infty} (x') + b  \, \}  \subset \Omega_{\infty} $}
\end{cases}
\eeq
and so $v_{\infty} >0$ on $\Omega_{\infty}$ again by the strong maximum principle. 

Moreover, combining \eqref{stima-punt2} and the definition of $\Omega_{\infty}$ we get :
\beq\label{pos.origine}
\text {the origin} \quad  0 = (0^{'}, 0) \in \R^N  \quad \text {belongs to} \quad {\overline {\Omega_{\infty}}}
\eeq
which in turn implies that
\beq\label{pos.origine2}
\text {the point} \quad  (0^{'}, T) \in \R^N  \quad \text {belongs to} \quad{\Omega_{\infty}}
\eeq
since $ T>0$ and $\Omega_{\infty}$ is an epigraph!

On the other hand, \eqref{vT-sopra-w} provides
\beq\label{vT-inf-sopra-w-inf}
(v_{\infty})_T \ge w_{\infty} \qquad \text{on \quad $\Omega_{\infty}$}
\eeq
 and \eqref{cons-ipot-contrad} yields

\beq\label{cons-ipot-contrad-inf}
(v_{\infty})_T (0) = \lim_{n\to+\infty} v_{t_n} (x^n) \le \lim_{n \to +\infty} w (x^n) = w_{\infty}(0).
\eeq 
But  $(v_{\infty})_T (0) = v_{\infty} (0^{'}, T) >0$, by \eqref{pos.origine2} and $v_{\infty} >0$ on $\Omega_{\infty}$, thus also $w_{\infty}(0)>0$ and so also
\beq\label{pos.origine3}
\text {the origin} \quad  0 = (0^{'}, 0) \in \R^N  \quad \text {belongs to} \quad{\Omega_{\infty}}.
\eeq

To achieve the desired contradiction we consider the smooth function $\theta := w_{\infty} - (v_{\infty})_T $ which fulfills, thanks to \eqref{vT-inf-sopra-w-inf} and \eqref{cons-ipot-contrad-inf}, 
\beq
\begin{cases}
-\Delta \theta  = c \theta & \text{in \quad $\Omega_{\infty}$},\\
\theta \le 0 & \text{on \quad $\Omega_{\infty}$},\\
\theta (0) = 0 &
\end{cases}
\eeq
where 
\beq
c(x) = \begin{cases}
\frac{(w_{\infty}-(w_{\infty})^3) - ((v_{\infty})_T-((v_{\infty})_T)^3) }{w_{\infty}-(v_{\infty})_T} & \text{if \quad $ w_{\infty}(x) \neq (v_{\infty})_T(x)$},\\
0 & \text{otherwise},
\end{cases}
\eeq
is well-defined and belongs to $L^{\infty}(\Omega_{\infty}) $. 
In view of \eqref{pos.origine3}, the strong maximum principle gives $ \theta \equiv 0$ on ${\overline {\Omega_{\infty}}}$, that is $ w_{\infty} \equiv  (v_{\infty})_T$ on ${\overline {\Omega_{\infty}}}$. 

But for $x \in \partial \Omega_{\infty}$ we have 
\beq\label{fin-contrad}
0 = w_{\infty}(x) =  (v_{\infty})_T(x) =  (v_{\infty})_T(x^{'}, \varphi_{\infty}(x^{'})) = v_{\infty}(x^{'}, \varphi_{\infty}(x^{'})  + T) >0
\eeq
since $ T>0 $ and so the point $(x^{'}, \varphi_{\infty}(x^{'})+ T)$ must belong to $\Omega_{\infty}.$ \eqref{fin-contrad} is the desired contradiction, which concludes the proof of the Theorem.

\end{proof}

\bigskip

We are now ready to prove Theorem \ref{Teo-grafico=sim}. 

\medskip

\begin{proof} [Proof of Theorem \ref{Teo-grafico=sim}]  
As discussed at the beginning of the present section we can assume that $\Omega^{+}_u := \{ u > 0 \} =  \{(x',x_N) \in \R^{N-1} \times \R \, : \, x_N > \varphi (x') \}$, for some $\varphi \in \mathcal{C}^0(\R^{N-1},\R)$ having a  modulus of uniform continuity. 

Now we prove that 

\beq\label{mezza-monot}
\frac{\partial u}{\partial x_N} (x) > 0 \qquad \quad \forall \, x \in \Omega^{+}_u.  
\eeq
To see this, it is enough to apply Theorem \ref{Teo-v-sopra-w} with $ w = u$, $ v = u_t $ (the translated of $u$ with $t>0$) and $\Omega =  \Omega^{+}_u. $ Indeed, proceeding in this way we obtain that $ u_t \ge u $ on $\Omega^{+}_u$, for every $t >0$. Thus $\frac{\partial u}{\partial x_N} \ge 0$ on $ \Omega^{+}_u$ and the strong maximum principle applied to the non-negative function $\frac{\partial u}{\partial x_N}$ provides \eqref{mezza-monot}. 
Now we set $ w (x) = - u(x',-x_N)$. $w$ is still an entire solution of the Allen-Cahn equation whose set of zeros is given by a function having a modulus of uniform continuity. Therefore we can apply the above argument to $w$ to deduce that $\frac{\partial w}{\partial x_N} > 0 $ on $\Omega^{+}_w.$ By definition of $w$ this means that $\frac{\partial u}{\partial x_N} > 0$ on $ \Omega^{-}_u$. The latter and \eqref{mezza-monot} yield
\beq
\frac{\partial u}{\partial x_N} (x) > 0 \qquad \quad \forall \, x \in \{ u \neq 0 \}.
\eeq
The latter and the strong maximum principle then imply that $u$ is monotone in the variable $x_N$, i.e., 
\beq\label{tutta-monot}
\frac{\partial u}{\partial x_N} (x) > 0 \qquad \quad \forall \, x \in \R^N.
\eeq
On the other hand, we already know (by Lemma \ref{comp-infinito}) that 
\begin{equation}
\lim_{x_N\rightarrow-\infty} u(x',x_n)=-1
\quad {\mbox{ and }} \quad
\lim_{x_N\rightarrow+\infty} u(x',x_n)=1,
\end{equation}
therefore $u$ is a non-constant local minimizer and to conclude we can argue as in the proof of Theorem \ref{monot-cono-cilindro}.  We consider the subgraphs 

\beq
\{ \,  u < 0 \, \} = \{ (x^{'},x_N) \in \R^{N-1} \times \R \, : \, x_N < \varphi(x^{'}) \, \} 
\eeq
and 
\beq\label{sottografici-risc-non}
\{ \,  u_{\eps_n} < 0 \, \} = \{ (x^{'},x_N) \in \R^{N-1} \times \R \, : \, x_N < \varphi_n(x^{'}) \, \} : = U_n
\eeq
where 
\beq
\varphi_n (x^{'}) : = \eps_n  \varphi \Big (\frac{x^{'}}{\eps_n} \Big).  
\eeq

By item i) of Theorem \ref{Teo-minloc} and Lebesgue dominated convergence theorem we immediately have that

\beq\label{cv-sottografici-non}
 {\bf 1}_{U_n} \longrightarrow {\bf 1}_{\R^N \setminus C} \quad  {\text {in}} \quad L^1_{loc}(\R^N).
\eeq
The latter and Lemma 16.3 of \cite{Giusti} tell us that $\R^N \setminus C$ is itself a subgraph of a measurable function 
$ \varphi_{\infty} : \R^{N-1} \rightarrow [- \infty, + \infty]$ which is the a.e.-limit of $\varphi_n$ (up to a subsequence) . In particular, according to \cite{Giusti, GMM}, we have that $ \varphi_{\infty} $ is a quasi-solution of the minimal surface equation, since $C$ is a minimal cone.  On the other hand, $\varphi$ has at most linear growth at infinity (since it is uniformly continuous) and so, also the rescaled functions $\varphi_n$ share this property. This yields that the limit function $\varphi_{\infty}$ is locally  bounded and therefore it is a classical solution of the minimal surface equation. 
Now, by proceeding exactly as at the end of the proof of Theorem \ref{monot-cono-cilindro} we can prove  that $C$ is a half-space. Therefore, Theorem \ref{1D-minloc} tell us that $u$ is necessarily a 1D monotone profile. 

\end{proof}

\appendix

\section{Complete classification in the one-dimensional case}\label{sec : appendix A}

\medskip

In this section we classify all the entire solutions of the Allen-Cahn equation (not necessarily of finite energy) in the one-dimensional case .

\medskip

\begin{theorem}\label{Teo-clas1D}
Let $u \in \mathcal{C}^2(\R)$ be a  solution of the Allen-Cahn equation 

\beq\label{eqAC1D}
- u^{''} = u - u^3  \qquad \text{in} \qquad \R.
\eeq
Then one of the following alternatives occurs : 

\smallskip

\item[i)] $u$ is constant (and necessarily $ u \equiv 0, \pm1$).  

\item[ii)] $u = u(s) = \pm \tanh \Big ( \frac{s +\alpha}{\sqrt{2}} \Big )$, for every $s \in \R$ and some $\alpha \in \R$. 

\item[iii)] $u$ is a non-constant periodic solution. 

\medskip

\noindent Furthermore, if $u$ has finite energy then one of the following alternatives occurs : 

\smallskip

\item[a)] $u \equiv \pm 1$ (and in this case $ L_u =0$).  

\item[b)] $u = u(s) = \pm \tanh \Big ( \frac{s +\alpha}{\sqrt{2}} \Big )$, for every $s \in \R$ and some $\alpha \in \R$. In this case $L_u = \frac{2\sqrt{2}}{3}$.

\end{theorem}

\medskip

\begin{proof}
The above result follows from a standard ode analysis. We give here a short proof for the reader's convenience. 

The only constant solutions of \eqref{eqAC1D} are $-1, 0$ and $1$. Let us suppose that $u$ is not constant, then $ \vert u \vert < 1 $ everywhere on $ \R $ and thus $u$ must change sign on $\R$ (by a standard convexity argument). 
If $u$ has no critical points then, up to replace $u$ by $-u$, we can suppose that $u^{'}>0$ on $\R$. Thus, either $\lim_{s\to -\infty}u(s) = -1 $ or $\lim_{s\to +\infty} u(s)= 1 $. Multiplying \eqref{eqAC1D} by $u^{'}$ and integrating we get $ \frac{(u^{'})^2}{2} = \frac{(1 - u^2)^2}{4} + c$ on $ \R$, where $c$ is a constant. The behaviour at infinity of $u$ then yields $c=0$ and so $u$ solves $ u^{'} = \frac{1 - u^2}{\sqrt 2} $ on $ \R$. The Cauchy uniqueness theorem then implies that $ u(s) = \tanh \Big ( \frac{s +\alpha}{\sqrt{2}} \Big )$, with $ \alpha =  \sqrt{2}\tanh^{-1}(u(0))$. To continue, we observe that if $u$ has a critical point then $u$ must be symmetric (even) with respect to this point, again by Cauchy uniqueness theorem. Thus, either $u$ has more the one critical point and then it is periodic (non constant) or $u$ has exactly one critical point $s_0$, it is even with respect to $s_0$ and it is strictly monotone for $s \not = s_0$. In particular, $s_0$ is either a minimum point or a maximum point and, once again up to replace $u$ by $-u$, we can suppose that $s_0$ is a minimum point for $u$. Therefore, equation \eqref{eqAC1D} gives $u(s_0) < 0 $ and so $\lim_{s\to +\infty} u(s)= 1 $, since $u$ must change sign, and thus $u$ solves $ \frac{(u^{'})^2}{2} = \frac{(1 - u^2)^2}{4} $ on $ \R$. The latter yields $ 0 = \frac{(u^{'}(s_0))^2}{2} = \frac{(1 - u^2(s_0))^2}{4}$, that is $ u(s_0) = - 1$. Then, the strong minimum principle gives $ u \equiv -1$, which contradicts the fact that $u$ was supposed to be non constant. This concludes the proof of the first part of the theorem.  

Observe that $u$ has finite energy if and only if $ \frac{(u^{'})^2}{2} + \frac{(1 - u^2)^2}{4} \in L^{1}(\R)$.  This immediately excludes $u \equiv 0$ and every non constant periodic solution. Since $ u \equiv \pm 1$ clearly satisfies $ L_u= 0$, it remains to treat the case $u(s) = \tanh \Big ( \frac{s}{\sqrt{2}} \Big )$. Since $ u^{'} = \frac{1 - u^2}{\sqrt 2} $ on $ \R$ we have 

\beq\begin{split}
\int_{\R}  \Big [ \frac{(u^{'}(s))^2}{2} + \frac{(1 - u^2(s))^2}{4} \Big ] ds & = \int_{\R}  (u^{'}(s))^2 ds = \int_{\R}  \frac{(1 - u^2(s))}{\sqrt 2} u^{'}(s) ds \\
& = \int_{-1}^1  \frac{(1 - t^2)}{\sqrt 2} dt = \frac{2\sqrt{2}}{3}. 
\end{split} \eeq

\end{proof}

Next we prove that every one-dimensional monotone profile is a local minimizer in any dimension and we calculate its energy.   

\medskip

\begin{theorem} \label{Energia1D} 
Assume $ N \ge 1$.  Let $\alpha \in \R$ and $ \nu \in \R^N$ be a unit vector. Every one-dimensional monotone profile  $v=v(x) = \pm \tanh\big (\frac{(\nu\cdot x) +\alpha}{\sqrt 2} \big)$ is a local minimizer and $ L_v = \frac{2\sqrt{2}}{3} \omega_{N-1}$ (where we have set $ \, \omega_0 = 1$). 

\end{theorem}{

\begin{proof}
The local minimality is an immediate consequence of the following minimality property demostrated by 
G. Carbou \cite{Car} (see also \cite{Fa3, AlbRiar}). The function $u = u(s) = \tanh\big (\frac{s}{\sqrt 2} \big)$ is a minimizer of the functional  $\int_{\R}  \Big [ \frac{(v^{'})^2}{2} + \frac{(1 - v^2)^2}{4} \Big ] ds$ in 

\beq
{\cal K} : =  \{ v \in H^{1}_{loc}(\R) \, : \, \lim_{s\to-\infty} v(s) = -1, \quad \lim_{s\to+\infty} v(s)= 1 \}
\eeq
and that $u=u(x_1,...,x_N) = \tanh\big (\frac{x_1}{\sqrt 2} \big)$ is a minimizer of the functional

\beq
\int_{\R \times \omega} \frac{\vert \nabla v\vert^2}{2} + \frac{(1-v^2)^2}{4}
\eeq
in
\beq
{\cal K}_{\omega} : =  \{ v \in H^{1}_{loc}(\R^N) \cap L^{\infty}(\R^N)\, : \lim_{x_1 \to\pm\infty} v= \pm 1 \, {\text {uniformly in}} \,\,  x_2,...,x_N  \in \omega\}
\eeq
for every open and bounded set $ \omega \subset \R^{N-1}$, $ N\ge2$. 


\medskip

For $N=1$ the energy has already been calculated in Theorem \ref{Teo-clas1D}, so we can consider $ N \ge 2$.  
In view of item iv) of Theorem \ref{Energia} it is enough to calculate $ L_v $ for $v=v(x) = \tanh\big (\frac{x_N }{\sqrt 2} \big)$. To this end, set $h = h(s) : = \tanh( \frac{s}{\sqrt{2}})$ and consider $ R>0$, then 
\beq\begin{split}
&\int_{B^N(0,R)} \frac{\vert \nabla v\vert^2}{2} + \frac{(1-v^2)^2}{4} dx = \int_{-R}^R \int_{B^{N-1}\Big (0^{'}, {\sqrt{R^2 - x_N^2}} \Big )} \frac{\vert \nabla v\vert^2}{2} + \frac{(1-v^2)^2}{4} dx^{'}dx_N \\
& = \int_{-R}^R \omega_{N-1} (R^2 - x_N^2)^{\frac{N-1}{2}} \Big (\frac{\vert h^{'}(x_N) \vert^2}{2} + \frac{(1-h^2(x_N))^2}{4} \Big) dx_N,
\end{split}\eeq
hence
\beq\begin{split}
 \frac{1}{R^{N-1}} \int_{B^N(0,R)} \frac{\vert \nabla v\vert^2}{2} + \frac{(1-v^2)^2}{4} dx = \\
\omega_{N-1} \int_{\R}  \Big (1 - \Big (\frac{x_N}{R} \Big )^2 \Big )^{\frac{N-1}{2}} \Big (\frac{\vert h^{'}(x_N) \vert^2}{2} + \frac{(1-h^2(x_N))^2}{4} \Big) {\bf 1}_{(-R,R)}dx_N,
\end{split}\eeq
so that 
\beq\begin{split}
 \frac{1}{R^{N-1}} \int_{B^N(0,R)} \frac{\vert \nabla v\vert^2}{2} + \frac{(1-v^2)^2}{4} dx \longrightarrow 
 \omega_{N-1} \int_{\R} \frac{\vert h^{'}\vert^2}{2} + \frac{(1-h^2)^2}{4} = \frac{2\sqrt{2}}{3} \omega_{N-1}
\end{split}\eeq
by Lebesgue's dominated convergence. 

\end{proof}

We conclude this section with the stability analysis of every solution to the ode \eqref{eqAC1D}.

\medskip

\begin{theorem} \label{stab-1D} 
The functions $u \equiv \pm 1$ and $u(s) = \pm \tanh \Big ( \frac{s +\alpha}{\sqrt{2}} \Big )$ are stable, while $u\equiv 0$ and every periodic solution are unstable outside every compact set of $ \R$.  
\end{theorem}

\begin{proof}
The functions $u \equiv \pm 1$ and $u(s) = \pm \tanh \Big ( \frac{s +\alpha}{\sqrt{2}} \Big )$ are local minimizers of the functional  $\int_{\R}  \Big [ \frac{(v^{'})^2}{2} + \frac{(1 - v^2)^2}{4} \Big ] ds$, thus they are stable. 
To study the instability of the others solutions we consider the second variation of the energy functional  $\int_{\R}  \Big [ \frac{(v^{'})^2}{2} + \frac{(1 - v^2)^2}{4} \Big ] ds$ at $u$, that is, the quadratic form 

\beq\label{SecVar-1D}
{\cal Q}_u (\psi) : = \int_{\R}\vert \psi^{'} \vert^2 - (1-3u^2) \psi^2 \qquad \forall \psi \in \mathcal{C}^1_c(\R).
\eeq
For $u \equiv 0$, \eqref{SecVar-1D} is   
\beq
{\cal Q}_0 (\psi) : = \int_{\R}\vert \psi^{'} \vert^2 - \psi^2 \qquad \forall \psi \in \mathcal{C}^1_c(\R),
\eeq  
for which it is well-known the existence of $\phi \in \mathcal{C}^1_c(\R)$ such that $ {\cal Q}_0 (\phi) <0$. The translation invariance of ${\cal Q}_0 $ immediately implies the instability of $u \equiv 0$ outside every compact set of $\R$. 

Let $u$ be any periodic (non-constant) solution to the ode \eqref{eqAC1D} and let $M_u $ be the unbounded and countable set of its maxima. Let $ m_1, m_2 \in M_u $ with $ m_1< m_2$, then $v : = u^{'}$ solves the linearized equation $ - v^{''} = (1-3u^2) v $ and an integration by parts, after multiplying the latter equation by $v$,  gives 

\beq\label{ener-lin0}
\int_{m_1}^{m_2} \vert v^{'} \vert^2 - (1-3u^2) v^2 = 0.
\eeq 
We claim that the latter implies the existence of $\phi \in H^1_0((m_1, m_2))$ such that $ \int_{m_1}^{m_2} \vert \phi^{'} \vert^2 - (1-3u^2) \phi^2 <0$, which immediately gives the desired conclusion (since the set $M_u $ is unbounded and countable and $\mathcal{C}^1_c((m_1, m_2))$ is dense in $H^1_0((m_1, m_2))$). 

Suppose, for contradiction that  
\beq\label{ener-lin}
\int_{m_1}^{m_2} \vert \psi^{'} \vert^2 - (1-3u^2) \psi^2 \ge 0  \qquad \forall \psi \in H^1_0((m_1, m_2))
\eeq
then, thanks to \eqref{ener-lin0}, the function $ \phi : = \vert u^{'} \vert \in H^1_0((m_1, m_2))$ also minimizes the functional in \eqref{ener-lin} on $H^1_0((m_1, m_2)).$ Therefore, $ \phi : = \vert u^{'} \vert$ is a {\it non-negative} solution of the linearized equation $ - v^{''} = (1-3u^2) v $ on the interval $(m_1, m_2)$, which vanishes somewhere in $(m_1, m_2)$ (since between two maxima of a periodic non-constant solution $u$ there is a minimum point of $u$). 
The strong maximum principle then yields $  \phi : = \vert u^{'} \vert\equiv 0$ on the interval $(m_1, m_2)$. A contradiction. 

\end{proof}

\section{Radial solutions}\label{sec : appendix B}

\medskip

In this section we are concerned with entire radial solutions of the Allen-Cahn \eqref{eqAC}, i.e., solutions satisfying 
$ u(x) = u (\vert x \vert )$  for every $ x \in \R^N$, $ N \ge2$. 

\medskip

\begin{theorem} \label{Solradiali} 
Assume $ N \ge 2$ and let $u$ be a radial solution to the Allen-Cahn \eqref{eqAC}. Then

\item i) $u$ has finite energy if and only if $ u \equiv \pm 1$.

\item ii) if $u$ is non constant, then $u$ is unstable outside every compact set.

\end{theorem}{

\begin{proof}
Apply Corollary \ref{nonEsist-radiali1} to obtain item i). 

To prove item ii) we follow the idea used to study stability for solutions of the Allen-Cahn equation in the one-dimensional case (see Theorem \ref{stab-1D}). Since $u$ is a non-constant radial solution of \eqref{eqAC}, we know from item iii) of Theorem \ref{Reg-univBound} that $u$ must be an oscillating solution, that is $u$ must vanish on infinitely many spheres centered at the origin, whose radii go to $+\infty$. Therefore also the radial derivative $u$ must vanish on infinitely many spheres centered at the origin, whose radii go to $+\infty$. Let us denote by $v$ the radial derivative of $u$ and let $\{ r_n\}_{n \ge1}$ be a strictly increasing sequence of positive radii such that $v = 0$ on the sphere centered at the origin and of radius $r_n$. 

Differentiating the equation \eqref{eqAC} with respect to the radial direction and recalling that $u$ is a radial function we immediately get that $v$ solves 

\beq\label{eq-lin-rad}
 -  \Delta v + (N-1) \frac{v}{\vert x \vert^2} = (1-3u^2) v \qquad {\rm on} \qquad \R^N \setminus \{0\}. 
\eeq
Multiply the equation \eqref{eq-lin-rad} by $v$ and integrate by parts on the annulus $A_n : = B(0,r_{n+1}) \setminus {\overline {B(0,r_{n})}}$, by definition of $v$ and $\{ r_n\}_{n \ge1}$, we get

\beq
\int_{A_n} \vert \nabla v \vert^2 - (1-3u^2)v^2 \le - \int_{A_n} (N-1)\vert \frac{v^2}{\vert x \vert^2} < 0,
\eeq 
since $v$ is not identically 0 on the annulus $A_n$ and $ N \ge 2$. 

Since $v$ is smooth in $ {\overline {A_n}}$, a standard density argument provides a function $\varphi \in \mathcal{C}^1_c (\R^N)$, with $supp(\varphi) \subset A_n$, such that 

\beq
\int_{\R^N} \vert \nabla \varphi \vert^2 - (1-3u^2){\varphi}^2  = \int_{A_n} \vert \nabla \varphi \vert^2 - (1-3u^2){\varphi}^2   < 0. 
\eeq
The latter being valid for any $n\ge1$, the claim is proven.  

\end{proof}

\medskip

\section{Saddle-type solutions : the two-dimensional case}\label{sec : appendix C}

In this section we prove that the saddle-type solutions in $\R^2$, constructed in \cite{DFP, ACM, FLP} have finite energy and we also calculate the exact value of their energy. Theorem \ref{Teo-finite-energy}, about the general structure of finite-energy solutions, will play a crucial role in the analysis. Also, the final result will be achieved through some intermediate results which we believe are of independent interest. 

\bigskip

Let us recall that, for every integer $ k \ge 2$, the saddle-type solutions $u_k$ built in \cite{DFP, ACM, FLP} are entire solutions of \eqref{eqAC} with dihedral symmetry and such that their nodal set is exactly 
equal to the union of $k$ straight lines through the origin forming angles equal to $\frac{\pi}{k}$.  More precisely, the nodal set is exactly the union of the $k$ straight lines $\{ (x_1,x_2) \in \R^2 : \cos \Big (\frac{j\pi}{k}\Big) x_2 = \sin \Big (\frac{j\pi}{k} \Big) x_1 \}$, $ j=0,...,k-1$.  
All these solutions are furthermore characterized by the fact that along directions not 
parallel to that lines they tend to $\pm 1$ in an alternate way. 

\medskip

The main result of this section is

\medskip

\begin{theorem} \label{Teo-Energia-saddle-2D} 
Assume $ k \ge 2$ and let $u_k$ be the corresponding saddle-type solution. Then $u_k$ has finite energy and

\beq\label{valore-energia-saddle}
L_{u_k} =  2 k {c_{\star}} = k  {\ell}_{min}.
\eeq

\end{theorem}{

\medskip

To prove the above Theorem \ref{Teo-Energia-saddle-2D} we shall need the following refinement of Lemma 9.1. of \cite{FV} (which in turn, was inspired by the proof of Theorem 1.3 of \cite{AC}). 

\medskip

Given a unit vector $\eta \in \R^N$, we say that a domain $\Omega \subset \R^N$ is {\it $\eta-$invariant } if 
$\Omega + t \eta \subseteq \Omega$ for every $ t>0$. 

Typical examples of invariant domains are provided by \emph{epigraphs}. Let us recall that a domain $\Omega \subset \R^N$, $N \ge 2$, is a \emph{locally lipschitz epigraph} if, up to a rotation of the space, there exists $g \in \mathcal{C}^{0,1}_{loc}(\R^{N-1},\R)$ such that $ \Omega = \{(x',x_N) \in \R^{N-1} \times \R \, : \, x_N > g(x') \}$.

\medskip

\begin{lemma} \label{Refin-Lemma} 
Let $ N \ge 2$ and let $\eta$ be a unit vector of $\R^N$. Let $\Omega \subset \R^N$ be a $\eta-$invariant Lipschitz domain. Let $v \in \mathcal{C}^2(\overline{\Omega})$ 
be a solution of 

\beq\label{pbl-in-Omega}
\begin{cases}
-\Delta v = v-v^3& \text{in \quad $\Omega$},\\
\frac{\partial v}{\partial \eta}  > 0& \text{in \quad $\Omega$},\\
\sup_{x \in \Omega} \vert \nabla v(x) \vert < + \infty.  &
\end{cases}
\eeq

\noindent For $ t \ge 0$ and $ x \in \Omega$ set $ v^{t} (x) := v(x + t \eta)$. For $ x_0 \in \partial \Omega$ and $ R>0$ set $ \Omega_R = \Omega \cap B(x_0,R)$ and 

\beq \label{funz-ener-Omega}
{\mathcal{E}}(v^{t}, x_0, R):=\int_{\Omega_R} \frac{|\nabla v^{t}|^2}{2}
+  \frac{(1-(v^{t})^2)^2}{4}.
\eeq
Then
\beq \label{stima-ener-Omega}
{\mathcal{E}}(v, x_0, R) \le {\mathcal{E}}(v^{t}, x_0, R)  + \Big (\sup_{x \in \Omega} \vert \nabla v(x) \vert \Big)  \int_{\partial (\Omega_R)} (v^{t}-v)  d\sigma.
\eeq

\end{lemma}

\begin{proof} Since $\Omega $ is $\eta-$invariant, the functions $v^t$ are well-defined on $ \Omega$ for every $t \ge 0$. Furthermore, they solve the Allen-Cahn equation on $\Omega$.  

Arguing as in the proof of Lemma 9.1. of \cite{FV} we obtain 

\beq\begin{split}
& \frac{d}{dt} {\mathcal{E}}(v^{t}, x_0, R)  = \int_{\Omega_R}  \nabla v^{t} \cdot \nabla (\partial_t v^{t}) 
+  [(v^{t})^3 - v^t]\partial_t v^{t} \\
& = \int_{\partial (\Omega_R)} \partial_t v^{t} (\nabla v^{t} \cdot \nu)  d\sigma \ge - \Big (\sup_{x \in \Omega} \vert \nabla v(x) \vert \Big) \int_{\partial (\Omega_R)} \partial_t v^{t}  d\sigma,
\end{split}
\eeq
where in the latter we have used the monotonicity assumption $ \frac{\partial v}{\partial \eta}  > 0$ on $ \Omega$. Therefore, integrating the above inequality over $ [0,T]$, $T >0$, and using Tonelli's Theorem we get
\beq\begin{split}
& {\mathcal{E}}(v^{T} , x_0, R) -  {\mathcal{E}}(v, x_0, R) = \int_0^T \frac{d}{dt} {\mathcal{E}}(v^{t}, x_0, R)  dt \\
& \ge -  \Big(\sup_{x \in \Omega} \vert \nabla v(x) \vert \Big)  \int_0^T \int_{\partial (\Omega_R)} \partial_t v^{t}  d\sigma dt = - \Big(\sup_{x \in \Omega} \vert \nabla v(x) \vert \Big) \int_{\partial (\Omega_R)}  \int_0^T  \partial_t v^{t}  dt d\sigma \\
& = - \Big(\sup_{x \in \Omega} \vert \nabla v(x) \vert \Big) \int_{\partial (\Omega_R)} (v^{T}-v)  d\sigma.
\end{split}
\eeq
Which concludes the proof.
\end{proof}

\medskip

\begin{theorem} \label{Teo-Stima-Energia-da-sopra} 
Assume $ k \ge 2$ and let $u_k$ be the corresponding saddle-type solution. Then $u_k$ has finite energy and

\beq\label{valore-energia-saddle-da-sopra}
L_{u_k} \le  {c_{\star}} \Big (\frac{3\pi}{2} + 3k \Big).
\eeq

\end{theorem}{

\medskip

\begin{proof}
By Theorem \ref{DefEnergia}, the energy of any solution is rotation-invariant and also it remains unchanged when we exchange $u$ with $-u$. Therefore, after a counterclockwise rotation of $ \frac{\pi}{2} - \frac{\pi}{2k}$ (and up to exchange $u_k$ with $-u_k$ if necessary) we can suppose that $u_k$ is positive in the domain given by the (open) infinite sector $\Omega $ defined, in polar coordinates $(\rho, \theta)$, by $ \rho> 0 $ and $ \frac{\pi}{2} - \frac{\pi}{2k} <  \theta < \frac{\pi}{2} + \frac{\pi}{2k}$ and that $ u_k =0 $ on $\partial \Omega$.

We observe that $ \Omega $ is the globally lipschitz epigraph defined by the function $ \varphi(x_1) = \tan (\frac{\pi}{2} - \frac{\pi}{2k}) \vert x_1\vert $ (actually a cone) and that $u_k \in \mathcal{C}^{2}({\overline{\Omega}})$ satisfies 

\beq\label{pbl-in-Omegak}
\begin{cases}
-\Delta u_k = u_k-{u_k}^3& \text{in \quad $\Omega$},\\
u_k >0 & \text{in \quad $\Omega$},\\
u_k =0 & \text{on \quad $\partial \Omega$},
\end{cases}
\eeq
therefore a standard moving plane procedure (cf. for instance \cite{EL,BCN3}) provides 

\beq\label{monot-u-k}
\frac{\partial u_k}{\partial x_2}  > 0\qquad  \text{in \quad $\Omega$}.
\eeq 

Problem \eqref{pbl-in-Omegak} has a unique solution (see Theorem 1.2. of \cite{BCN3} for instance) and so $u_k$ is even, since the function $ u_k(-x_1,x_2)$ is also a solution to  \eqref{pbl-in-Omegak}. This proves that $u_k$ has  dihedral symmetry and that, in each nodal domain $u_k$ tends to $\pm 1$ in an alternate way (actually it tends to $\pm 1$ along every direction not parallel to the boundary). 

The dihedral symmetry then implies 

\beq\label{2k-energia}
{\mathcal{E}}(u_k , B(0,R)) = 2k {\mathcal{E}}(u_k , B(0,R) \cap \Omega) = 2k {\mathcal{E}}(u_k , 0, R) 
\eeq 

We can thus focus on the estimate of the r-h-s of \eqref{2k-energia}. In view of the monotonicity property \eqref{monot-u-k}, we can apply Lemma \ref{Refin-Lemma} to get, from \eqref{stima-ener-Omega}, that $\forall R>0, t>0$, 

\beq
{\mathcal{E}}(u_k, 0, R) \le {\mathcal{E}}(u_k^{t}, 0, R)  + \Big (\sup_{x \in \Omega} \vert \nabla u_k(x) \vert \Big)  \int_{\partial (\Omega_R)} (u_k^{t}-u_k)  d\sigma
\eeq
and passing to the limit, $ t \to +\infty$, in the latter we have for every $R>0$ 

\beq \label{quasi-stima}
{\mathcal{E}}(u_k, 0, R) \le 0 + \Big (\sup_{x \in \Omega} \vert \nabla u_k(x) \vert \Big)  \int_{\partial (\Omega_R)} (1-u_k)  d\sigma,
\eeq
since by construction, $ u_k \to 1$ in $\Omega$ along the direction $x_2$. 

On the other hand $ u_k >0$ on $\Omega$ by \eqref{pbl-in-Omegak} and 

\beq \label{Modica-gradbound}
\forall x \in \R^2 \qquad \frac{|\nabla u_k(x)|^2}{2}  \le \frac{(1-(u_k)(x)^2)^2}{4} 
\eeq
Then $ 0 < 1 - u_k \le 1 $ in $ {\overline \Omega} $ and  $\sup_{x \in \Omega} \vert \nabla u_k(x) \vert \le \frac{\sqrt2}{2}$. The insertion of all these informations into \eqref{quasi-stima} then yields
\beq \label{quasi-stima2}
\forall R>0 \qquad {\mathcal{E}}(u_k, 0, R) \le \frac{\sqrt2}{2} {\cal H}^1(\partial (\Omega_R)) = 
\frac{\sqrt2}{2} \Big( \frac{\pi R}{k} + 2R \Big).
\eeq
From the latter and \eqref{2k-energia} we infer 
\beq \label{quasi-stima3}
\forall R>0 \qquad {\mathcal{E}}(u_k, B(0,R)) \le \sqrt2 \Big( \pi  + 2k \Big) R = c_{\star} \Big( \frac{3 \pi}{2} + 3k \Big) R
\eeq
which gives the desired conclusion.

\end{proof}

\smallskip

We can now turn to the proof of Theorem \ref{Teo-Energia-saddle-2D}

\medskip

\begin{proof}[Proof of Theorem \ref{Teo-Energia-saddle-2D}]




From Theorem \ref{Teo-Stima-Energia-da-sopra} we know that $u_k$ has finite energy, therefore we can apply Theorem \ref{Teo-finite-energy}. Let us now identify the set $E$ and the varifold $V$ given by Theorem \ref{Teo-finite-energy}. The asymptotics of $u_k$ imply that the set $E$ is (equivalent to) the cone $ \{ \, x \in \R^2 \, :  u_k(x) >0 \} $. Moreover, the support of the limit varifold $V$ is exactly the set of zeros of $u_k$, i.e., $\partial E$. To see this we observe that, for every $\eps>0$, the set of zeros of the rescaled function $(u_k)_{\eps}$ coincides with the set of zeros of $u_k$ and thus the claim follows by the fact that the nodal set of $(u_k)_{\eps}$ converges to $ supp \Vert V \Vert $ locally in the Hausdorff distance sense (by item (3) of Theorem 1 of \cite{HT}). We also note that the convergences in item i) of Theorem \ref{Teo-finite-energy} take place for all the family $(u_k)_{\eps}$ and not only for a subsequence. 

Having identified $E$ and $V$ we now prove that {all} the $2k$ open half-lines forming the $ supp \Vert V \Vert \setminus \{0\}$ have the same multiplicity. This result follows immediately from Proposition \ref{Sim1} and the fact that $u_k$ has  dihedral symmetry.

By Theorem \ref{Teo-finite-energy-2D}  we also know that $V$ is integral and conical about the origin, thus 

\beq\label{dens-semiretta}
2 \Theta (V,0) = 2k m 
\eeq 
where $m$ is the multiplicity of the open half-line $ \{ (x_1,x_2) \in \R^2 \, : \, x_1 >0 \}$ (and actually of any other open half-line forming the $supp \Vert V \Vert$, as proved before).  This and  \eqref{energiaV} yield 

\beq\label{quasi-energia-saddle}
L_{u_k}= 2k m c_{\star}
\eeq 
where $m$ is an {\it odd} integer since the open half-line $ \{ (x_1,x_2) \in \R^2 \, : \, x_1 >0 \}$ is a {\it transition line}. 

Now we prove that $ m=1$ in \eqref{quasi-energia-saddle}, which will give the desired conclusion \eqref{valore-energia-saddle}.  By comparing \eqref{quasi-energia-saddle} and \eqref{valore-energia-saddle-da-sopra} we find
$2k m \le \frac{3\pi}{2} + 3k$ and so $m \le \frac{3\pi}{8} + \frac{3}{2}$, since $ k \ge 2$. Then, $ m \le \frac{3\pi}{8} + \frac{3}{2}< 3 $, which immediately gives $m=1$, since $m$ is an odd integer. 

\end{proof}

\section{Saddle-type solutions : the higher-dimensional case}\label{sec : appendix D}

This section is devoted to the study of saddle-type solutions in dimension $N\ge3$. 

The first result concerns the existence of finite energy solution with arbitrary large energy in any dimension. In particular we shall prove that {\it any integer multiple} of the {\it minimum positive energy level} $ \ell_{min}$ is a positive energy level and that it is occupied by an explicit saddle-type solution. 

\medskip

\begin{theorem} \label{energia-arbitr-grande} 
Assume $ N \ge 3$. For every integer $k \ge 2$ there exists a solution $u$ of the Allen-Cahn \eqref{eqAC} on $\R^N$ 
such that 

\beq\label{energia-arbitraria}
L_u =  k \omega_{N-1} \frac{2\sqrt{2}}{3}  = k \ell_{min}.
\eeq
In particular, there are finite energy solutions with arbitrary large energy.
\end{theorem}

\medskip

\begin{proof}
For $ N\ge 3$ and $ x \in \R^N$ we set $ u(x) : = u_k(x_1,x_2)$, where $u_k$ is the two-dimensional saddle-type solution considered in Appendix \ref{sec : appendix C}. Since $u_k$ is a finite-energy  entire solution of the Allen-Cahn equation in $\R^2$, we can apply Proposition \ref{Sim2} about cylindrical solutions to get that $ supp (V_u) = supp (V_{u_k}) \times \R^{N-2}$ and $ \Theta (V_u,0) = \Theta (V_{u_k}, 0)$. On the other hand (see \eqref{dens-semiretta}) in the course of the proof of Theorem \ref{Teo-Energia-saddle-2D} we proved that $ \Theta (V_{u_k}, 0) = k$. The conclusion then follows by applying \eqref{energiaV}. 
\end{proof}

\medskip

In view of the results established in Section \ref{sec : finite-energy-2D} for the two-dimensional case and the above Theorem \ref{energia-arbitr-grande}, one might ask whether or not a quantization effect holds true also in the higher-dimensional case. To investigate this question we analyse the saddle-type solutions built in every even dimension $N=2m\ge 4$ in \cite{CT1}. 

For $ m \ge 2$ and $(x,y) \in \R^m \times \R^m = \R^{2m}$, the saddle-type solution $u$ built in \cite{CT1} is an entire solution of \eqref{eqAC} which depends only on the radial variables $\vert x \vert $  and $\vert y \vert $, is positive in the open (solid) cone $ \{  \vert x \vert  > \vert y \vert \} $ and is odd with respect to $ \{  \vert x \vert  = \vert y \vert \}$, i.e. $ u(x,y) = u(\vert x \vert, \vert y \vert \} = - u( \vert y \vert,   \vert x \vert)$.  In particular, $u$ vanishes exactly on the Simons cone
\beq\label{Cono-Cabré}
{\cal C} : = \{ (x,y) \in \R^{2m} \, : \,  \vert x \vert^2  = \vert y \vert^2 \} = \partial  \{  \vert x \vert  > \vert y \vert \}.
 \eeq

\noindent Furthermore, $u$ has finite energy. 


\bigskip

In the second result of this section we compute the energy of this solution. More precisely we have 

\medskip

\begin{theorem} \label{energia-saddle-Cabré} 
Assume $ N = 2m \ge 4$ and let $u$ be the corresponding saddle-type solution.   Then 
\beq\label{valore-energia-saddle-Cabré}
L_{u} =   {c_{\star}} P(E,B(0,1)) = \frac{P(E,B(0,1))}{\omega_{N-1}} \, {\ell}_{min} 
\eeq
\noindent where $P(E,B(0,1))$ denotes the perimeter of the solid cone $ E = \{  \vert x \vert  > \vert y \vert \} $ in the open unit ball $B(0,1).$
\end{theorem}



\begin{proof}
Since $u$ has finite energy, we can apply Theorem \ref{Teo-finite-energy} to evaluate the energy $L_u$. To this end we observe that the support of the limit varifold $V$ is exactly the set of zeros of $u$, i.e., the Simons cone ${\cal C}$. To see this we note that, for every $\eps>0$, the set of zeros of the rescaled function $u_{\eps}$ coincides with the set of zeros of $u$ (since this set is a cone) and thus the claim follows by the fact that the nodal set of $u_{\eps}$ converges to $ supp \Vert V \Vert $ locally in the Hausdorff distance sense (by item (3) of Theorem 1 of \cite{HT}). From the asymptotics and the monotonicity properties of $u$ proved in \cite{CT1, CT2, CT3}, we infer that the set of locally finite perimeter $E$ coincides with the set $ \{  \vert x \vert  > \vert y \vert \}.$ We also note that the convergences in item i) of Theorem \ref{Teo-finite-energy} take place for all the family $u_{\eps}$ and not only for a subsequence. 

Since the Simons cone ${\cal C}$ is smooth outside the origin and that the solution $u$ is $O(m)\times O(m)-$invariant and also satisfies $ u(\vert x \vert, \vert y \vert ) = - u( \vert y \vert,   \vert x \vert)$, an application of Proposition \ref{Sim1} gives that the density function is constant on ${\cal C} \setminus \{0\}$.  Let us denote by $M$ this integer.  Thus 

\beq
{\omega_{N-1}} \Theta(V,0) =  \Vert V \Vert ( B(0,1)) = M P(E,B(0,1)).
\eeq

By combining  the latter and \eqref{energiaV} we get 
\beq\label{energiaM}
L_{u} =  {c_{\star}}{\omega_{N-1}} \Theta(V,0) = {c_{\star}} M P(E,B(0,1)) =  {\ell}_{min} M \frac{P(E,B(0,1))}{\omega_{N-1}}. 
\eeq

Since $ supp \Vert V \Vert = supp \Vert \partial E \Vert =\partial E = {\cal C}$ and it is smooth outside the origin, 
$M$ must be an odd integer by item (4) of  Theorem 1 of \cite{HT}. 

\medskip

Next we plan to prove that $M$ is equal to 1. This will provide the desired conclusion.

\medskip

To this end, we exploit once again the symmetries enjoyed by $u$, as well as the asymptotics and the monotonicity properties of $u$ proved in \cite{CT1, CT2, CT3}. Since $ u(\vert x \vert, \vert y \vert ) = - u( \vert y \vert,   \vert x \vert)$ and $u$ is even with respect to the variable $x_1$, arguing as in the proof of Theorem \ref{Teo-Energia-saddle-2D} 
we get  
\beq\label{energia-risc-2-palla}
{\mathcal{E}}(u, B(0,R) ) =  2 {\mathcal{E}}(u, B(0,R) \cap E),
\eeq 
\beq\label{energia-risc-1-semipalla}
{\mathcal{E}}(u, B(0,R) \cap E )= {\mathcal{E}}(u, E_R^{+}) + {\mathcal{E}}(u, E_R^{-}), 
\eeq 

\noindent where $ E_R^{+} : = B(0,R) \cap E \cap \{x_1 >0\}$ and $ E_R^{-} : = B(0,R) \cap E \cap \{x_1 <0\}$.  

To proceed further we observe that $E^{+} : = E \cap \{x_1 >0\}$ is an $e_1-$invariant domain, while $ E^{-} := E \cap \{x_1 <0\}$ is a $(-e_1)-$invariant domain, where as usual, we have denoted by $e_1$ the first vector of the canonical basis of $\R^N$. We also note that, thanks to the monotonicity properties of $u$ established in \cite{CT1, CT2, CT3}, we have that $\frac{\partial u}{\partial x_1}  > 0$ on $E^{+} $ and that $\frac{\partial u}{\partial x_1}  < 0$ on $E^{-} $.
We can apply the integration-by-parts argument from Lemma \ref{Refin-Lemma} to the locally Lipschitz domains  $E^{\pm} \setminus \overline{B(0,\eta)}$, $\eta>0$. The translation in the direction $\pm e_1$ preserves those domains and the contribution from the inner boundary is $O(\eta^{N-1})$, as $\eta\longrightarrow 0$.  
Therefore, for every $t \ge 0$ and every $ R>0$,  
\beq\label{energia-risc-3-semipalla} 
{\mathcal{E}}(u, E_R^{+}) = {\mathcal{E}}(u, 0, R) \le {\mathcal{E}}(u^{t}, 0, R)  + \Big (\sup_{x \in E} \vert \nabla u(x) \vert \Big)  \int_{\partial (E_R^{+})} (u^{t}-u)  d\sigma. 
\eeq
Passing to the limit, as $ t \to +\infty$, in the latter we have for every $R>0$ 

\beq \label{quasi-stima+}
{\mathcal{E}}(u, E_R^{+}) \le 0 + \Big (\sup_{x \in E} \vert \nabla u(x) \vert \Big)  \int_{\partial (E_R^{+})}  (1-u)  d\sigma,
\eeq
since by construction, $ u \to 1$ in $E$ along the direction $e_1$.  

On the other hand $ u >0 $ on $E$ and 

\beq \label{Modica-gradbound2}
\frac{|\nabla u|^2}{2}  \le \frac{(1-u^2)^2}{4} \qquad {\rm on} \qquad \R^N 
\eeq
Then $ 0 < 1 - u \le 1 $ in $ {\overline E}$, $\sup_{x \in E} \vert \nabla u(x) \vert \le \frac{\sqrt2}{2}$ and so

\beq\begin{split} \label{quasi-stima+2}
& {\mathcal{E}}(u, E_R^{+}) \le  \frac{\sqrt2}{2} \int_{\partial (E_R^{+})}  (1-u)  d\sigma \\
& \le \frac{\sqrt2}{2} \Big [ \int_{E^{+} \cap \partial B(0,R)}  (1-u)  d\sigma  + \int_{\partial E \cap B^{+}(0,R)}  (1-u) d\sigma  + \int_{X_R}  (1-u)  d\sigma  \Big ] \\
&\le \frac{\sqrt2}{2} \Big [ \int_{E^{+} \cap \partial B(0,R)}  (1-u)  d\sigma + P(E,B^{+}(0,R)) + {\cal H}^{N-1} (X_R) \Big ]. 
\end{split}
\eeq
where $  B^{+}(0,R) : = B(0,R) \cap \{x_1 >0\}$ and $ X_R : = \{x_1 = 0\} \cap E \cap B(0,R)$.   

\medskip

In the same way one obtains, for every $R>0$,  

\beq\begin{split} \label{quasi-stima-2}
& {\mathcal{E}}(u, E_R^{-}) \le  \frac{\sqrt2}{2} \int_{\partial (E_R^{-})}  (1-u)  d\sigma \\
&\le \frac{\sqrt2}{2} \Big [ \int_{E^{-} \cap \partial B(0,R)}  (1-u)  d\sigma + P(E,B^{-}(0,R)) + {\cal H}^{N-1} (X_R) \Big ]. 
\end{split}
\eeq
where $  B^{-}(0,R) : = B(0,R) \cap \{x_1 <0\}$.

Putting together \eqref{energia-risc-2-palla}, \eqref{energia-risc-1-semipalla}, \eqref{quasi-stima+2} and \eqref{quasi-stima-2} we get, for every $R>0$,  

\beq \begin{split}\label{quasi-stima-finale}
{\mathcal{E}}(u, B(0,R)) \le \sqrt2 \Big [ P(E,B(0,R)) + \int_{E \cap \partial B(0,R)}  (1-u)  d\sigma + 2 {\cal H}^{N-1} (X_R) \Big ]
\end{split}
\eeq

\medskip

Next we estimate $\int_{E \cap \partial B(0,R)}  (1-u)  d\sigma $ for large values of $R$. 

To this end, for every $ \delta >0$, we consider a compact set $ \omega_{\delta} \subset E \cap \partial B(0,1) $ such that ${\cal H}^{N-1} ((E \cap \partial B(0,1)) \setminus \omega_{\delta}) \le \delta$. Since  $u_{\eps} \to 1$ uniformly on compact sets contained in the open cone $E$ (by Theorem 1 of \cite{HT}, see also \cite{CT1, CT2, CT3}) there is $ R_{\delta} > 0 $ such that 

\beq\label{conv1Simons}
\forall \, R > R_{\delta} \qquad 0 \le 1 - u \le \delta \qquad {\rm on} \qquad R\omega_{\delta} \subset E \cap \partial B(0,R). 
\eeq
Therefore, for every $\delta >0$ and every $R > R_{\delta}$, 

\beq\begin{split}\label{stima-sfera-Simons}
0 \le &\int_{E \cap \partial B(0,R)}  (1-u)  d\sigma = \int_{(E \cap \partial B(0,R)) \setminus R\omega_{\delta}}  (1-u)  d\sigma + \int_{R\omega_{\delta}}  (1-u)  d\sigma\\
& \le \int_{(E \cap \partial B(0,R)) \setminus R\omega_{\delta}}  1 d\sigma + \int_{R\omega_{\delta}}  \delta d\sigma\\
& \le {\cal H}^{N-1} ((E \cap \partial B(0,R)) \setminus R\omega_{\delta}) + \delta {\cal H}^{N-1} (R\omega_{\delta})\\
& \le R^{N-1} {\cal H}^{N-1} ((E \cap \partial B(0,1)) \setminus \omega_{\delta}) + \delta R^{N-1}  {\cal H}^{N-1} (\omega_{\delta})\\
& \le R^{N-1} \delta + R^{N-1} \delta  {\cal H}^{N-1} (\omega_{\delta}),
\end{split}
\eeq
where in the penultimate line we have used the fact that $E$ is a cone. 

Since $ \omega_{\delta} \subset E \cap \partial B(0,1) $ for every $ \delta >0$, we see that ${\cal H}^{N-1} (\omega_{\delta}) \le C(N)$, for some positive constant depending only on $N$, thus combining \eqref{quasi-stima-finale} and \eqref{stima-sfera-Simons} we are led to : $ \forall \, \delta>0, \, \forall \, R > R_{\delta} \quad $ 

\beq \label{quasi-stima-finale2}
{\mathcal{E}}(u, B(0,R)) \le \sqrt2 \Big [ P(E,B(0,1))+ 
\delta (1 + C(N)) + 2 {\cal H}^{N-1} (X_1)\Big ] R^{N-1}.
\eeq
Letting $R\to +\infty$ in the latter and recalling the definition of energy $L_u$, we immediately get
\beq \label{quasi-stima-finale3}
\forall \, \delta>0, \qquad L_u \le \sqrt2\Big [ P(E,B(0,1))+ 
\delta (1 + C(N)) + 2 {\cal H}^{N-1} (X_1)  \Big ],
\eeq
letting $\delta \to 0$ in the latter we finally obtain
\beq \begin{split}\label{quasi-stima-finale4}
&L_u \le {\sqrt2}\Big [ P(E,B(0,1))+ 2 {\cal H}^{N-1} (X_1)  \Big ] 
\end{split}
\eeq
Set $x = (x_1, x')$ and $y=(y_1,y')$, with $x',y' \in \R^{m-1}$, then
\beq \begin{split}\label{quasi-stima-finale4-bis}
 X_1 & = \left\lbrace (x',y) \in  \R^{m-1} \times \R^m \, : \, \vert x'\vert^2 > y_1^2 + \vert y'\vert^2 \right\rbrace \cap B^{N-1}(0,1) \\
& \subset \left\lbrace (x',y) \in  \R^{m-1} \times \R^m \, : \, \vert x'\vert^2 > \vert y'\vert^2 \right\rbrace \cap B^{N-1}(0,1) =: \mathcal{A} \subset \R^{N-1}
\end{split}
\eeq
and 
\beq \label{quasi-stima-finale4-tris}
{\cal H}^{N-1} (\mathcal{A}) \leq \frac{\omega_{N-1}}{2}.
\eeq
Indeed, the isometry $T(x',y_1,y') = (y',y_1,x')$ is such that $\mathcal{A} \cap T(\mathcal{A}) = \emptyset $ and $ \mathcal{A} \cup T(\mathcal{A}) \subset B^{N-1}(0,1)$, therefore $ 2 {\cal H}^{N-1} (\mathcal{A}) = {\cal H}^{N-1} (\mathcal{A}) + {\cal H}^{N-1} T(\mathcal{A}) \leq {\cal H}^{N-1} (B^{N-1}(0,1)) = \omega_{N-1}$. 

\medskip

From \eqref{quasi-stima-finale4-bis} and  \eqref{quasi-stima-finale4-tris} we deduce that $2 {\cal H}^{N-1} (X_1) \leq \omega_{N-1}$, whose combination with \eqref{quasi-stima-finale4} leads to 
\beq \begin{split}\label{quasi-stima-finale5}
L_u &\le {\sqrt2}\Big [ P(E,B(0,1))+ \omega_{N-1}  \Big ] \\ 
& = {\sqrt2}\Big [ P(E,B(0,1)) +  \frac{\omega_{N-1}}{P(E,B(0,1))} P(E,B(0,1)) \Big ]\\
& = {\sqrt2}  P(E,B(0,1)) \Big [1 + \frac{\omega_{N-1}}{P(E,B(0,1))}\Big ] \\
& = {\sqrt2}  P(E,B(0,1)) \Big [1 + \frac{1}{\Theta({\cal C},0)}\Big ] 
\end{split}
\eeq
where $\Theta({\cal C},0)$ is the density of the Simons cone ${\cal C}$ at the origin. \\
Comparing the latter with \eqref{energiaM} we obtain
\beq
M \leq \frac{3}{2} \Big [1 + \frac{1}{\Theta({\cal C},0)}\Big ] < 3
\eeq
since $\Theta({\cal C},0) >1$.  Hence $M=1$, since $M$ is an odd integer. \end{proof}

\medskip

As a consequence of Theorem \ref{energia-saddle-Cabré} and of Proposition \ref{Teo-dens-coni-Simons} we see that, in every dimension $ N \ge 4$, there are energy levels which are not an integer multiple of the minimum positive energy level $\ell_{min}$. Therefore a quantization effect like the one demonstrated in the two-dimensional case, does not hold in higher dimensions. 

This is the content of the next

\medskip

\begin{corollary}\label{CorCabré}
For every $N \ge 4$ there exists a solution $u$ of the Allen-Cahn \eqref{eqAC} on $\R^N$ 
such that 

\beq\label{energia-non-intera}
L_u = \frac{\pi}{2} \ell_{min}.
\eeq
In particular, there are energy levels which are not an integer multiple of the minimum positive energy level $\ell_{min}$.

\end{corollary}

\begin{proof}

In view of formula \eqref{valore-energia-saddle-Cabré} and Proposition \ref{Teo-dens-coni-Simons}, the saddle-type solution $u=u_2$ built in \cite{CT1} for $ m = 2$ and $(x,y) \in \R^2 \times \R^2 = \R^{4}$, is an entire solution of \eqref{eqAC} satisfying $L_u = \frac{\pi}{2} \ell_{min}$. This proves the claim for $N=4$. For $ N\ge 5$ and $ x \in \R^N$ we set $ u(x) : = u_2(x_1,x_2,x_3,x_4)$, where $u_2$ is the 4-dimensional saddle-type solution just considered. The desired conclusion then follows by applying \eqref{energiaV} and Proposition \ref{Sim2} about cylindrical solutions. 
\end{proof}

\medskip

Now, for every saddle-type solution considered in the previous theorem, we compute the exact value of its energy, as well as the asymptotics behaviour of the energy, as the Euclidean dimension goes to infinity.

\medskip

\begin{proposition} \label{Energia-saddle+asint} 
Assume $ N = 2n \ge 4$ and let $u = u_n$ be the corresponding saddle-type solution.  Then 

\beq\label{val-energ-sad}
L_{u} =  \Theta({\cal C}_{n,n},0) \, {\ell}_{min} 
\eeq
where

\beq\label{valore-dens-Simons}
\Theta({\cal C}_{n,n},0) = \begin{cases}
\frac{n^2 (2n-3) !! }{(n !!)^2 \, 2^{n-2}}& \text{if $n$ is odd},\\
\quad  & \quad \\
\pi \, \frac{n^2 (2n-3) !! }{({\frac{n}{2}} ! )^2 \, 2^{2n-1}}& \text{if $n$ is even}.
\end{cases}
\eeq
Moreover, we have that 

\beq \label{asint-dens-Simons}
\lim_{n\to + \infty} \Theta({\cal C}_{n,n},0) = \sqrt{2}. 
\eeq

\beq \label{asint-energia-saddles}
\lim_{n\to + \infty} \frac{L_{u_n}(2n)}{\ell_{min}(2n)} = \sqrt{2}. 
\eeq

\end{proposition}

\medskip

\begin{proof}
By Proposition \ref{Teo-dens-coni-Simons} we have that 

\beq\label{Theta-Simons}
\Theta({\cal C}_{n,n},0) = \frac{(\omega_n)^2}{\omega_{2n-1}} \frac{n^2}{(2n-1) 2^{n-1}}
\eeq
and recalling that 

\beq
\omega_n = \begin{cases}
\frac{2^{\frac{n+1}{2}} \pi^{\frac{n-1}{2}}}{n !!}& \text{if $n$ is odd},\\
\quad  & \quad \\
\frac{\pi^{\frac{n}{2}}}{({\frac{n}{2}} ) ! }& \text{if $n$ is even}.
\end{cases}
\eeq
we obtain

\beq
\Theta({\cal C}_{n,n},0) = \begin{cases}
\frac{2(2n-1)!!}{(n!!)^2} \frac{n^2}{(2n-1) 2^{n-1}}& \text{if $n$ is odd},\\
\quad  & \quad \\
\pi \frac{(2n-1)!!}{({\frac{n}{2}} ! )^2 \, 2^n} \frac{n^2}{(2n-1) 2^{n-1}}&  \text{if $n$ is even}.
\end{cases}
\eeq
which yields \eqref{valore-dens-Simons}. 

To prove \eqref{asint-dens-Simons} (and thus also \eqref{asint-energia-saddles}) we recall that $\omega_n = \frac {\pi^{\frac{n}{2}}}{\Gamma(\frac{n}{2} + 1)} $, where $ \Gamma$ denotes the {\it gamma function}, therefore, from \eqref{Theta-Simons} we get that

\beq\label{asint2}
\Theta({\cal C}_{n,n},0) = \sqrt{\pi} \, \frac{\Gamma(n + \frac{1}{2})}{\Gamma^2(\frac{n}{2} + 1)} \frac{n^2}{(2n-1) 2^{n-1}}
\eeq

On the other hand : $ \Gamma (z) \sim \sqrt{\frac{2 \pi }{z}} \big (\frac{z}{e} \big)^z$ as $z \to \infty$ and so, 

\beq
\frac{\Gamma(n + \frac{1}{2})}{\Gamma^2(\frac{n}{2} + 1)}  \sim \frac{1}{\sqrt{2 \pi}} \frac{\frac{n}{2} + 1}{\sqrt{n + \frac{1}{2}}}  \frac{\sqrt{n}}{n^2 \, 2^{-(n+2)}}
\eeq
which, combined with \eqref{asint2}, immediately gives the desired conclusion. 

\end{proof}

\begin{remark}\label{Rk10}
AGGIUNGERE ANALISI stabilit\`a,stabilit\`a  fuori dal compatto/minimizzione per la saddle.
\end{remark}

\medskip

\section{A genuine three-dimensional saddle-type solution}\label{sec : appendix D3}

In the previous sections we considered finite-energy solutions in even dimensions $N =2n$. Clearly, those solutions provide also higher dimensional finite-energy solutions of cylidrical type, by adding redundant variables, as done in Theorem \ref{energia-arbitr-grande}. It is then natural to ask whether or not a genuine three-dimensional finite-energy solution exists. Here by {\it genuine}, we mean a solution which truly depends on 3 variables, i.e. it does not come from a 2D solution. In the present section we shall prove that the three-dimensional saddle-type solution built in \cite{AM} fulfills all the desired requirements. To this end we recall that the solution built in \cite{AM}, and denoted by $\theta_3$, satisfies the following properties (see Theorem 1.1. in \cite{AM}) : 
its nodal set coincides with the union of the three cartesian planes, $\theta _3 >0$ on the octant $ \{ \, (x,y,z) \in \R^3\, \ : \, x,y,z>0 \, \}$, it is increasing in each variable on $ \{ \, (x,y,z) \in \R^3\, \ : \, x,y,z>0 \, \}$, it is odd with respect to each variable and symmetric with respect to the diagonal planes $x = y, x = z$ and $y = z$ Also, $\theta_3$ is asymptotic to $\pm 1$ along directions parallel to the space diagonals. 

Also, as pointed out in \cite{AM}, $\theta_3$ has infinite Morse index (actually it is unstable outside every compact set of $\R^3$), since it converges to a 2D saddle-type solution when $z \to +\infty$ (which is unstable).

With all these informations we can prove 

\bigskip

\begin{theorem}\label{genuina} The solution  $\theta_3$ has finite energy and

\beq\label{valore-energia-saddle-3D}
L_{\theta_3} =  3 {\ell}_{min}.
\eeq
\end{theorem}

\medskip

\begin{proof}
We follow the proof of Theorem \ref{Teo-Stima-Energia-da-sopra} to demonstrate  that $\theta_3$ has finite energy and then we conclude as in the proof of  Thereom \ref{Teo-Energia-saddle-2D}.  Since the octant $ \{ \, (x,y,z) \in \R^3\, \ : \, x,y,z>0 \, \}$ is a lipschitz epigraph with respect to $ \eta = (\frac{1}{\sqrt{3}}, \frac{1}{\sqrt{3}},\frac{1}{\sqrt{3}})$ 
and $ \frac{\partial {\theta_3}}{\partial \eta}  > 0 $ on the octant, we can apply Lemma \ref{Refin-Lemma} to $\theta_3$ with $ \Omega = \{ \, (x,y,z) \in \R^3\, \ : \, x,y,z>0 \, \}$. Proceeding as in the proof of Theorem \ref{Teo-Stima-Energia-da-sopra} we get, $\forall \, R>0, $

\beq\label{8-energia}
{\mathcal{E}}(\theta_3 , B(0,R)) = 8 {\mathcal{E}}(\theta_3, 0, R) \le 8 \frac{\sqrt{2}}{2} \Big (\frac{3}{4} \pi R^2 + \frac{1}{8} s_2 R^2 \Big) = 8 \frac{\sqrt{2}}{2} \Big (\frac{3}{4} \pi R^2 + \frac{1}{8} 4 \pi R^2 \Big)
\eeq 
here we have used the oddness of $\theta_3$ with respect to each variable. The latter yields that $\theta_3$ has finite energy and also that

\beq\label{energ-oct-1}
L_{\theta_3} \le  8 \frac{\sqrt{2}}{2}  \Big(\frac{5}{4} \pi \Big) = \frac{15}{2}\ell_{min}.  
\eeq 

Since $\theta_3$ has finite energy, it is odd with respect to each variable and symmetric with respect to the diagonal planes $x = y, x = z$ and $y = z$, we can exploit Proposition \ref{Sim1} to obtain 

\beq\label{dens-oct}
\Theta (V,0) = 3 m 
\eeq 
where $m$ is the integer multiplicity of each connected component of the set $ (\{ \, x= 0\} \cup \{ \, y= 0\} \cup \{ \, z= 0\}) \setminus S $, and $S$ is the union of the three coordinate axes. Since the previous set is smooth, $m$ must be an odd integer by item (4) of  Theorem 1 of \cite{HT}. 

Now we apply the Structure Theorem \ref{Teo-finite-energy} to infer that

\beq\label{energ-oct-2}
L_{\theta_3} =  3 m \ell_{min}.  
\eeq 
The comparison of \eqref{energ-oct-1} and \eqref{energ-oct-2} gives that $ 1 \le m \le \frac{15}{6} < 3 $. Which implies $m =1$, since $m$ is odd, and concludes the proof. 

\end{proof}

\medskip

\section{The energy and the skeleton of the solution built by {\sl Del Pino, Kowalczyk and Wei}} \label{sec : appendix E}

\medskip

\bigskip

\begin{theorem}\label{energia-DKW}
Assume $ N=9$ and let $U$ be the solution of the Allen-Cahn equation \eqref{eqAC} built by M. del Pino, M. Kowalczyk and J. Wei in  \cite{DKW}.

\noindent Let $(\eps_n)_{n\ge1}$ be a sequence of positive real numbers such that $ \eps_n \to 0$. Then, the limit minimal cone $C$ provided by Theorem \ref{Teo-minloc} is the cylinder 

\beq\label{conf-cilindro-Simons}
C = E_{4,4} \times \R
\eeq
where $E_{4,4}  = \{ (x,y) \in \R^4\times \R^4  \,\, : \,\, \vert x \vert  > \vert y \vert \} $ is the solid minimal Simons cone.

\smallskip

\noindent Furthermore, the energy of $U$ is given by

\beq\label{valore-energia-DKW}
L_U = \frac{15}{32} \pi \, \ell_{min}  
\eeq
\smallskip

\noindent  Finally, the limit cone $C$ is unique, i.e. it does not dependent on the sequence $(\eps_n)_{n\ge1}$. 

\end{theorem}

\medskip

\begin{proof}  By item ii) of  Theorem \ref{monot-con-limiti-cost}, $U$ is a local minimizer and the associated limit minimal cone $C$ is a cylinder of the form $C' \times \R$ where $C'$ is a cone singular at the origin. Moreover, $ C' \subset \R^{N-1} \equiv \R^8$ and thus $C'$ has only the origin as singular point (i.e. it is smooth outside the origin) by the well-known regularity theory for minimal sets (see for  instance \cite{Giusti}). In particular, $C$ has exactly  the line $l \, := \{ \, (0',x_9) \, : \, x_9 \in \R \, \}$ as set of singular points. 
  
Now we use in parallel Theorem \ref{Teo-minloc} and Theorem \ref{Teo-finite-energy}.  Firts we note that $C$ and $E$ in Theorem \ref{Teo-minloc} and Theorem \ref{Teo-finite-energy} are equivalent and then from the minimality of $C$ we infer that $\Vert V \Vert = \Vert \partial C \Vert = Per (C, \cdot)$. Indeed, as pointed out in \cite{AAC} and already used in the proof of Theorem \ref{Teo-minloc}, we know that $ {\mathcal{E}_{\eps_n}} (u_{\eps_n}, \cdot) \longrightarrow c_{\star} Per(C, \cdot) $ locally weakly-$\star$ as measures. Therefore, $\Vert V \Vert = \Vert \partial C \Vert = Per (C, \cdot)$ by item i) and item ii) of Theorem \ref{Teo-finite-energy}. In particular we have that  $ supp \Vert V \Vert = supp \Vert \partial C \Vert = \partial C' \times \R$. 

Now we recall that, by construction, $U$ satisfies 

\beq\label{sim-DKW}
U(\vert x \vert, \vert y \vert, x_9) = - U(\vert y \vert, \vert x \vert, - x_9) \qquad \forall \, (x,y,x_9) \in \R^9
\eeq
and in particular $U$ vanishes on the set 
\beq
{\cal C}_{4,4} \times \{ 0\} \subset \R^8 \times \R = \R^9
\eeq
where ${\cal C}_{4,4}: = \{ (x,y) \in \R^4 \times \R^4 \,\, : \,\, \vert x \vert^2 = \vert y \vert^2  \} $ is the Simons' cones (see Appendix \ref {sec : appendix F}) and so, also $U_{\eps_n}$ vanishes on ${\cal C}_{4,4} \times \{ 0\} $. The latter and the fact that the sequence $(u_{{\epsilon}_n})_{n\ge 1}$ (given by Theorem \ref{Teo-finite-energy}) converges locally uniformly to $1$ (resp. $-1$) in each connected component of $(\R^9 \setminus supp \Vert V \Vert) $ (thanks to Theorem 1 of \cite{HT}) immediately yield 
\beq
{\cal C}_{4,4} \times \{ 0\} \subset supp \Vert V \Vert = \partial C' \times \R
\eeq
and so
\beq\label{inclusione}
{\cal C}_{4,4} \times \R \subset supp \Vert V \Vert = \partial C' \times \R
\eeq

To proceed further we invoke the unique continuation result for stationary codimension 1 integral varifolds established by N. Wickramasekera in \cite{W1} (see Theorem 1.2. therein).  We apply this theorem to the stationary codimension 1 integral (N-1)-varifolds (= 8-varifolds) ${\cal C}_{4,4} \times \R$ and $V$ (note that ${\cal C}_{4,4} \times \R $ is the boundary of $F \times \R$ which is a minimal cylinder). To this end we observe that all the assumptions for the application of the above mentioned theorem are satisfied, indeed $ supp \,  {\cal C}_{4,4} \times \R$ and $supp V$ are connected, ${\cal H}^{N-2}(Sing({\cal C}_{4,4} \times \R)) = {\cal H}^{N-2}(Sing(V)) = {\cal H}^{N-2} (l) = 0$ (since $l $ is a line and $ N=9$).  Finally, ${\cal H}^{N-1}(Supp \, ({\cal C}_{4,4} \times \R) \cap  supp \Vert V \Vert ) >0 $ since \eqref{inclusione} is in force. Therefore,  also the last assumption of Theorem 1.2 of \cite{W1} is satisfied by choosing $ \gamma =1$ therein.  The application of Theorem 1.2 of \cite{W1} then provides 

\beq
{\cal C}_{4,4} \times \R =  supp \Vert V \Vert = \partial C' \times \R,
\eeq
which establishes \eqref{conf-cilindro-Simons} and thus

\beq\begin{split}
& Per (C, B^9(0,1)) = \int_{-1}^1 {\cal H}^7 [(r + {\cal C}_{4,4} ) \cap B^9(0,1) ] dr = \\
& \int_{-1}^1 {\cal H}^7[ {\cal C}_{4,4}  \cap B^8(0,{\sqrt{1-r^2}})]  dr  =   \int_{-1}^1 \Theta( {\cal C}_{4,4} ,0') \omega_{7}  ({\sqrt{1-r^2}})^7dr = \\
&\Theta( {\cal C}_{4,4} ,0')  \int_{-1}^1 \omega_{7} ({\sqrt{1-r^2}})^7dr = \Theta( {\cal C}_{4,4} ,0') \omega_8. 
\end{split}
\eeq
where in the second line we have used that  ${\cal C}_{4,4} $ is a cone. The desired conclusion \eqref{valore-energia-DKW} then follows by applying Proposition \ref{Teo-dens-coni-Simons} and item ii) of Theorem \ref{Teo-minloc}. 
Since all the above arguments apply to any sequence $(\eps_n)_{n\ge1}$, the last part of the claim is also proven.

\end{proof}

\medskip

\section{Cones }\label{sec : appendix F}

\medskip

In this section we are concerned with the following family of cones 

\beq\label{coni-gen}
{\cal C}_{n,m}^{\gamma} : = \{ (x,y) \in \R^n \times \R^m \quad : \quad \vert x \vert^2 = \gamma^2 \vert y \vert^2  \}
\eeq
where $ n, m \ge1$ and $ \gamma >0$.   

The cone ${\cal C}_{n,m}^{\gamma}$ can be seen as the zero level set of the function $F(x,y) = \vert x \vert^2 - \gamma^2 \vert y \vert^2$ and so it is smooth outside the origin. Thus, its mean curvature at every point $(x,y) \neq (0,0)$, is given by  

\beq\label{mean-curv-coni}
 H({\cal C}_{n,m}^{\gamma}) : = div \Big ( \frac{ \nabla F(x,y)}{\vert \nabla F(x,y) \vert} \Big)
\eeq
and a straightforward calculation yields

\beq\label{mean-curv-coni-2}
 H({\cal C}_{n,m}^{\gamma}) : = div \Big ( \frac{ \nabla F(x,y)}{\vert \nabla F(x,y) \vert} \Big) = \frac{(n-\gamma^2 m)(\vert x \vert^2 + \gamma^4 \vert y \vert^2) - \vert x \vert^2 + \gamma^6 \vert y \vert^2}{(\vert x \vert^2 + \gamma^4 \vert y \vert^2)^{\frac{3}{2}}}.
\eeq
Using $\vert x \vert^2 = \gamma^2 \vert y \vert^2$ in the latter we get
\beq\begin{split}\label{mean-curv-coni-3}
 H({\cal C}_{n,m}^{\gamma}) =  \frac{ \gamma^2 \vert y \vert^2 [ (n-\gamma^2 m)(1 + \gamma^2) - 1 + \gamma^4] }{(\vert x \vert^2 + \gamma^4 \vert y \vert^2)^{\frac{3}{2}}} = & \\
 \frac{ (\gamma^2 \vert y \vert^2)(1 + \gamma^2)  [(n-\gamma^2 m) - (1 - \gamma^2)]}{(\vert x \vert^2 + \gamma^4 \vert y \vert^2)^{\frac{3}{2}}} 
\end{split}
\eeq

It follows that the cone ${\cal C}_{n,m}^{\gamma}$ has mean curvature $  H({\cal C}_{n,m}^{\gamma}) \equiv 0$, if and only if $ (n-\gamma^2 m) = (1 - \gamma^2)$ that is, if and only if $ n-1 = \gamma^2(m-1)$, i.e. if and only if

\beq
\begin{cases}
n=m=1 \qquad \text{\it and} \qquad \gamma >0,& \\
\text { }&\\
\text { \it or}&\\
\text { }&\\
n, m >1 \qquad \text{\it and} \qquad \gamma^2 = \frac{n-1}{m-1}. &
\end{cases}
\eeq

In the sequel, for $n,m>1$ and $ \gamma^2 = \frac{n-1}{m-1}$,  we shall denote the cone ${\cal C}_{n,m}^{\gamma} $  simply  by ${\cal C}_{n,m}. $ The cones

\beq\label{coni-Simons}
{\cal C}_{n,m}: = \{ (x,y) \in \R^n \times \R^m \quad : \quad (m-1)\vert x \vert^2 = (n-1)\vert y \vert^2  \}
\eeq
are known as Simons' cones and in the next result we calculate the exact value of their density at the origin. More precisely we have 

\medskip

\begin{proposition} \label{Teo-dens-coni-Simons} 
Assume $n,m>1$ (and thus $N=n+m \ge 4$). Then 
\beq\begin{split}\label{valore-dens-coni-Simons}
\omega_{n+m-1} \Theta({\cal C}_{n,m},0) = P(E_{n,m},B(0,1)) = & \\ s_{n-1}s_{m-1} \frac{(n-1)^{\frac{n-1}{2}} (m-1)^{\frac{m-1}{2} }}{(n+m-1) (n+m-2)^{\frac{n+m-2}{2}}} = & \\ \omega_{n}\omega_{m} \frac{nm (n-1)^{\frac{n-1}{2}} (m-1)^{\frac{m-1}{2} }}{(n+m-1) (n+m-2)^{\frac{n+m-2}{2}}} 
\end{split}
\eeq
where by $P( E_{n,m},B(0,1))$ we have denoted the perimeter of the open (solid) cone $E_{n,m} : = \{ (x,y) \in \R^n \times \R^m \, : \, (m-1)\vert x \vert^2 > (n-1)\vert y \vert^2  \} $ in the open unit ball $B(0,1)$, $ \Theta({\cal C}_{n,m},0) $ is the density of the cone ${\cal C}_{n,m} = \partial E_{n,m}$ at the origin, $ \omega_{k} $ is the volume of the unit ball in $\R^{k}$ and $s_{k-1}$ is the surface area of the unit sphere of $\R^k$.    

\end{proposition}

\begin{proof}
For any $r>0$ and any integer $k\ge2$ we denote by $ S^{k-1}_{r}$ the sphere of radius $r$ of the Euclidean space $\R^k$. Since the intersection of the Simons' cone ${\cal C}_{n,m}$ with $ S^{n+m-1}_{r}$ 
is $ S^{n-1}_{r \sqrt{\frac{n-1}{n+m-2}}} \times S^{m-1}_{r \sqrt{\frac{m-1}{n+m-2}}} $ we see that 

\beq\begin{split}\label{calc-dens-coni-Simons}
P(E_{n,m},B(0,1)) = \int_0^1 s_{n-1} \Big (r \sqrt{\frac{n-1}{n+m-2} }\Big )^{n-1}  s_{m-1} \Big (r \sqrt{ \frac{m-1}{n+m-2}} \Big )^{m-1} \, dr = & \\ 
s_{n-1}s_{m-1} \frac{(n-1)^{\frac{n-1}{2}} (m-1)^{\frac{m-1}{2} }}{(n+m-2)^{\frac{n+m-2}{2}}} \int_0^1 r^{n+m-2} dr = & \\
s_{n-1}s_{m-1} \frac{(n-1)^{\frac{n-1}{2}} (m-1)^{\frac{m-1}{2} }}{(n+m-1)(n+m-2)^{\frac{n+m-2}{2}}} = &\\
\omega_{n}\omega_{m} \frac{nm (n-1)^{\frac{n-1}{2}} (m-1)^{\frac{m-1}{2} }}{(n+m-1) (n+m-2)^{\frac{n+m-2}{2}}} 
\end{split}
\eeq
where in the latter we have used the relation $ k \omega_k = s_{k-1}$. 

\end{proof}

\medskip

\section{Open questions}\label{sec : appendix F1}

\medskip

\begin{Oquestion} Assume $ N\ge 9$. In view of the gap phenomena for local minimizers (cf. Theorem \ref{Gap-minloc}) and the local minimizer built by M. del Pino, M. Kowalczyk and J. Wei in \cite{DKW}, it is easily seen that the {\it second energy level for local minimizers} is well defined. As for $ \ell_2$, the second positive energy level for finite energy solutions, it is defined by  

\beq\label{defLiv2-min}
{\ell}_{2}( {\cal M}) : = \inf \{ \, L_u > \ell_{min} \,\, : \,\, u \in {\cal M}\, \}. 
\eeq

\bigskip

\noindent {\it Is ${\ell}_{2}( {\cal M})$ achieved by a member of ${\cal M}$ when $N \ge 9 $ ? In the affirmative, compute the exact value. }

\bigskip

\noindent A possible candidate for $N=9$, is the local minimizer built in \cite{DKW}. 

\end{Oquestion}

\medskip

\begin{Oquestion} {\it Does there exist ${\ell}_{2}( {\cal M})$ in $ \R^8$ ? }

\end{Oquestion}

\medskip

\begin{Oquestion} 

\noindent {\it Does there exist a local minimizer in $\R^8$ which is not 1D ?}

\medskip

\noindent A possible candidate is the saddle-type solution built by X. Cabr\'e and J. Terra in \cite{CT1} and studied in detail in Appendix \ref {sec : appendix D}. 

\end{Oquestion}

\medskip

\begin{Oquestion} 

\noindent {\it Prove (or disprove) that there are only a finite number of energy levels in the set ${\cal M}$ of the local minimizers. }

\medskip

\noindent For $ N\le 7$, there is only one energy level, the minimum one $\ell_{min}$, thanks to item iii) of Theorem \ref{Teo-minloc}. The upper bound proved in item iv) of Theorem \ref{Teo-minloc} suggests that this might be true also for $ N \ge 8$. 

\end{Oquestion}

\medskip

\begin{Oquestion} 

\noindent {\it Prove (or disprove) that there are only a finite number of energy levels in the set ${\cal M}_m$ of the monotone local minimizers. }

\medskip

\noindent For $ N\le 8$, thanks to item iii) of Theorem \ref{Teo-minloc} and to Theorem \ref{monot-1D-minloc}, there is only one energy level, the minimum one $\ell_{min}$.

\end{Oquestion}

\medskip

\begin{Oquestion} 

\noindent {\it Do there exist nonconstant, non-monotone local minimizers for $N \ge 9$ ?} 

\medskip

\noindent A possible candidate is the saddle-type solution built in \cite{CT1}, at least in even dimensions $N.$  

\end{Oquestion} 

\medskip

\begin{Oquestion} 

\noindent {\it Assume $ N \ge 4$. Prove (or disprove) that every monotone solution has finite energy. }

\medskip

\noindent As pointed out in Remark \ref{Rk-monot-finite-energy}, this is true for $ N\le 3$.

\end{Oquestion}

\medskip

\begin{Oquestion} 

\noindent {\it Assume $ N \ge 2$ and let $u$ be a non constant finite-energy solution of the Allen-Cahn equation \eqref{eqAC}. Prove (or disprove) that

\[
u \, \textit{ is 1D} \quad \Longleftrightarrow \quad L_u = \ell_{min}. 
\]

}

\medskip

\noindent Note that this is true when $u$ is also a local minimizer, as proved in Theorem \ref{1D-minloc}. It is also true for stable solutions when $ N =2,3$ (cf. Theorem \ref{Teo-stabili-folklore}).

\end{Oquestion}

\medskip

\begin{Oquestion} 

\noindent {\it Assume $ N \ge 2$ and let $u$ be a non constant {\sl monotone} solution of the Allen-Cahn equation \eqref{eqAC}. 
Prove (or disprove) that
}
\[
u \, \textit{ is 1D} \quad \Longleftrightarrow \quad L_u = \ell_{min}. 
\]

\medskip

\noindent Note that this is true when $u$ is also a local minimizer, as proved in Theorem \ref{1D-minloc}. 

\end{Oquestion}

\medskip

\begin{Oquestion}
\noindent {\it Assume $ N=2$. Do there exist entire solutions of the Allen-Cahn equation \eqref{eqAC} having a Pythagorean conical limit interface $V$ ? }
\end{Oquestion}

\medskip

\begin{Oquestion} 
\noindent {\it Assume $ N\ge 5$. Is the second positive energy level ${\ell}_{2}$ achieved by a member of ${\cal F}_p$ ? In the affirmative, compute the exact value. }

\bigskip

\noindent A good candidate is the saddle-type solution built in \cite{CT1}, at least in even dimensions $N.$  

\end{Oquestion}

\bigskip

\vspace{0.5cm}

\noindent \textbf{Acknowledgements: } This work was started at the beginning of December 2014, when the author was visiting the Scuola Normale Superiore di Pisa, Italy. The author is grateful to Scuola Normale Superiore di Pisa for the kind invitation and the warm hospitality. The author is indebted to Carlo Mantegazza 
for many interesting and stimulating discussions on the subject of the present work as well as on geometric measure theory. He also thanks Luigi Ambrosio, Andrea Malchiodi and G\'erard Tronel for reading earlier versions of this work and for their constant support.

The author is partially supported by the ERC grant EPSILON ({\it Elliptic Pde's and Symmetry of Interfaces and Layers for Odd Nonlinearities}) and by the ERC grant COMPAT ({\it Complex Patterns for Strongly Interacting Dynamical Systems}).

\end{document}